\documentclass[12pt,leqno]{amsart}
\usepackage{amssymb}
\usepackage{euscript}
 \input xy
   \xyoption{all}
\newcommand{\lon}{\longrightarrow}
\newcommand{\rar}{\rightarrow}
\newcommand{\ot}{\otimes}
\newcommand{\tl}{\tilde}
\newcommand{\sip}{\smallskip}
\newcommand{\bip}{\bigskip}
\newcommand{\Z}{{\mathbb Z}}
\newcommand{\fg}{{\mathfrak g}}
\newcommand{\fG}{{\mathfrak G}}
\newcommand{\tln}{{\tilde{1}}}
\newcommand{\tlo}{{\tilde{0}}}
\newcommand{\p}{\partial}
\newcommand{\cT}{{\mathcal T}}
\newcommand{\cM}{{\mathcal M}}
\newcommand{\Id}{{\mathrm I\mathrm d}}
\newcommand{\f}{\mathcal O}
\newcommand{\al}{\alpha}
 \newcommand{\be}{\beta}

\newcommand{\bS}{{\Sigma}}
\newcommand{\fA}{{\mathfrak A}}
\newcommand{\sG}{{\mathsf G}}
\newcommand{\one}{{1\hspace{-1.3mm}1}}

 \newcommand{\Beq}{\begin{equation}}
 \newcommand{\Eeq}{\end{equation}}
 \newcommand{\Beqr}{\begin{eqnarray}}
 \newcommand{\Eeqr}{\end{eqnarray}}
 \newcommand{\Beqrn}{\begin{eqnarray*}}
 \newcommand{\Eeqrn}{\end{eqnarray*}}
 \newcommand{\Ba}{\begin{array}}
 \newcommand{\Ea}{\end{array}}
 \newcommand{\Bi}{\begin{itemize}}
 \newcommand{\Ei}{\end{itemize}}
 \newcommand{\Bc}{\begin{center}}
 \newcommand{\Ec}{\end{center}}

 \newcommand{\cE}{{\mathcal E}}
 \newcommand{\cF}{{\mathcal F}}

 \newcommand{\cW}{{\mathcal W}}
\newcommand{\R}{{\mathbb R}}

\swapnumbers
\newtheorem{theorem}{Theorem}[subsection]

\newtheorem{lemma}[theorem]{Lemma}
\newtheorem{proposition}[theorem]{Proposition}

\newtheorem{prop-def}[theorem]{Proposition-definition}

\newtheorem*{theorem-A}{Theorem~A}
\newtheorem*{theorem-B}{Theorem~B}
\newtheorem*{theorem-C}{Theorem~C}
\newtheorem*{corollary-D}{Corollary~D}
\newtheorem*{theorem-E}{Theorem~E}

\theoremstyle{definition}
\newtheorem{example}[theorem]{Example}
\newtheorem{remark}[theorem]{Remark}
\newtheorem{definition}[theorem]{Definition}

\newtheorem{CY}[theorem]{Calabi-Yau manifolds}
\newtheorem{MES}[theorem]{Master equation on supermanifolds and complexes}
\newtheorem{fact}[theorem]{Fact}

\def\UPQ#1#2{{}_{\raisebox{-.11em}{\scriptsize \it#1}}\hskip .05em%
              \fA^\circlearrowright_{#2}}
\def\Cof#1#2{{}_{#1}U_{#2}}
\def\tCof#1#2{{}_{#1}{\widetilde U}_{#2}}
\def\tM#1#2{{}_{#1}{\widetilde M}_{#2}}

\def\MM#1#2{{}_{#1}{M}_{#2}}
\def\NN#1#2{{}_{#1}{N}_{#2}}
\def\OO{{}_{r}{O}_{s}}
\def\KK{{}_{r}{K}_{s}}
\def\Srs{{}_{\raisebox{-.11em}{\scriptsize \it r}}\hskip .05emS_s^1}
\def\Cyc{{\it Cyc}}

\def\brackett{{
\unitlength=.05em
\begin{picture}(24.00,20.00)(-2.00,2.00)
\bezier{40}(10.00,10.00)(15.00,5.00)(20.00,0.00)
\bezier{40}(10.00,10.00)(5.00,5.00)(0.00,0.00)
\put(10.00,20.00){\line(0,-1){10.00}}
\put(10.00,9.00){\makebox(0,0){$\bullet$}}
\end{picture}}
}

\def\krouzek#1{
\setlength{\unitlength}{#1em}
\bezier{100}(-2,0)(-2,.5)(-1.73,1)
\bezier{100}(-1.73,1)(-1.48,1.48)(-1,1.73)
\bezier{100}(-1,1.73)(-.5,2)(0,2)
\bezier{100}(2,0)(2,.5)(1.73,1)
\bezier{100}(1.73,1)(1.48,1.48)(1,1.73)
\bezier{100}(1,1.73)(.5,2)(0,2)
\bezier{100}(2,0)(2,-.5)(1.73,-1)
\bezier{100}(1.73,-1)(1.48,-1.48)(1,-1.73)
\bezier{100}(1,-1.73)(.5,-2)(0,-2)
\bezier{100}(-2,0)(-2,-.5)(-1.73,-1)
\bezier{100}(-1.73,-1)(-1.48,-1.48)(-1,-1.73)
\bezier{100}(-1,-1.73)(-.5,-2)(0,-2)
}

\def\dotkrouzek#1{
\setlength{\unitlength}{#1em}
\bezier{100}(-2,0)(-2,.5)(-1.73,1)
\bezier{5}(-1.73,1)(-1.48,1.48)(-1,1.73)
\bezier{100}(-1,1.73)(-.5,2)(0,2)
\bezier{5}(2,0)(2,.5)(1.73,1)
\bezier{100}(1.73,1)(1.48,1.48)(1,1.73)
\bezier{5}(1,1.73)(.5,2)(0,2)
\bezier{100}(2,0)(2,-.5)(1.73,-1)
\bezier{5}(1.73,-1)(1.48,-1.48)(1,-1.73)
\bezier{100}(1,-1.73)(.5,-2)(0,-2)
\bezier{100}(-2,0)(-2,-.5)(-1.73,-1)
\bezier{100}(-1.73,-1)(-1.48,-1.48)(-1,-1.73)
\bezier{5}(-1,-1.73)(-.5,-2)(0,-2)
}

\def\hkrouzek#1{
\setlength{\unitlength}{#1em}
\bezier{100}(-2,0)(-2,.5)(-1.73,1)
\bezier{100}(-1.73,1)(-1.48,1.48)(-1,1.73)
\bezier{100}(-1,1.73)(-.5,2)(0,2)
\bezier{100}(2,0)(2,.5)(1.73,1)
\bezier{100}(1.73,1)(1.48,1.48)(1,1.73)
\bezier{100}(1,1.73)(.5,2)(0,2)
}

\def\dkrouzek#1{
\setlength{\unitlength}{#1em}
\bezier{100}(2,0)(2,-.5)(1.73,-1)
\bezier{100}(1.73,-1)(1.48,-1.48)(1,-1.73)
\bezier{100}(1,-1.73)(.5,-2)(0,-2)
\bezier{100}(-2,0)(-2,-.5)(-1.73,-1)
\bezier{100}(-1.73,-1)(-1.48,-1.48)(-1,-1.73)
\bezier{100}(-1,-1.73)(-.5,-2)(0,-2)
}

\def\hgw{{\mbox {\large $\circlearrowright$}}}

\def\catPROP{{\tt PROP}} \def\calP{{\mathcal P}}
\def\catWPROP{\mbox{${\tt PROP}^\circlearrowright$}}
\def\wc{\circlearrowright} \def\Ker{\mbox{\it Ker}}
 \def\Coim{\mbox{\it Coim}}
\def\susp{\uparrow\!} \def\desusp{\downarrow\!}
\def\bcalP{{\overline {\calP}}}
\def\freeOP (#1){\Gamma_{\hskip -.2em\it op} \langle #1 \rangle}
\def\freeMproperad (#1){{\rm MF}^\wc \langle #1 \rangle}
\def\freeMproperadname {{\rm MF}^\wc}
\def\susps{\uparrow\!\!} \def\desusps{\downarrow\!\!}
\def\Tr{{\it Tr}}\def\calC{{\mathcal C}}\def\ext{\mbox{\large$\land$}}
\def\Aut{{\it Aut\/}}
\def\bsfC{{\overline{\sf C}}}\def\rmF(#1){{\rm F\langle #1 \rangle}}
\def\WrmF(#1){{\rm F^\circlearrowright}\langle #1 \rangle}\def\fW{\WrmF}
\def\WrmFname{{\rm F^\circlearrowright}}
\def\rmFname{{\rm F}}
\def\cW(#1){{\rm C^\circlearrowright}\langle #1 \rangle}
\def\cWname{{\rm C^\circlearrowright}}
\def\wB{{B^\circlearrowright}}
\def\wO{{\Omega^\circlearrowright}}
\def\whO{{\widehat{\Omega}^\circlearrowright}}
\def\bcW(#1){{\rm \overline{C}^\circlearrowright}\langle #1 \rangle}
\def\bWrmF(#1){{\rm \overline{F}^\circlearrowright}\langle#1\rangle}
\def\DGrw#1#2{{{\tt G}^\circlearrowright}(#1,#2)}
\def\bsfP{{\overline{\sfP}}}\def\sfC{{\sf C}}
\def\Lie{{\sf Lie}} \def\Com{{\sf Com}}
\def\Ass{{\sf Ass}} \def\frakk{{\mathfrak k}}
\def\sgn{{\rm sgn}} \def\ww{{\mathbf w}}
\def\leg{{\it Leg}}\def\bfk{{\mathbf k}}
\def\Vert{{\it Vert}}\def\vert{{\it vert}}
\def\Flag{{\it Flag}}\def\Edg{{\it Edg}}
\def\rada#1#2{#1,\ldots,#2}
\def\Hom{{\it Hom}}
\def\DGrwc#1#2{{{\tt G}_c^\circlearrowright}(#1,#2)}
\def\bDGrwc#1#2{{\overline{\tt G}_c^\circlearrowright}(#1,#2)}
\def\DGrwcnic{{{\tt G}_c^\circlearrowright}}

\def\id{{\mbox{1 \hskip -8pt 1}}}

\def\PROP{{\small PROP}}

\def\bfk{{\mathbf k}}
\def\ot{\otimes}

\def\sfP{{\sf P}} \def\pa{\partial}

\def\Magogh#1#2#3{{\tt G}_c^\wc(#1;#2,#3)}
\def\bicol#1#2{\{#1(m,n)\}_{m,n \geq #2}}

\def\End{{\mathcal E} \hskip -.1em{\it nd\/}}
\def\Lin{{\it Lin\/}}
\def\otexp#1#2{#1^{\otimes #2}}
\def\sigmabimod{{\mbox{$\Sigma$-{\tt bimod}}}}

\newcommand{\Bij}{{\it Bij}}
\def\DGr#1#2{{\tt G}^\uparrow(#1,#2)}
\def\edge{{\it edge}}\def\DGrc#1#2{{\tt G}^\uparrow_c(#1,#2)}
\def\bfreeWPROP (#1){{\overline{\Gamma}^\circlearrowright\langle #1 \rangle}}
\def\In{{\it in}} \def\Out{{\it out}}
\def\freePROP (#1){{\Gamma\langle #1 \rangle}}
\def\freeWPROP (#1){{\Gamma^\circlearrowright\langle #1 \rangle}}
\def\cfreePROP (#1){{\Gamma_c\langle #1 \rangle}}
\def\cfreeWPROP (#1){{\Gamma_c^\circlearrowright \langle #1 \rangle}}
\def\bcfreeWPROP (#1){{\overline{\Gamma}_c^\circlearrowright\langle#1\rangle}}

\def\colim#1{\mathop{{\rm colim}}%
             \limits_{\rule{0em}{1em}\mbox{\scriptsize $#1$}}}

\def\cases#1#2#3#4{
                  \left\{
                         \begin{array}{ll}
                           #1,\ &\mbox{#2}
                           \\
                           #3,\ &\mbox{#4}
                          \end{array}
                   \right.
}

\begin{document}
\pagestyle{myheadings}
\bibliographystyle{plain}
\baselineskip18pt plus 1pt minus 1pt
\parskip3pt plus 1pt minus .5pt

\title{Wheeled PROPs, graph complexes and the master equation}

\author{ M.~Markl, S.~Merkulov and S.~Shadrin}
\thanks{The first author was supported by the grant GA \v CR 201/05/2117 and by
   the Academy of Sciences of the Czech Republic,
   Institutional Research Plan No.~AV0Z10190503. The second author was
supported by the G\"oran Gustafsson foundation.
The third author was supported by the grants RFBR-04-02-17227-a,
RFBR-05-01-02806-CNRS-a, NSh-1972.2003.1, MK-5396.2006.1,
NWO-RFBR-047.011.2004.026 (RFBR-05-02-89000-NWO-a),
by the G{\"o}ran Gustafsson foundation, and by
Pierre Deligne's fund based on his 2004 Balzan prize in mathematics.
}

\address{M.~Markl: Math. Inst. of the Academy, {\v Z}itn{\'a} 25,
         115 67 Prague 1, The Czech Republic}
\email{markl@math.cas.cz}
\address{S.~Merkulov: Department of Mathematics, Stockholm University, 10691 Stockholm, Sweden}
\email{sm@math.su.se}
\address{S.~Shadrin: Department of Mathematics, Stockholm University, 10691 Stockholm, Sweden
and Institute of System Research, Nakhimovskii prospekt~36-1, 117218 Moscow, Russia}
\email{shadrin@math.su.se, shadrin@mccme.ru}

\begin{abstract}
We introduce and study wheeled PROPs, an extension of the theory of PROPs which can treat
traces and, in particular, solutions to the master equations which involve divergence
operators. We construct a dg free wheeled PROP whose representations are in one-to-one
correspondence with formal germs of SP-manifolds, key geometric objects in the theory
of Batalin-Vilkovisky quantization. We also construct minimal wheeled resolutions of
classical operads $\Com$ and $\Ass$ as rather non-obvious extensions of
$\Com_\infty$ and $\Ass_\infty$,
involving, e.g., a mysterious mixture of associahedra with cyclohedra. Finally, we
apply the above results to a computation of cohomology of a directed
version of Kontsevich's complex of ribbon graphs.
\end{abstract}

\maketitle

\tableofcontents

\section{Introduction}

\noindent
Theory of operads and, more generally, of PROPs undergoes a rapid
development in recent years
\cite{handbook, markl-shnider-stasheff:book}. Originated in algebraic
topology and homotopy theory, it has provided us with new
powerful means to attack problems in many other areas of mathematics
including deformation theory and differential geometry. Rather
surprisingly, such classical systems of non-linear differential
equations as Poisson equations, Nijenhuis integrability equations
for an almost complex structure, or even Hochschild equations for
star products on germs of smooth functions,
--- they all can be understood as representations of certain dg
PROPs~\cite{me1,me2,me3}. The common feature of all these equations
is that
 they admit an interpretation as Maurer-Cartan equations in certain
dg Lie algebras. However, in theoretical physics one often
encounters not the Maurer-Cartan equations but the {\em master
equations}\, which involve a divergence operator. Such equations can
{\em not}\, be understood in terms of ordinary PROPs and operads
as the latter have no room to encode such a basic operation in
geometry and physics as {\em trace.}

In this paper we introduce and study an extension of the theory of
PROPs which can treat traces and hence can be used to describe
solutions to master equations. This new theory turns out to be in a
sense simpler than the original one. We call it the theory of {\em
wheeled}\, PROPs and operads.
The motivation for the terminology is that
free objects in this
new category are based on directed graphs with may contain wheels,
that is, directed {\em closed}\, paths of oriented
 edges. Such graph complexes as well as some methods of computing
their cohomology have been studied in \cite{me3}.

There is a canonical forgetful functor, $\catPROP^\circlearrowright
\rar \catPROP$, from the category of wheeled PROPs to the category
of PROPs
 which admits a left adjoint,
$(-)^\circlearrowright$:~$\catPROP\rightarrow
\catPROP^\circlearrowright$. This functor is neither exact nor full, and its study is another main purpose of
 our work.

This means, in particular, that a resolution,
${\mathsf A}_\infty$ of, say,
 an operad $\mathsf A$ within the category of ordinary operads may {not},
in general, produce via the wheeled completion $\circlearrowright$
 a resolution, $({\mathsf A}^\circlearrowright)_\infty$,
of its wheeled completion, $\mathsf A^\circlearrowright$, in the
category of wheeled PROPs, i.e., in general
$$
({\mathsf A}^\circlearrowright)_\infty \neq \ ({\mathsf
A}_\infty)^\circlearrowright.
$$
This phenomenon is studied in detail for classical
operads, ${\Ass}$ and ${\Com}$, for associative and commutative algebras:

$\bullet$
first, we compute cohomology groups
\[
\bigoplus_{n\leq 0}H^n((\mbox{\Ass}_\infty)^\circlearrowright, \p)\ \ \ \ \mbox{and}\ \ \ \ \
\bigoplus_{n\leq 0}H^n((\mbox{\Com}_\infty)^\circlearrowright, \p)
\]
of the wheeled completions of their well-known
minimal resolutions, $(\mbox{\sf Ass}_\infty, \p)$ and
$(\mbox{\sf Com}_\infty, \p)$, and show that these groups are
non-zero for {\em all}\, $n\leq 0$
(with
$H^0((\mbox{\Ass}_\infty)^\circlearrowright, \p)={\Ass}^\circlearrowright$ and
$H^0((\mbox{\Com}_\infty)^\circlearrowright,
\p)={\Com}^\circlearrowright$
as expected);

$\bullet$
second, by adding new generators we construct their {\em
wheeled}\, minimal resolutions,
$$
{\Ass}^\circlearrowright_\infty :=({\Ass}^\circlearrowright)_\infty
\ \ \ \mbox{and}\ \ \
{\Com}^\circlearrowright_\infty:=
(\Com^\circlearrowright)_\infty,
$$
 as rather non-trivial extensions of
$(\mbox{\sf Ass}_\infty)^\circlearrowright$ and, respectively,
 $(\mbox{\sf Com}_\infty)^\circlearrowright$, involving, e.g., a mysterious mixture of associahedra
with cyclohedra.

More precisely, we prove the following two theorems:

\begin{theorem-A}

The minimal wheeled resolution of the
operad  of associative algebras, $\Ass$,
is the free wheeled dg PROP, $(\Ass_\infty^\circlearrowright, \pa)$,
 generated by

\begin{itemize}
\item[(i)] planar (1,n)-corollas in degree $2-n$,
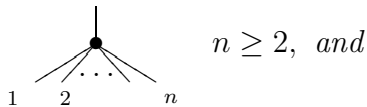
\[
\begin{xy}
 <0mm,0mm>*{\bullet};<0mm,0mm>*{}**@{},
 <0mm,0mm>*{};<-8mm,-5mm>*{}**@{-},
 <0mm,0mm>*{};<-4.5mm,-5mm>*{}**@{-},
 <0mm,0mm>*{};<0mm,-4mm>*{\ldots}**@{},
 <0mm,0mm>*{};<4.5mm,-5mm>*{}**@{-},
 <0mm,0mm>*{};<8mm,-5mm>*{}**@{-},
   <0mm,0mm>*{};<-11mm,-7.9mm>*{^{1}}**@{},
   <0mm,0mm>*{};<-4mm,-7.9mm>*{^{2}}**@{},
   <0mm,0mm>*{};<10.0mm,-7.9mm>*{^{n}}**@{},
 <0mm,0mm>*{};<0mm,5mm>*{}**@{-},
 \end{xy}\ \ \ n\geq 2,\ \mbox {and}
\]
\item[(ii)] planar $(0,m+n)$-corollas in degree $-m-n$
\[
\begin{xy}
 <0mm,-0.5mm>*{\blacktriangledown};<0mm,0mm>*{}**@{},
 <0mm,0mm>*{};<-16mm,-5mm>*{}**@{-},
 <0mm,0mm>*{};<-12mm,-5mm>*{}**@{-},
 <0mm,0mm>*{};<-3.5mm,-5mm>*{}**@{-},
 <0mm,0mm>*{};<-6mm,-5mm>*{...}**@{},
   <0mm,0mm>*{};<-21mm,-7.9mm>*{^{1}}**@{},
   <0mm,0mm>*{};<-14mm,-7.9mm>*{^{2}}**@{},
   <0mm,0mm>*{};<-4mm,-7.9mm>*{^{m}}**@{},
 <0mm,0mm>*{};<16mm,-5mm>*{}**@{-},
 <0mm,0mm>*{};<12mm,-5mm>*{}**@{-},
 <0mm,0mm>*{};<3.5mm,-5mm>*{}**@{-},
 <0mm,0mm>*{};<6.6mm,-5mm>*{...}**@{},
   <0mm,0mm>*{};<22mm,-7.9mm>*{^{m\hspace{-0.5mm}+\hspace{-0.5mm}n}}**@{},
   <0mm,0mm>*{};<6mm,-7.9mm>*{^{m\hspace{-0.5mm}+\hspace{-0.5mm}1}}**@{},
 \end{xy}\ \ \ m,n\geq 1,
\]
\end{itemize}
having the cyclic skew-symmetry
\[
\hskip 3em
\begin{xy}
 <0mm,-0.5mm>*{\blacktriangledown};<0mm,0mm>*{}**@{},
 <0mm,0mm>*{};<-16mm,-5mm>*{}**@{-},
 <0mm,0mm>*{};<-12mm,-5mm>*{}**@{-},
 <0mm,0mm>*{};<-3.5mm,-5mm>*{}**@{-},
 <0mm,0mm>*{};<-6mm,-5mm>*{...}**@{},
   <0mm,0mm>*{};<-19mm,-7.9mm>*{^{1}}**@{},
   <0mm,0mm>*{};<-13mm,-7.9mm>*{^{2}}**@{},
   <0mm,0mm>*{};<-4mm,-7.9mm>*{^{m}}**@{},
 <0mm,0mm>*{};<16mm,-5mm>*{}**@{-},
 <0mm,0mm>*{};<12mm,-5mm>*{}**@{-},
 <0mm,0mm>*{};<3.5mm,-5mm>*{}**@{-},
 <0mm,0mm>*{};<6.6mm,-5mm>*{...}**@{},
   <0mm,0mm>*{};<19mm,-7.9mm>*{^{m\hspace{-0.5mm}+\hspace{-0.5mm}n}}**@{},
   <0mm,0mm>*{};<5mm,-7.9mm>*{^{m\hspace{-0.5mm}+\hspace{-0.5mm}1}}**@{},
 \end{xy}
\hskip -.5em
= (-1)^{\sgn(\zeta)}\hspace{-13mm}
\begin{xy}
 <0mm,-0.5mm>*{\blacktriangledown};<0mm,0mm>*{}**@{},
 <0mm,0mm>*{};<-16mm,-5mm>*{}**@{-},
 <0mm,0mm>*{};<-12mm,-5mm>*{}**@{-},
 <0mm,0mm>*{};<-3.5mm,-5mm>*{}**@{-},
 <0mm,0mm>*{};<-6mm,-5mm>*{...}**@{},
   <0mm,0mm>*{};<-21mm,-7.9mm>*{^{\zeta(1)}}**@{},
   <0mm,0mm>*{};<-14mm,-7.9mm>*{^{\zeta(2)}}**@{},
   <0mm,0mm>*{};<-4mm,-7.9mm>*{^{\zeta(m)}}**@{},
 <0mm,0mm>*{};<16mm,-5mm>*{}**@{-},
 <0mm,0mm>*{};<12mm,-5mm>*{}**@{-},
 <0mm,0mm>*{};<3.5mm,-5mm>*{}**@{-},
 <0mm,0mm>*{};<6.6mm,-5mm>*{...}**@{},
   <0mm,0mm>*{};<20mm,-7.9mm>*{^{m\hspace{-0.5mm}+\hspace{-0.5mm}n}}**@{},
   <0mm,0mm>*{};<6mm,-7.9mm>*{^{m\hspace{-0.5mm}+\hspace{-0.5mm}1}}**@{},
 \end{xy}
\hskip -.5em
 =(-1)^{\sgn(\xi)}\hspace{-11mm}
\begin{xy}
 <0mm,-0.5mm>*{\blacktriangledown};<0mm,0mm>*{}**@{},
 <0mm,0mm>*{};<-16mm,-5mm>*{}**@{-},
 <0mm,0mm>*{};<-12mm,-5mm>*{}**@{-},
 <0mm,0mm>*{};<-3.5mm,-5mm>*{}**@{-},
 <0mm,0mm>*{};<-6mm,-5mm>*{...}**@{},
   <0mm,0mm>*{};<-21mm,-7.9mm>*{^{1}}**@{},
   <0mm,0mm>*{};<-14mm,-7.9mm>*{^{2}}**@{},
   <0mm,0mm>*{};<-4mm,-7.9mm>*{^{m}}**@{},
 <0mm,0mm>*{};<16mm,-5mm>*{}**@{-},
 <0mm,0mm>*{};<12mm,-5mm>*{}**@{-},
 <0mm,0mm>*{};<3.5mm,-5mm>*{}**@{-},
 <0mm,0mm>*{};<6.6mm,-5mm>*{...}**@{},
   <0mm,0mm>*{};<20mm,-7.9mm>*{^{\xi(m\hspace{-0.5mm}+\hspace{-0.5mm}n)}}**@{},
   <0mm,0mm>*{};<6mm,-7.9mm>*{^{\xi(m\hspace{-0.5mm}+\hspace{-0.5mm}1)}}**@{},
 \end{xy}
\]
with respect to the cyclic permutations $\zeta=(12\ldots m)$
and $\xi=((m+1)(m+2)\ldots(m+n))$.
The differential is given on generators as
\begin{eqnarray}
\label{ass1}
\partial
\begin{xy}
 <0mm,0mm>*{\bullet};<0mm,0mm>*{}**@{},
 <0mm,0mm>*{};<-8mm,-5mm>*{}**@{-},
 <0mm,0mm>*{};<-4.5mm,-5mm>*{}**@{-},
 <0mm,0mm>*{};<0mm,-4mm>*{\ldots}**@{},
 <0mm,0mm>*{};<4.5mm,-5mm>*{}**@{-},
 <0mm,0mm>*{};<8mm,-5mm>*{}**@{-},
   <0mm,0mm>*{};<-9mm,-7.9mm>*{^1}**@{},
   <0mm,0mm>*{};<-5mm,-7.9mm>*{^2}**@{},
   <0mm,0mm>*{};<10.0mm,-7.9mm>*{^n}**@{},
 <0mm,0mm>*{};<0mm,5mm>*{}**@{-},
 \end{xy}
&=&\sum_{k=0}^{n-2}\sum_{l=2}^{n-k}
(-1)^{k+l(n-k-l)+1}
\begin{xy}
<0mm,0mm>*{\bullet},
<0mm,5mm>*{}**@{-},
<4mm,-7mm>*{^{1\  \dots \  k\qquad\ \ k+l+1\ \ \dots\ \  n }},
<-14mm,-5mm>*{}**@{-},
<-6mm,-5mm>*{}**@{-},
<20mm,-5mm>*{}**@{-},
<8mm,-5mm>*{}**@{-},
<0mm,-5mm>*{}**@{-},
<0mm,-5mm>*{\bullet};
<-5mm,-10mm>*{}**@{-},
<-2mm,-10mm>*{}**@{-},
<2mm,-10mm>*{}**@{-},
<5mm,-10mm>*{}**@{-},
<0mm,-12mm>*{_{k+1\ \, \dots\ \, k+l}},
\end{xy},
\\
\label{ass2}
\hskip 2em \partial \hskip -1em
\begin{xy}
 <0mm,-0.5mm>*{\blacktriangledown};<0mm,0mm>*{}**@{},
 <0mm,0mm>*{};<-16mm,-5mm>*{}**@{-},
 <0mm,0mm>*{};<-12mm,-5mm>*{}**@{-},
 <0mm,0mm>*{};<-3.5mm,-5mm>*{}**@{-},
 <0mm,0mm>*{};<-6mm,-5mm>*{...}**@{},
   <0mm,0mm>*{};<-19mm,-7.9mm>*{^{1}}**@{},
   <0mm,0mm>*{};<-13mm,-7.9mm>*{^{2}}**@{},
   <0mm,0mm>*{};<-4mm,-7.9mm>*{^{m}}**@{},
 <0mm,0mm>*{};<16mm,-5mm>*{}**@{-},
 <0mm,0mm>*{};<12mm,-5mm>*{}**@{-},
 <0mm,0mm>*{};<3.5mm,-5mm>*{}**@{-},
 <0mm,0mm>*{};<6.6mm,-5mm>*{...}**@{},
   <0mm,0mm>*{};<19mm,-7.9mm>*{^{m\hspace{-0.5mm}+\hspace{-0.5mm}n}}**@{},
   <0mm,0mm>*{};<5mm,-7.9mm>*{^{m\hspace{-0.5mm}+\hspace{-0.5mm}1}}**@{},
 \end{xy}
\hskip -1.2em
& = &
\sum_{i=0}^{m-1}\left((-1)^{m+1}\zeta\right)^i\sum_{j=1}^{n-1}
\left(
(-1)^{n+1}\xi\right)^j \hskip -.2em \left(\rule{0em}{2.3em}\right.
\begin{xy}
 <0mm,0mm>*{\bullet};<0mm,0mm>*{}**@{},
 <0mm,0mm>*{};<-11mm,-5mm>*{}**@{-},
 <0mm,0mm>*{};<-8.5mm,-5mm>*{}**@{-},
 <0mm,0mm>*{};<-5mm,-5mm>*{...}**@{},
 <0mm,0mm>*{};<-2.5mm,-5mm>*{}**@{-},
   <0mm,0mm>*{};<-13mm,-7.9mm>*{^1}**@{},
   <0mm,0mm>*{};<-9mm,-7.9mm>*{^2}**@{},
   <0mm,0mm>*{};<-3mm,-7.9mm>*{^m}**@{},
<0mm,0mm>*{};<11mm,-5mm>*{}**@{-},
 <0mm,0mm>*{};<8.5mm,-5mm>*{}**@{-},
 <0mm,0mm>*{};<5.5mm,-5mm>*{...}**@{},
 <0mm,0mm>*{};<2.5mm,-5mm>*{}**@{-},
   <0mm,0mm>*{};<14mm,-7.9mm>*{^{m\hspace{-0.5mm}+\hspace{-0.5mm}n}}**@{},
   <0mm,0mm>*{};<3mm,-7.9mm>*{^{m\hspace{-0.5mm}+\hspace{-0.5mm}1}}**@{},
 <0mm,-10mm>*{};<0mm,5mm>*{}**@{-},
(0,5)*{}
   \ar@{->}@(ur,dr) (0,-10)*{}
 \end{xy}
 \\
 &&
\hskip -10em
+ \sum_{k=2}^{m}
(-1)^{k(m+n)} \hskip -1em
\begin{xy}
 <0mm,-0.5mm>*{\blacktriangledown};<0mm,0mm>*{}**@{},
 <0mm,0mm>*{};<-16mm,-5mm>*{}**@{-},
 <0mm,0mm>*{};<-11mm,-5mm>*{}**@{-},
 <0mm,0mm>*{};<-3.5mm,-5mm>*{}**@{-},
 <0mm,0mm>*{};<-6mm,-5mm>*{...}**@{},
   <0mm,0mm>*{};<-10mm,-7.9mm>*{^{k\hspace{-0.5mm}+\hspace{-0.5mm}1}}**@{},
   <0mm,0mm>*{};<-4mm,-7.9mm>*{^{m}}**@{},
 <0mm,0mm>*{};<16mm,-5mm>*{}**@{-},
 <0mm,0mm>*{};<12mm,-5mm>*{}**@{-},
 <0mm,0mm>*{};<3.5mm,-5mm>*{}**@{-},
 <0mm,0mm>*{};<6.6mm,-5mm>*{...}**@{},
   <0mm,0mm>*{};<19mm,-7.9mm>*{^{m\hspace{-0.5mm}+\hspace{-0.5mm}n}}**@{},
   <0mm,0mm>*{};<5mm,-7.9mm>*{^{m\hspace{-0.5mm}+\hspace{-0.5mm}1}}**@{},
 <-16mm,-5.5mm>*{\bullet};<0mm,0mm>*{}**@{},
 <-16mm,-5.5mm>*{};<-20mm,-11mm>*{}**@{-},
 <-16mm,-5.5mm>*{};<-12mm,-11mm>*{}**@{-},
 <-16mm,-5.5mm>*{};<-18mm,-11mm>*{}**@{-},
 <-16mm,-5.5mm>*{};<-14mm,-11mm>*{}**@{-},
 <-16mm,-13mm>*{_{1\ \, \dots\ \ k}},
 \end{xy}
\nonumber
+ \sum_{k=2}^{n-2}
(-1)^{m+k+nk+1}
\begin{xy}
<0mm,-0.5mm>*{\blacktriangledown};<0mm,0mm>*{}**@{},
 <0mm,0mm>*{};<-16mm,-5mm>*{}**@{-},
 <0mm,0mm>*{};<-12mm,-5mm>*{}**@{-},
 <0mm,0mm>*{};<-3.5mm,-5mm>*{}**@{-},
 <0mm,0mm>*{};<-6mm,-5mm>*{...}**@{},
   <0mm,0mm>*{};<-18mm,-7.9mm>*{^{1}}**@{},
   <0mm,0mm>*{};<-13mm,-7.9mm>*{^{2}}**@{},
   <0mm,0mm>*{};<-4mm,-7.9mm>*{^{m}}**@{},
 <0mm,0mm>*{};<17mm,-5mm>*{}**@{-},
 <0mm,0mm>*{};<7mm,-5mm>*{}**@{-},
 <0mm,0mm>*{};<3.5mm,-5mm>*{}**@{-},
 <0mm,0mm>*{};<11.6mm,-5mm>*{...}**@{},
   <0mm,0mm>*{};<22mm,-7.9mm>*{^{m\hspace{-0.5mm}+\hspace{-0.5mm}n}}**@{},
   <0mm,0mm>*{};<11mm,-7.9mm>*{^{m\hspace{-0.5mm}+\hspace{-0.5mm}k+\hspace{-0.5mm}1}}**@{},
 <3.5mm,-5.5mm>*{\bullet};<0mm,0mm>*{}**@{},
 <3.5mm,-5.5mm>*{};<-0.5mm,-11mm>*{}**@{-},
 <3.5mm,-5.5mm>*{};<1.5mm,-11mm>*{}**@{-},
 <3.5mm,-5.5mm>*{};<5.5mm,-11mm>*{}**@{-},
 <3.5mm,-5.5mm>*{};<7.5mm,-11mm>*{}**@{-},
 <3.5mm,-13mm>*{_{m+1\, \dots\ m+k}},
 \end{xy}\left. \rule{0em}{2.3em}\right).
\nonumber
\end{eqnarray}
\end{theorem-A}

In Theorem~B below, $\Sigma(k,n)$ denotes the set of all
$(k,n-k)$-unshuffles,
\[
\Sigma(k,n):=\{\tau\in\bS_n\,|\,\tau(1)<\dots<\tau(k),\,
\tau(k+1)<\dots<\tau(n)\}.
\]

\begin{theorem-B}
The minimal wheeled resolution of the operad
of commutative associative
algebras, $\Com$,
is the free wheeled dg PROP, $(\Com_\infty^\circlearrowright, \pa)$,
generated by
\begin{itemize}
\item[(i)]
planar $(1,n)$-corollas in degree $2-n$,
\[
\begin{xy}
 <0mm,0mm>*{\bullet};<0mm,0mm>*{}**@{},
 <0mm,0mm>*{};<-8mm,-5mm>*{}**@{-},
 <0mm,0mm>*{};<-4.5mm,-5mm>*{}**@{-},
 <0mm,0mm>*{};<0mm,-4mm>*{\ldots}**@{},
 <0mm,0mm>*{};<4.5mm,-5mm>*{}**@{-},
 <0mm,0mm>*{};<8mm,-5mm>*{}**@{-},
   <0mm,0mm>*{};<-11mm,-7.9mm>*{^{1}}**@{},
   <0mm,0mm>*{};<-4mm,-7.9mm>*{^{2}}**@{},
     <0mm,0mm>*{};<10.0mm,-7.9mm>*{^{n}}**@{},
 <0mm,0mm>*{};<0mm,5mm>*{}**@{-},
 \end{xy}\ \ \ n\geq 2,
\]
modulo the shuffle relations:
\begin{equation}
\sum_{\tau\in\Sigma(k,n)}
(-1)^{\sgn(\tau)}
\begin{xy}
 <0mm,0mm>*{\bullet};<0mm,0mm>*{}**@{},
 <0mm,0mm>*{};<-8mm,-5mm>*{}**@{-},
 <0mm,0mm>*{};<-4.5mm,-5mm>*{}**@{-},
 <0mm,0mm>*{};<0mm,-4mm>*{\ldots}**@{},
 <0mm,0mm>*{};<4.5mm,-5mm>*{}**@{-},
 <0mm,0mm>*{};<8mm,-5mm>*{}**@{-},
   <0mm,0mm>*{};<-11mm,-7.9mm>*{^{\tau(1)}}**@{},
   <0mm,0mm>*{};<-4mm,-7.9mm>*{^{\tau(2)}}**@{},
     <0mm,0mm>*{};<10.0mm,-7.9mm>*{^{\tau(n)}}**@{},
 <0mm,0mm>*{};<0mm,5mm>*{}**@{-},
 \end{xy}
=0, \hskip 2em
1 \leq k \leq n-1,\ \mbox {and}
\end{equation}
\item[(ii)]
 planar $(0,n)$-corollas in degree $-n$
\[
\begin{xy}
 <0mm,0mm>*{\bullet};<0mm,0mm>*{}**@{},
 <0mm,0mm>*{};<-8mm,-5mm>*{}**@{-},
 <0mm,0mm>*{};<-4.5mm,-5mm>*{}**@{-},
 <0mm,0mm>*{};<0mm,-4mm>*{\ldots}**@{},
 <0mm,0mm>*{};<4.5mm,-5mm>*{}**@{-},
 <0mm,0mm>*{};<8mm,-5mm>*{}**@{-},
   <0mm,0mm>*{};<-11mm,-7.9mm>*{^{1}}**@{},
   <0mm,0mm>*{};<-4mm,-7.9mm>*{^{2}}**@{},
     <0mm,0mm>*{};<10.0mm,-7.9mm>*{^{n}}**@{},
 \end{xy}\ \ \ {n\geq 2,\ \sigma\in \bS_n},
\]
which are cyclic skew-symmetric
\begin{equation}
\nonumber
\begin{xy}
 <0mm,0mm>*{\bullet};<0mm,0mm>*{}**@{},
 <0mm,0mm>*{};<-8mm,-5mm>*{}**@{-},
 <0mm,0mm>*{};<-4.5mm,-5mm>*{}**@{-},
 <0mm,0mm>*{};<0mm,-4mm>*{\ldots}**@{},
 <0mm,0mm>*{};<4.5mm,-5mm>*{}**@{-},
 <0mm,0mm>*{};<8mm,-5mm>*{}**@{-},
   <0mm,0mm>*{};<-9mm,-7.9mm>*{^{1}}**@{},
   <0mm,0mm>*{};<-4mm,-7.9mm>*{^{2}}**@{},
     <0mm,0mm>*{};<10.0mm,-7.9mm>*{^{n}}**@{},
 \end{xy}
=(-1)^{\sgn(\xi)}
\begin{xy}
 <0mm,0mm>*{\bullet};<0mm,0mm>*{}**@{},
 <0mm,0mm>*{};<-8mm,-5mm>*{}**@{-},
 <0mm,0mm>*{};<-4.5mm,-5mm>*{}**@{-},
 <0mm,0mm>*{};<0mm,-4mm>*{\ldots}**@{},
 <0mm,0mm>*{};<4.5mm,-5mm>*{}**@{-},
 <0mm,0mm>*{};<8mm,-5mm>*{}**@{-},
   <0mm,0mm>*{};<-11mm,-7.9mm>*{^{\xi(1)}}**@{},
   <0mm,0mm>*{};<-4mm,-7.9mm>*{^{\xi(2)}}**@{},
     <0mm,0mm>*{};<10.0mm,-7.9mm>*{^{\xi(n)}}**@{},
 \end{xy}, \ \mbox {for}\
\xi = (1,\ldots,n),
\end{equation}
\end{itemize}
with the differential given as
\begin{eqnarray}
\label{comm-diff}
\partial
\begin{xy}
 <0mm,0mm>*{\bullet};<0mm,0mm>*{}**@{},
 <0mm,0mm>*{};<-8mm,-5mm>*{}**@{-},
 <0mm,0mm>*{};<-4.5mm,-5mm>*{}**@{-},
 <0mm,0mm>*{};<0mm,-4mm>*{\ldots}**@{},
 <0mm,0mm>*{};<4.5mm,-5mm>*{}**@{-},
 <0mm,0mm>*{};<8mm,-5mm>*{}**@{-},
   <0mm,0mm>*{};<-9mm,-7.9mm>*{^1}**@{},
   <0mm,0mm>*{};<-5mm,-7.9mm>*{^2}**@{},
   <0mm,0mm>*{};<10.0mm,-7.9mm>*{^n}**@{},
 <0mm,0mm>*{};<0mm,5mm>*{}**@{-},
 \end{xy}
&=&\sum_{k=0}^{n-2}\sum_{l=2}^{n-k}
(-1)^{k+l(n-k-l)+1}
\begin{xy}
<0mm,0mm>*{\bullet},
<0mm,5mm>*{}**@{-},
<4mm,-7mm>*{^{1\  \dots \  k\qquad\ \ k+l+1\ \ \dots\ \  n }},
<-14mm,-5mm>*{}**@{-},
<-6mm,-5mm>*{}**@{-},
<20mm,-5mm>*{}**@{-},
<8mm,-5mm>*{}**@{-},
<0mm,-5mm>*{}**@{-},
<0mm,-5mm>*{\bullet};
<-5mm,-10mm>*{}**@{-},
<-2mm,-10mm>*{}**@{-},
<2mm,-10mm>*{}**@{-},
<5mm,-10mm>*{}**@{-},
<0mm,-12mm>*{_{k+1\ \, \dots\ \, k+l}},
\end{xy},
\\
\label{wheeled-delta}
\partial
\begin{xy}
 <0mm,0mm>*{\bullet};<0mm,0mm>*{}**@{},
 <0mm,0mm>*{};<-8mm,-5mm>*{}**@{-},
 <0mm,0mm>*{};<-4.5mm,-5mm>*{}**@{-},
 <0mm,0mm>*{};<0mm,-4mm>*{\ldots}**@{},
 <0mm,0mm>*{};<4.5mm,-5mm>*{}**@{-},
 <0mm,0mm>*{};<8mm,-5mm>*{}**@{-},
   <0mm,0mm>*{};<-8mm,-7.9mm>*{^{1}}**@{},
   <0mm,0mm>*{};<-4mm,-7.9mm>*{^{2}}**@{},
     <0mm,0mm>*{};<10.0mm,-7.9mm>*{^{n}}**@{},
 \end{xy}
&=&
\sum_{i=1}^{n} (-1)^{(n-1)(i-1)}
\begin{xy}
<0mm,0mm>*{\bullet},
<-5mm,-5mm>*{}**@{-},
<-2mm,-5mm>*{}**@{-},
<2mm,-5mm>*{}**@{-},
<5mm,-5mm>*{}**@{-},
<0mm,-7mm>*{_{\hskip 1mm  i \hskip 3mm \cdots \hskip 3mm i-1}},
(0,0)*{} \ar@{}@(u,r) (0,0)*{}
\end{xy}
+\sum_{i=1}^{n}\sum_{l=1}^{n-2} (-1)^{(n-1)(i-1)+l}
\begin{xy}
<0mm,0mm>*{\bullet},
<0mm,-7mm>*{_{\hskip 3mm i \hskip 3mm \cdots \hskip 3mm i+l}\qquad},
<-8mm,-5mm>*{}**@{-},
<-5mm,-5mm>*{}**@{-},
<-1mm,-5mm>*{}**@{-},
<2mm,-5mm>*{}**@{-},
<8mm,-5mm>*{}**@{-},
<8mm,-5mm>*{\bullet};
<6mm,-10mm>*{}**@{-},
<9mm,-10mm>*{}**@{-},
<13mm,-10mm>*{}**@{-},
<16mm,-10mm>*{}**@{-},
<8mm,-12mm>*{_{\hskip 4mm i+l+1 \hskip 2mm \cdots \hskip 2mm i-1}},
\end{xy}.
\end{eqnarray}
\end{theorem-B}
The labels in
the right hand side of~(\ref{wheeled-delta})
denote cyclic permutations, so that for example
$i \cdots i-1$ in the first term means $i, \ldots
n,1, \ldots ,i-1$.

By contrast, the operad \mbox{\sf Lie} is rigid
with respect to the wheelification: it was proven in \cite{me3} that
$(\mbox{\sf Lie}^\circlearrowright)_\infty
 = \ (\mbox{\sf Lie}_\infty)^\circlearrowright$.

A conceptual understanding of homological properties of a wide
class of wheeled \PROP{s} is provided by
wheeled quadratic duality and wheeled Koszulness which we set up in
Sections~\ref{pekelna_hudba} and~\ref{Tomicek}. This theory generalizes
its ``unwheeled'' precursor developed
in~\cite{ginzburg-kapranov:DMJ94}.

As an application of the above results we compute cohomology of a
directed version of Kontsevich's complex of ribbon graphs.  Let
$\fG_g$ be the linear span of (not necessary connected) directed
ribbon graphs of genus $g$ such that (i) each vertex has at least
three attached internal edges of which at least one is incoming and
at least one is outgoing, (ii) vertices in a closed path all have
either precisely one incoming edge, or precisely one outgoing edge,
and (iii) there are no input and output legs (i.e.\ every edge of
the graph is internal)\footnote{The space $\fG_g$ is obviously zero
in the class of directed graphs {\em without wheels}. However it is
highly non-trivial with {\em wheels allowed}. This simple fact
provides us with one more motivation for introducing the category of
wheeled PROPs.}. For example,
$$
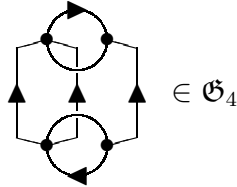
\xy
*\xycircle(4,4){-};
<0mm,4mm>*{\blacktriangleright}; <-4mm,0mm>*{\bullet};
<-4mm,0mm>*{};<0mm,-1mm>*{}**@{-},
<0mm,-1mm>*{};<0mm,-3.5mm>*{}**@{-},
<0mm,-4.5mm>*{};<0mm,-9.4mm>*{}**@{-},
<0mm,-6.5mm>*{};<0mm,-9.4mm>*{}**@{-}, <0mm,-7mm>*{\blacktriangle};
<-4mm,0mm>*{};<-8mm,-1mm>*{}**@{-},
<-8mm,-1mm>*{};<-8mm,-13mm>*{}**@{-},
<-8mm,-13mm>*{};<-4mm,-14mm>*{}**@{-},
<-8mm,-7mm>*{\blacktriangle};
<0mm,-14mm>*{}*\xycircle<11.5pt>{-};
<0mm,-18mm>*{\blacktriangleleft}; <-4mm,-14mm>*{\bullet};
<-4mm,-14mm>*{};<0mm,-13mm>*{}**@{-},
<0mm,-13mm>*{};<0mm,-10.5mm>*{}**@{-},
 <4mm,0mm>*{\bullet};
  <4mm,-14mm>*{\bullet};
  <4mm,0mm>*{};<8mm,-1mm>*{}**@{-},
  <4mm,-14mm>*{};<8mm,-13mm>*{}**@{-},
  <8mm,-7mm>*{\blacktriangle};
 <8mm,-13mm>*{};<8mm,-1mm>*{}**@{-},
 <17mm,-7mm>*{\in \fG_4};
\endxy
$$
Orientation on an element $G\in \fG_g$ is, by definition, an
orientation of the vector space, $\R^{v(G)}$, spanned by the set,
$v(G)$, of vertices of $G$ (which is in fact the same as the
orientation of $\R^{e(G)}\oplus H_\bullet(|G|,\R)$, where $e(G)$ is
the set of internal edges of $G$ and $H_\bullet(|G|,\R)$ is the
homology of $G$ viewed as a 1-dimensional CW complex). The space
$\fG_g$ is naturally a cochain complex, $(\fG_g, \p)$, with respect
to the grading,
$$
\fG_g=\bigoplus_{n\geq 0} \fG_g^n,
$$
by the number, $n=|v(G)|$, of vertices of its elements $G$
\cite{konts, MV}.

In the following theorem, $\uparrow^k$ denotes, for $k \geq 0$,
the suspension of a graded vector space iterated $k$-times.

\begin{theorem-C}
Let $\Gamma\langle E \rangle=\{\Gamma\langle E\rangle
(m,n)\}_{m,n\geq 0}$ be the
free (ordinary) PROP generated by the $\bS$-bimodule $E=\{E(m,n)\}_{m,n\geq 0}$,
with
\[
E(m,n)=\left\{\Ba{ll}
\underset{p,q\geq 0\atop p+q=m}{\bigoplus}\uparrow^m G(p,q)
\ \oplus\   \underset{p,q\geq 1\atop p+q=m}
{\bigoplus}\uparrow^{m-1}G(p,q) & \mbox{for}\
m\geq 1, n=0, \vspace{2mm}\\
\underset{p,q\geq 0\atop p+q=n}{\bigoplus}\uparrow^{n}G(p,q)
\ \oplus\   \underset{p,q\geq 1\atop p+q=n}
{\bigoplus}\uparrow^{n-1}G(p,q) & \mbox{for}\
n\geq 1, m=0, \vspace{2mm} \\
0 & \mbox{otherwise.}
\Ea
\right.
\]
where $G(p,q):= \R[\bS_{p+q}]^{C_p\times C_q}$ is the space of coinvariants
with respect to the product of cyclic subgroups
$C_p\times C_q\subset \bS_{p+q}$ generated by permutations $(12\ldots p)$
and $(p+1\ldots p+q)$.

Then the graded vector space $H^\bullet(\fG_g)$ is isomorphic to the subspace of the component
$\Gamma\langle E\rangle(0,0)$ generated by graphs with $g-1$ internal edges.
In particular, $\bigoplus_g H^\bullet(\fG_g)
\simeq \Gamma\langle E\rangle(0,0)$.
\end{theorem-C}

\begin{corollary-D}
The vector space
$H^n(\fG_g)$ is nonzero only for $n$ in the range
\ $
g-\frac{1}{2}(1-(-1)^g)\leq n \leq 2(g-1)$.
\end{corollary-D}

It is worth  noting another interesting phenomenon of
substantial and sometimes highly non-trivial change of the set of
{\em morphisms}\, between ordinary PROPs under their wheeled
completions. For example, it was shown in
\cite{me3} that deformation quantization can be understood as a
morphism,
$
\mbox{{\sf DefQ}} \longrightarrow
\mbox{\sf PolyV}^\circlearrowright,
$
between the dg
PROP, { \sf DefQ}, of star products and the {\em wheeled completion}\,
of the dg PROP, {\sf PolyV}, of polyvector fields. {\em No}\, such
morphism (satisfying the quasi-classical limit condition),
$\mbox{{\sf DefQ}} \longrightarrow \mbox{\sf PolyV}$,
exists for the original dg PROPs { \sf DefQ} and {\sf PolyV} within
the category of ordinary PROPs.

\sip

Here is an itemized and more detailed  list of main results of our paper:

In \S 2 we construct a triple over the category of $\bS$-bimodules using
directed graphs with directed cycles,  and then define wheeled PROPs as algebras over
that triple. Modifications of this notion (such as wheeled properads and modular wheeled properads)
are also given.

In \S 3 we construct a particular example of dg wheeled PROP
whose representations are in 1-1 correspondence with the set of solutions of master equations
which describe formal germs of so called $SP$-manifolds, key
geometric objects
in the theory of Batalin-Vilkovisky quantization (cf.~\cite{Sch}).

In \S 4 we introduce the notion of wheeled coproperads,
define wheeled bar and cobar functors (which, rather
surprisingly, turn out to be much simpler than their analogues for
(co)properads) and prove Theorems~\ref{cobar} and~\ref{ll}
on bar+cobar resolutions of wheeled properads and coproperads.

In \S 5 we study quadratic wheeled operads and introduce the notion of wheeled
Koszulness. Theorems~A and B~above
imply that classical wheeled operads $\Com^\circlearrowright$,
$\Ass^\circlearrowright$ are wheeled Koszul (while the wheeled Koszulness of
$\Lie^\circlearrowright$ was proved in~\cite{me3}).

In \S 6 we prove Theorems A, C and D (in fact we prove a stronger
Theorem~\ref{hg} of which Theorem C is
a corollary),  and, as a concrete and natural example of
$\Ass_\infty^\circlearrowright$-structure  we compute  explicitly new trace-type Massey operations
on the homology of an arbitrary finite-dimensional dg associative algebra. We also find  a new
{\em cyclic characteristic class}\, of an arbitrary finite-dimensional $\Ass_\infty$-algebra and prove
that its vanishing is a necessary and sufficient condition for extendability of the given
$\Ass_\infty$-structure to $\Ass_\infty^\circlearrowright$ one.

Finally, in \S 7 we compute cohomology of
the dg properad $(\Com_\infty)^\circlearrowright$ and prove Theorem~B.

\section{Wheeled PROPs}
\label{Chicago}

\subsection{Basic definitions.}
Let $\bfk$ denote a ground field which will always be assumed of
characteristic zero.  Recall that a dg (differential graded) \PROP\ is
a collection $\sfP = \{\sfP(m,n)\}$, $m,n \geq 0$, of dg
$(\Sigma_m,\Sigma_n)$-bimodules (left $\Sigma_m$- right
$\Sigma_n$-modules such that the left action commutes with the right
one), together with two types of compositions, {\em horizontal\/}
\begin{equation}
\label{hor}
\ot : \sfP(m_1,n_1) \ot \cdots \ot \sfP(m_s,n_s) \to
\sfP(m_1+\cdots +m_s,n_1+\cdots +n_s),
\end{equation}
defined for all $m_1,\ldots,m_s,n_1,\ldots,n_s \geq 0$, and {\em vertical\/}
\begin{equation}
\label{vert}
\circ : \sfP(m,n) \ot \sfP(n,k) \to \sfP(m,k),
\end{equation}
defined for all $m,n,k \geq 0$. These compositions are compatible with the dg
structures. One also assumes the existence of
a {\em unit\/} $\id \in \sfP(1,1)$.

\PROP{s} should satisfy axioms which could be
read off from the example of the {\em
endomorphism \PROP\/} $\End_V$ of a vector space $V$, with
$\End_V(m,n)$ the space of linear maps
$\Hom(\otexp Vn,\otexp Vm)$ with $n$ `inputs' and $m$ `outputs,'
$\id \in \End_V(1,1)$ the identity map,
horizontal composition given by the tensor product of linear maps, and
vertical composition by the ordinary composition of linear
maps. For a precise definition
see~\cite{maclane:RiceUniv.Studies63,markl:JPAA96}.

Recall also that a {\em $\sfP$-algebra\/} is a \PROP\
homomorphism $\rho : \sfP \to \End_V$. It is determined by a system
\[
\alpha: \sfP(m,n) \otimes \otexp Vn \to \otexp Vm,\ m,n, \geq 0,
\]
of linear maps satisfying appropriate axioms. $\sfP$-algebras are
sometimes called {\em representations\/} of~$\sfP$.

\PROP{s} are devices that describe structures consisting of operations
with several inputs and several outputs. Therefore various
bialgebras (associative, Lie, infinitesimal) are \PROP{ic}
algebras. In this section we introduce wheeled \PROP{s} that
generalize \PROP{s} in that they describe structures whose axioms involve
also ``contraction of indices'' or ``traces.''

Before we give a precise definition of wheeled \PROP{s},
we present the following fundamental example of this kind of structure.

\begin{example}
\label{fund}
{\rm
Let $V$ be a finite-dimensional vector space, with a basis
$\{e_s\}_{1 \leq s \leq d}$. Consider a linear map
$f: \otexp Vn \to \otexp Vm \in
\End_V(m,n)$, $m,n \geq 1$,  that acts on
products of basic elements~by
\[
f(e_{\alpha_1} \ot \cdots \ot e_{\alpha_n})
= \sum M^{\rada {\beta_1}{\beta_m}}_{\rada {\alpha_1}{\alpha_n}}
e_{\beta_1} \ot \cdots \ot e_{\beta_m},
1 \leq \alpha_s \leq d,\ 1 \leq s \leq n,
\]
where $M^{\rada {\beta_1}{\beta_m}}_{\rada {\alpha_1}{\alpha_n}} \in
\bfk$ are scalars and the sum is taken over all $1 \leq
\beta_t \leq d$, $1 \leq t \leq m$.
For any pair of indices $i$ and $j$, $1 \leq i \leq m$, $1 \leq j \leq n$,
define a multilinear map $\xi^i_j(f) : \otexp V{{(n-1)}} \to \otexp
V{{(m-1)}}$~by
\begin{equation}
\label{psano_na_''novem''_pocitaci}
\xi^i_j(f)(e_{\mu_1} \ot \cdots \ot e_{\mu_{n-1}})
:=
\sum M^{
\rada {\nu_1}{\nu_{i-1}},
\gamma,
\rada{\nu_i}{\nu_{m-1}}
}_{\rada {\mu_1}{\mu_{j-1}},\gamma,\rada{\mu_j}{\mu_{n-1}}}
e_{\nu_1} \ot \cdots \ot e_{\nu_{m-1}},
\end{equation}
where $1 \leq \mu_u \leq d$, $1 \leq u \leq n-1$, and
the summation runs over all $1 \leq \nu_v,\gamma \leq d$,
$1 \leq v \leq m-1$. Remarkably,
the above definition of $\xi^i_j(f)$ does not depend on the
basis. Formula~(\ref{psano_na_''novem''_pocitaci}) therefore defines
a linear map
\begin{equation}
\label{zase_nefunguje}
\xi_j^i : \End_V(m,n) \to \End_V(m-1,n-1).
\end{equation}

A ``coordinate-free'' definition can
be given as follows. Using the duality in the category of
finite-dimensional vector spaces, one can associate to any map $f: \otexp
Vn \to \otexp Vm$ a map
\[f^i_j : \Hom(\otexp V{(n-1)},\otexp
V{(m-1)})  \to \Hom(V,V),
\]
by singling out the $i$th output and the
$j$th input of $f$. The composition with the ordinary trace $\Tr :\Hom(V,V)
\to \bfk$ is a map  $\Tr (f^i_j) :
\Hom(\otexp V{(n-1)},\otexp
V{(m-1)})  \to \bfk$ which in turn corresponds, via the duality, to the
contraction $\xi^j_i(f) : \otexp V{(n-1)} \to \otexp
V{(m-1)}$ constructed above.
}
\end{example}

Wheeled \PROP{s} are \PROP{s} that, besides the horizontal and vertical
compositions~(\ref{hor}) and ~(\ref{vert}) admit also {\em contractions\/}
\[
\xi_j^i : \sfP(m,n) \to \sfP(m-1,n-1),\ m,n \geq 1, 1\leq i \leq m,\ 1
\leq j \leq n,
\]
that generalize~(\ref{zase_nefunguje}) in Example~\ref{fund}.
A precise definition of wheeled \PROP{s} will be a
modification of the unbiased definition of ordinary
\PROP{s} given in~\cite[Section~8]{handbook} which we review below. The
difference between biased and unbiased definitions is explained
in~\cite{leinster:book}, see also a remark in Section~3 of~\cite{handbook}.

Recall that
a {\em $\Sigma$-bimodule\/} is a system $E= \bicol E0$ such that each
$E(m,n)$ is a left $\bfk[\Sigma_m]$- right $\bfk[\Sigma_n]$-bimodule. Let
$\sigmabimod$ denote the category of $\Sigma$-bimodules.
For $E \in \sigmabimod$ and finite sets $Y,X$ with $m$ resp.~$n$
elements put
\[
E(Y,X) := \Bij(Y,[m]) \times_{\Sigma_m} E(m,n) \times_{\Sigma_n}
\Bij([n],X),\ m,n \geq 0,
\]
where
\[
\Bij(T,S) := \{\vartheta : S \stackrel{\cong}{\longrightarrow} T\}
\]
is the set of bijections between finite sets $S$ and $T$, and
$[m] := \{\rada 1m\}$, $[n] := \{\rada 1n\}$ as
usual. We are going to define
an endofunctor $\Gamma$ on the category $\sigmabimod$
that assigns to each $\Sigma$-bimodule $E$ the
$\Sigma$-bimodule $\freePROP(E)$
of $E$-decorated graphs. Unfortunately, the naive concept of graph is
not refined enough for our purposes and we must recall the following
more sophisticated concept taken from~\cite{kontsevich-manin:CMP94}.

\begin{definition}
\label{graph-def}
A {\em graph\/} $G$ is a finite set $\Flag(G)$
(whose elements are called {\em flags\/} or {\em half-edges\/})
together with an involution
$\sigma$ and a partition $\lambda$.
The {\em vertices\/} $\vert(G)$
of a graph $G$ are the blocks of the partition
$\lambda$; we assume that the number of these blocks is finite.
The {\em edges\/}
$\Edg(G)$ are pairs of flags forming a two-cycle of $\sigma$. The
{\em legs\/} $\leg(G)$ are the fixed points
of $\sigma$.
\end{definition}

We also denote by $\edge(v)$ the flags belonging to the block $v$ or,
in common speech, half-edges adjacent to the vertex $v$.  We say that
graphs $G_1$ and $G_2$ are {\em isomorphic\/} if there
exists a set isomorphism $\varphi : \Flag(G_1) \to
\Flag(G_2)$ that preserves the partitions and commutes with the
involutions. We may associate to a graph $G$ a finite
one-dimensional cell complex $|G|$, obtained by taking one copy
of $[0,\frac 12]$ for each flag, a point for each block of the
partition, and imposing the following equivalence
relation: The points $0\in [0,\frac 12]$ are identified for all flags
in a block of the partition $\lambda$ with the point corresponding to
the block, and the points $\frac 12 \in
[0,\frac 12]$ are identified for pairs of flags exchanged by the
involution $\sigma$.

\begin{figure}[t]
\begin{center}
{
\unitlength=.800000pt
\begin{picture}(240.00,170.00)(0.00,0.00)
\thicklines
\put(220.00,160.00){\makebox(0.00,0.00){$5$}}
\put(20.00,60.00){\makebox(0.00,0.00){$4$}}
\put(80.00,160.00){\makebox(0.00,0.00){$3$}}
\put(50.00,160.00){\makebox(0.00,0.00){$2$}}
\put(160.00,170.00){\makebox(0.00,0.00){$1$}}
\put(160.00,70.00){\makebox(0.00,0.00){$3$}}
\put(200.00,70.00){\makebox(0.00,0.00){$2$}}
\put(100.00,80.00){\makebox(0.00,0.00){$1$}}
\put(200.00,30.00){\makebox(0.00,0.00){$\bullet$}}
\put(160.00,0.00){\makebox(0.00,0.00){$\bullet$}}
\put(130.00,30.00){\makebox(0.00,0.00){$\bullet$}}
\put(160.00,50.00){\makebox(0.00,0.00){$\bullet$}}
\put(160.00,100.00){\makebox(0.00,0.00){$\bullet$}}
\put(100.00,110.00){\makebox(0.00,0.00){$\bullet$}}
\put(220.00,120.00){\makebox(0.00,0.00){$\bullet$}}
\put(180.00,120.00){\makebox(0.00,0.00){$\bullet$}}
\put(200.00,100.00){\makebox(0.00,0.00){$\bullet$}}
\put(160.00,140.00){\makebox(0.00,0.00){$\bullet$}}
\put(60.00,50.00){\makebox(0.00,0.00){$\bullet$}}
\put(40.00,30.00){\makebox(0.00,0.00){$\bullet$}}
\put(60.00,110.00){\makebox(0.00,0.00){$\bullet$}}
\put(80.00,130.00){\makebox(0.00,0.00){$\bullet$}}
\put(40.00,130.00){\makebox(0.00,0.00){$\bullet$}}
\put(200.00,20.00){\vector(0,1){10.00}}
\qbezier(220.00,10.00)(200.00,10.00)(200.00,30.00)
\qbezier(240.00,30.00)(240.00,10.00)(220.00,10.00)
\qbezier(220.00,50.00)(240.00,50.00)(240.00,30.00)
\qbezier(200.00,30.00)(200.00,50.00)(220.00,50.00)
\put(160.00,0.00){\vector(4,3){40.00}}
\put(160.00,0.00){\vector(-1,1){30.00}}
\put(160.00,0.00){\vector(0,1){50.00}}
\put(130.00,30.00){\vector(3,2){30.00}}
\put(170.00,140.00){\vector(-1,0){10.00}}
\put(150.00,140.00){\vector(1,0){10.00}}
\put(200.00,80.00){\vector(0,1){20.00}}
\put(220.00,120.00){\vector(0,1){30.00}}
\put(200.00,100.00){\vector(1,1){20.00}}
\put(200.00,100.00){\vector(-1,1){20.00}}
\put(160.00,80.00){\vector(0,1){20.00}}
\put(160.00,140.00){\vector(0,1){20.00}}
\qbezier(160.00,100.00)(140.00,100.00)(140.00,120.00)
\qbezier(180.00,120.00)(180.00,100.00)(160.00,100.00)
\qbezier(160.00,140.00)(180.00,140.00)(180.00,120.00)
\qbezier(140.00,120.00)(140.00,140.00)(160.00,140.00)
\qbezier(80.00,70.00)(60.00,70.00)(60.00,50.00)
\qbezier(100.00,50.00)(100.00,70.00)(80.00,70.00)
\qbezier(80.00,30.00)(100.00,30.00)(100.00,50.00)
\qbezier(60.00,50.00)(60.00,30.00)(80.00,30.00)
\put(60.00,40.00){\vector(0,1){10.00}}
\put(0.00,140.00){\line(0,-1){120.00}}
\qbezier(20.00,0.00)(0.00,0.00)(0.00,20.00)
\qbezier(40.00,20.00)(40.00,0.00)(20.00,0.00)
\qbezier(20.00,160.00)(0.00,160.00)(0.00,140.00)
\qbezier(40.00,140.00)(40.00,160.00)(20.00,160.00)
\put(40.00,130.00){\vector(1,2){10.00}}
\put(100.00,90.00){\vector(0,1){20.00}}
\put(100.00,110.00){\vector(-1,1){20.00}}
\put(80.00,130.00){\vector(0,1){20.00}}
\put(80.00,130.00){\vector(0,1){0.00}}
\put(60.00,110.00){\vector(1,1){20.00}}
\put(60.00,50.00){\vector(0,1){60.00}}
\put(60.00,110.00){\vector(-1,1){20.00}}
\put(40.00,30.00){\vector(1,1){20.00}}
\put(40.00,30.00){\vector(-1,1){20.00}}
\put(40.00,130.00){\line(0,1){10.00}}
\put(40.00,20.00){\vector(0,1){10.00}}
\put(30.00,90.00){\makebox(0.00,0.00){$w_1$}}
\put(82.00,48.00){\makebox(0.00,0.00){$w_2$}}
\put(222.00,28.00){\makebox(0.00,0.00){$w_3$}}
\put(161.00,119.00){\makebox(0.00,0.00){$w_4$}}
\end{picture}}
\end{center}
\caption{\label{fg1}
A directed $(5,3)$-graph with four independent wheels
$w_1$, $w_2$, $w_3$ and $w_4$.
The graph is not connected and has a component with no legs.
}
\end{figure}
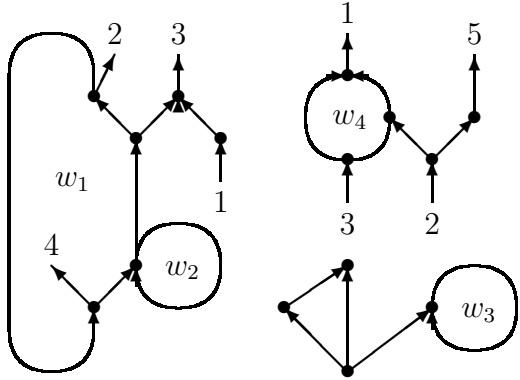

We call $|G|$ the {\em geometric
realization\/} of $G$.
Observe that empty blocks of the partition generate isolated vertices
in the geometric realization.
We will usually make no distinction
between the graph and its geometric realization.
See~\cite[Section~II.5.3]{markl-shnider-stasheff:book}
or~\cite[Section~7]{handbook} for more details.
A~graph $G$ as in Definition~\ref{graph-def} is a
{\em directed $(m,n)$-graph\/} if
\begin{itemize}
\item[(i)]
each edge of $G$ is equipped with a direction
\item[(ii)]
the set of legs of $G$ is divided into the set of inputs labeled by
$\{1, \dots, n\}$ and the set of outputs labeled by $\{1, \dots, m\}$.
\end{itemize}
The direction of edges together with the division
of legs into inputs and outputs determines at
each vertex $v \in \vert(G)$ of a directed graph $G$ a disjoint decomposition
\[
\edge(v) = \In(v) \sqcup \Out(v)
\]
of the set of edges adjacent to $v$ into the set $\In(v)$ of incoming
edges and the set $\Out(v)$ of outgoing edges.  The pair
$(\#(\Out(v)),\#(\In(v))) \in {\mathbb N} \times {\mathbb N}$ is
called the {\em biarity\/} of $v$.  By a {\em wheel\/} in a directed
graph we mean a directed cycle of edges, see Figure~\ref{fg1}. We
usually draw directed graphs in such a way that output edges of
vertices point upwards and input edges enter vertices from the bottom.
We denote by $\DGr mn$ the category of directed $(m,n)$-graphs
\underline{without} \underline{wheels} and
their isomorphisms.

To incorporate the \PROP{eradic} unit, we assume that
$\DGr mn$, for $m=n$, contains also the
{\em exceptional graph\/}
\[
\uparrow \ \uparrow \ \uparrow \ \cdots \uparrow\ \in \DGr nn,\ n \geq 1,
\]
with $n$ inputs, $n$ outputs and no vertices.
For a graph $G \in \DGr mn$ and a $\Sigma$-bimodule
$E$, let
\begin{equation}
\label{pomalu}
E(G) := \bigotimes_{v \in \vert(G)} E(\Out(v),\In(v))
\end{equation}
be the linear space of all decorations of vertices of the graph $G$ by
elements of $E$. Since the assignment $E \mapsto E(G)$
clearly defines a functor from the category $\DGr mn$ to the category of
vector spaces, it makes sense to define
\begin{equation}
\label{jsem_zvedav_jestli_se_ta_druha_Jitka_ozve}
\freePROP(E)(m,n) := \colim{{G \in \DGr mn}}{{E}}(G),\
m,n \geq 0.
\end{equation}

Denote finally $\freePROP(E)$ the $\Sigma$-bimodule $\bicol{\freePROP(E)}0$.
Our aim now is to explain
that the correspondence $E \mapsto \freePROP(E)$
defines a triple on the category of
$\Sigma$-bimodules such that  \PROP{s} are algebras over this
triple.

The concept of triples and their algebras is classical, so we recall it only
briefly.
Let ${\it End}({\calC})$ be the
strict symmetric monoidal category of endofunctors on a category
$\calC$ where multiplication is the composition of functors.

\begin{definition}
\index{triple}\label{deftriple}\index{monad}
A {\em triple} (also called a {\em monad\/}) $T$ on a category $\calC$ is an
associative and unital monoid $\left( T,\mu ,\upsilon \right) $ in
${\it End}(\calC)$. The multiplication $\mu : TT \to T$
and unit morphism $\upsilon:  {\it id} \to T$
satisfy the axioms given by commutativity of the diagrams in
Figure~\ref{triple}.
\end{definition}

Let us indicate how to construct transformations
\[
\mu  : \Gamma \Gamma \to \Gamma
\mbox { and }
\upsilon : \id \to \Gamma
\]
making $\Gamma$ a triple on the category of $\Sigma$-bimodules.
Let us start with the triple multiplication $\mu$.
It follows from
definition~(\ref{jsem_zvedav_jestli_se_ta_druha_Jitka_ozve}) that,
for each $\Sigma$-bimodule $E$,
\begin{equation}
\label{jeste}
\freePROP({\freePROP(E)})(m,n) :=
\colim{{G \in \DGr mn}}{\freePROP(E)}(G),\
m,n \geq 0.
\end{equation}
The elements in the right hand side are represented by directed graphs
with vertices decorated by elements of $\freePROP(E)$, while elements
of $\freePROP(E)$ are represented by directed graphs with vertices
decorated by $E$.  We may therefore imagine elements of
$\freePROP({\freePROP(E)})$ as `bracketed' or `nested' $E$-decorated
directed graphs $G$, with nests encompassing groups of vertices of $G$
that represent one vertex decorated by an element of $\freePROP(E)$.
See also~\cite[Section~5]{handbook}
where this nesting is described and
analyzed in detail.  The triple multiplication
$\mu_E : \freePROP({\freePROP(E)}) \to \freePROP(E)$ then simply forgets
the nests.  The triple unit $\upsilon_E : E \to \freePROP(E)$
identifies elements of $E$ with decorated corollas:
\[
\unitlength 4mm
\linethickness{0.4pt}
\begin{picture}(20,7.5)(10,16.8)
\put(14,20){\makebox(0,0)[cc]{$E(m,n)\ni e \hskip 10pt \longleftrightarrow$}}
\multiput(20,20)(-.0333333,-.0333333){60}{\line(0,-1){.0333333}}
\multiput(20,20)(-.033333,-.066667){30}{\line(0,-1){.066667}}
\multiput(20,20)(.0333333,-.0333333){60}{\line(0,-1){.0333333}}
\put(20,20){\makebox(0,0)[cc]{$\bullet$}}
\put(20.5,18){\makebox(0,0)[cc]{$\ldots$}}
\put(20,17){\makebox(0,0)[cc]{%
   $\underbrace{\rule{16mm}{0mm}}_{\mbox{\small $n$ inputs}}$}}
\put(0,40){
\multiput(20,-20)(-.0333333,.0333333){60}{\line(0,1){.0333333}}
\multiput(20,-20)(-.033333,.066667){30}{\line(0,1){.066667}}
\multiput(20,-20)(.0333333,.0333333){60}{\line(0,1){.0333333}}
\put(20.5,-18){\makebox(0,0)[cc]{$\ldots$}}
\put(20,-17){\makebox(0,0)[cc]{%
   $\overbrace{\rule{16mm}{0mm}}^{\mbox{\small $m$ outputs}}$}}
}
\put(20.5,20){\makebox(0,0)[l]{$e$}}
\put(22,20){\makebox(0,0)[cl]{$\in \freePROP(E)(m,n),\ m,n \geq 0$.}}
\end{picture}
\]

It is not difficult to verify that the above constructions indeed make
$\Gamma$ a
triple,
compare~\cite[Section~II.1.12]{markl-shnider-stasheff:book} or
~\cite[Section~5]{handbook}.
The last thing we need to recall is:

\begin{definition}
\label{snad_projdu}
A  {\em $T$-algebra\/} or {\em algebra over the triple $T$\/}
is an object $A$ of\/ $\calC$  together with a structure morphism
$\alpha :T(A)\rightarrow A$
satisfying
\[
\alpha (T(\alpha ))=\alpha(\mu_A)  \mbox{ and }
\alpha \upsilon_A=\id_{A},
\]
see Figure~\ref{Talgebra}.
\end{definition}
\begin{figure}
\begin{center}
{\setlength{\unitlength}{.66cm}
\begin{picture}(15,4)(2,0)
%
\put(2,3){\makebox(0,0){$TTT$}}
\put(6,3){\makebox(0,0){$TT$}}
\put(2,0){\makebox(0,0){$TT$}}
\put(6,0){\makebox(0,0){$T$}}
\put(1.7,1.5){\makebox(0,0)[r]{$\mu T$}}
\put(6.3,1.5){\makebox(0,0)[l]{$\mu$}}
\put(4.1,0.5){\makebox(0,0)[b]{$\mu$}}
\put(4.1,3.5){\makebox(0,0)[b]{$T\mu$}}
\put(3,0){\vector(1,0){2}}
\put(3.2,3){\vector(1,0){1.8}}
\put(2,2.5){\vector(0,-1){2}}
\put(6,2.5){\vector(0,-1){2}}
%
\put(-.5,0){
\put(9.6,3){\makebox(0,0){$T$}}
\put(12.3,3){\makebox(0,0){$TT$}}
\put(10,3){\vector(1,0){1.8}}
\put(11,3.5){\makebox(0,0){$T\upsilon$}}
\put(12,2.5){\vector(-1,-3){0.7}}
\put(10,2.5){\vector(1,-3){0.7}}
\put(9.5,1.5){\makebox(0,0){$\id$}}
\put(12,1.5){\makebox(0,0){$\mu$}}
\put(11,0){\makebox(0,0){$T$}}
}
%
\put(14.6,3){\makebox(0,0){$T$}}
\put(17.3,3){\makebox(0,0){$TT$}}
\put(15,3){\vector(1,0){1.8}}
\put(15.8,3.5){\makebox(0,0){$\upsilon T$}}
\put(17,2.5){\vector(-1,-3){0.7}}
\put(15,2.5){\vector(1,-3){0.7}}
\put(14.5,1.5){\makebox(0,0){$\id$}}
\put(17,1.5){\makebox(0,0){$\mu$}}
\put(16,0){\makebox(0,0){$T$}}
\end{picture}}
\end{center}
\caption{\label{triple}
Associativity and unit axioms for a triple.}
\end{figure}
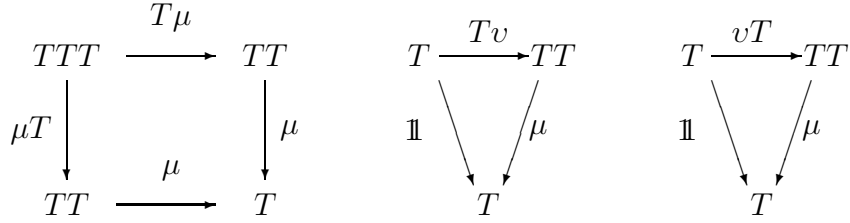

\begin{figure}
\begin{center}
{\setlength{\unitlength}{.66cm}
\begin{picture}(10,4)(2,0)
%
\put(2,3){\makebox(0,0){$T(T(A))$}}
\put(6,3){\makebox(0,0){$T(A)$}}
\put(2,0){\makebox(0,0){$T(A)$}}
\put(6,0){\makebox(0,0){$A$}}
\put(1.7,1.5){\makebox(0,0)[r]{$\mu$}}
\put(6.3,1.5){\makebox(0,0)[l]{$\alpha$}}
\put(4,0.5){\makebox(0,0)[b]{$\alpha$}}
\put(4.2,3.5){\makebox(0,0)[b]{$T(\alpha)$}}
\put(3,0){\vector(1,0){2}}
\put(3.2,3){\vector(1,0){1.8}}
\put(2,2.5){\vector(0,-1){2}}
\put(6,2.5){\vector(0,-1){2}}
%
\put(9.6,3){\makebox(0,0){$A$}}
\put(12.5,3){\makebox(0,0){$T(A)$}}
\put(10,3){\vector(1,0){1.6}}
\put(10.9,3.5){\makebox(0,0){$\upsilon_A$}}
\put(12,2.5){\vector(-1,-3){0.7}}
\put(10,2.5){\vector(1,-3){0.7}}
\put(9.5,1.5){\makebox(0,0){$\id$}}
\put(12.2,1.5){\makebox(0,0){$\alpha$}}
\put(11,0){\makebox(0,0){$A$}}
\end{picture}}
\end{center}
\caption{\label{Talgebra}
$T$-algebra structure.}
\end{figure}

The following proposition follows from~\cite[Section~8]{handbook},
see also~\cite{vallette:thesis,vallette:CMR04}.

\begin{proposition}
\label{unbiased-props}
\PROP{s} are algebras over the triple $\Gamma$.
\end{proposition}

As explained in~\cite{handbook}, the above proposition means that a
\PROP\ is given by coherent compositions along directed graphs
without wheels. Therefore, a \PROP\ is a $\Sigma$-bimodule $\sfP =
\bicol {\sfP}0$ equipped with a coherent system of linear maps
\begin{equation}
\label{7}
\alpha_G : \sfP(G) \to \sfP(m,n),\
G \in \DGr mn,\ m,n \geq 0,
\end{equation}
where $\sfP(G)$ is the space of $\sfP$-decorations of the graph $G$ introduced
in~(\ref{pomalu}).

Another important observations is that the triple multiplication $\mu
: \Gamma \Gamma \to \Gamma$ makes
$\freePROP(E)$ a \PROP, for each $\Sigma$-bimodule $E$ --
the maps~(\ref{7}) are the compositions
\[
\freePROP(E)(G) \stackrel{\iota_G}{\longrightarrow}
\freePROP({\freePROP(E)}) \stackrel{\mu_E}{\longrightarrow}
\freePROP(E),
\]
with $\iota_G$ the canonical map to the colimit in the right hand
side of~(\ref{jeste}).
It can be easily shown that the vertical composition in $\freePROP(E)$ is
given by the disjoint union of graphs,
the horizontal composition by grafting the
legs of graphs, and the unit is
the exceptional graph $\uparrow \hskip .2em \in
\freePROP(E)(1,1)$.  The following proposition follows from
general properties of algebras over triples~\cite{e-m}.

\begin{proposition}
\label{a_pratele_si_znechutim}
The \PROP\ $\freePROP(E)$ is the {\em free \PROP\/} generated by the
$\Sigma$-bimodule $E$.
\end{proposition}

Now we are finally ready to modify the above constructions and
introduce wheeled \PROP{s}, by allowing wheels in directed graphs.  We
start by denoting $\DGrw mn$, $m,n \geq 0$, the category of {\em
all\/} directed $(m,n)$-graphs and their isomorphisms. The little
oriented circle in $\DGrw mn$ indicates that wheels are allowed now.
The category $\DGrw nn$ contains, for each $n \geq 0$, the {\em
exceptional graphs\/}
\begin{equation}
\label{exceptional}
\uparrow \ \uparrow \ \uparrow \ \cdots \uparrow\
\hgw\hgw\cdots \hgw
\in \DGrw nn,\ n \geq 0,
\end{equation}

We denote
\[
\freeWPROP(E)(m,n) := \colim{{G \in \DGrw mn}}{{E}}(G),\
m,n \geq 0,
\]
with $E(G)$ given as in~(\ref{pomalu}). We then argue as before
that the above formula defines a triple $\Gamma^\circlearrowright$ on the
category of $\Sigma$-bimodules. The central definition of this section is:

\begin{definition}
{\em Wheeled \PROP{s}\/} are algebras over the triple $\Gamma^\circlearrowright :
\sigmabimod \to \sigmabimod$.
\end{definition}

The above definition can be reformulated by saying that
a wheeled \PROP\ is a $\Sigma$-bimodule
$\sfP = \bicol {\sfP}0$ with a coherent system of linear maps
\begin{equation}
\label{karta}
\alpha_G : \sfP(G) \to \sfP(m,n),\
G \in \DGrw mn,\ m,n \geq 0,
\end{equation}
where `coherent' means that the collection~(\ref{karta}) assembles into a
map $\alpha : \freeWPROP (\sfP) \to \sfP$ with the properties stated in
Definition~\ref{snad_projdu}.

As in the case of ordinary \PROP{s}, $\freeWPROP(E)$ carries the
`tautological' wheeled \PROP\ structure for an arbitrary
$\Sigma$-bimodule $E$. Another example of a wheeled \PROP\ is the
endomorphism \PROP\ $\End_V$ of a finite-dimensional vector space $V$
discussed in Example~\ref{fund}.

\begin{definition}
Let $\sfP$ be a wheeled \PROP\ and $V$ a finite-dimensional vector
space.  A {\em $\sfP$-algebra\/} (also called a {\em wheeled
representation\/} of $\sfP$) is a \PROP\ homomorphism $\rho : \sfP \to \End_V$,
where $\End_V$ is the wheeled endomorphism \PROP\ introduced in
Example~\ref{fund}.
\end{definition}

\begin{figure}[t]
\begin{center}
{
\unitlength=1.000000pt
\begin{picture}(190.00,80.00)(0.00,0.00)
\thicklines
\put(140.00,40.00){\vector(-1,1){0.00}}
\qbezier(140.00,40.00)(150.00,30.00)(160.00,30.00)
\qbezier(160.00,50.00)(170.00,50.00)(170.00,40.00)
\qbezier(160.00,30.00)(170.00,30.00)(170.00,40.00)
\qbezier(140.00,40.00)(150.00,50.00)(160.00,50.00)
\put(100.00,70.00){\makebox(0.00,0.00){$G/e$:}}
\put(180.00,40.00){\makebox(0.00,0.00){$f$}}
\put(60.00,40.00){\makebox(0.00,0.00){$f$}}
\put(140.00,40.00){\vector(0,1){40.00}}
\put(140.00,0.00){\vector(0,1){40.00}}
\put(0.00,70.00){\makebox(0.00,0.00){$G:$}}
\put(20.00,40.00){\makebox(0.00,0.00){$e$}}
\put(140.00,40.00){\makebox(0.00,0.00){$\bullet$}}
\put(40.00,60.00){\vector(-1,1){0.00}}
\put(40.00,60.00){\vector(1,1){0.00}}
\qbezier(50.00,40.00)(50.00,50.00)(40.00,60.00)
\qbezier(40.00,20.00)(50.00,30.00)(50.00,40.00)
\qbezier(30.00,40.00)(30.00,50.00)(40.00,60.00)
\qbezier(40.00,20.00)(30.00,30.00)(30.00,40.00)
\put(40.00,0.00){\vector(0,1){20.00}}
\put(40.00,60.00){\vector(0,1){20.00}}
\put(40.00,60.00){\makebox(0.00,0.00){$\bullet$}}
\put(40.00,20.00){\makebox(0.00,0.00){$\bullet$}}
\end{picture}}
\end{center}
\caption{\label{fig2}
A wheel created by collapsing an edge in an unwheeled graph.
}
\end{figure}
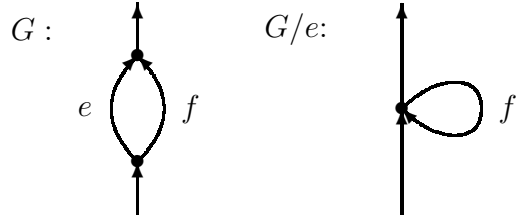

Surprisingly,  wheeled \PROP{s} are simpler
than (ordinary) \PROP{s}, because the category of all directed
graphs is, in contrast with the category of
graphs without wheels, closed under edge contractions.
This means that, given a directed graph $G \in \DGrw mn$ and an edge
$e$ of $G$, the quotient $G/e$ is again a directed graph, but $G/e$
may contain wheels although $G$ does not, see Figure~\ref{fig2}.
The consequence of this observation is that
the compositions~(\ref{karta}) in a wheeled \PROP\ are generated by
edge contractions and disjoint unions of graphs.

\begin{figure}[t]
\begin{center}
\thicklines
{
\unitlength=.7pt
\begin{picture}(120.00,140.00)(100.00,25.00)
\put(110,25){
\put(0,-22)
{
\put(45.00,107.00){\makebox(0.00,0.00)[b]{$\cdots$}}
\put(20.00,107.00){\makebox(0.00,0.00)[b]{$\cdot\, \cdot$}}
}
\put(0,-80)
{
\put(45.00,107.00){\makebox(0.00,0.00)[b]{$\cdots$}}
\put(20.00,107.00){\makebox(0.00,0.00)[b]{$\cdot\, \cdot$}}
}
\put(22.00,20.00){\makebox(0.00,0.00)[t]{$j$}}
\put(40.00,100.00){\makebox(0.00,0.00){$i$}}
\put(90.00,100.00){\line(0,-1){70.00}}
\qbezier(60.00,0.00)(90.00,0.00)(90.00,30.00)
\qbezier(30.00,30.00)(30.00,0.00)(60.00,0.00)
\qbezier(60.00,130.00)(90.00,130.00)(90.00,100.00)
\qbezier(30.00,90.00)(30.00,130.00)(60.00,130.00)
\put(30.00,60.00){\line(0,1){30.00}}
\put(30.00,30.00){\vector(0,1){30.00}}
\put(60.00,30.00){\vector(-1,1){30.00}}
\put(10.00,30.00){\vector(2,3){20.00}}
\put(0.00,30.00){\vector(1,1){30.00}}
\put(30.00,60.00){\vector(1,1){30.00}}
\put(30.00,60.00){\vector(-2,3){20.00}}
\put(30.00,60.00){\vector(-1,1){30.00}}
\put(30.00,60.00){\makebox(0.00,0.00){$\bullet$}}
}
\end{picture}}
\end{center}
\caption{\label{fig4}
The graph generating contractions~(\protect\ref{contr}).
Its vertex has biarity $(m,n)$.
}
\end{figure}
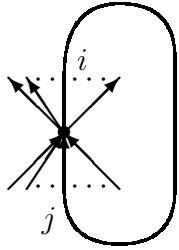

\begin{figure}[t]
\begin{center}
\thicklines
{
\unitlength=.7pt
\begin{picture}(120.00,170.00)(100.00,0.00)
\thicklines
\put(30.00,30.00){\makebox(0.00,0.00){$\bullet$}}
\put(30.00,140.00){\makebox(0.00,0.00){$\bullet$}}
\put(40.00,65.00){\makebox(0.00,0.00)[b]{$j$}}
\put(40.00,105.00){\makebox(0.00,0.00)[t]{$i$}}
\put(35.00,0.00){\makebox(0.00,0.00)[b]{$\cdots$}}
\put(35.00,170.00){\makebox(0.00,0.00)[t]{$\cdots$}}
\put(45.00,107.00){\makebox(0.00,0.00)[b]{$\cdots$}}
\put(20.00,107.00){\makebox(0.00,0.00)[b]{$\cdot\, \cdot$}}
\put(0,-45){
\put(45.00,107.00){\makebox(0.00,0.00)[t]{$\cdots$}}
\put(20.00,107.00){\makebox(0.00,0.00)[t]{$\cdot\, \cdot$}}
}
\put(30.00,30.00){\vector(1,1){30.00}}
\put(30.00,30.00){\vector(-2,3){20.00}}
\put(30.00,30.00){\vector(-1,1){30.00}}
\put(60.00,0.00){\vector(-1,1){30.00}}
\put(10.00,0.00){\vector(2,3){20.00}}
\put(0.00,0.00){\vector(1,1){30.00}}
\put(30.00,30.00){\vector(0,1){110.00}}
\put(60.00,110.00){\vector(-1,1){30.00}}
\put(10.00,110.00){\vector(2,3){20.00}}
\put(0.00,110.00){\vector(1,1){30.00}}
\put(30.00,140.00){\vector(1,1){30.00}}
\put(30.00,140.00){\vector(-2,3){20.00}}
\put(30.00,140.00){\vector(-1,1){30.00}}
\put(200,40){
\put(0.00,30.00){\vector(0,1){10.00}}
\put(0,-10){
\qbezier(130.00,20.00)(130.00,0.00)(110.00,0.00)
\qbezier(0.00,20.00)(0.00,0.00)(20.00,0.00)
}
\put(0,10){
\qbezier(110.00,80.00)(130.00,80.00)(130.00,60.00)
\qbezier(70.00,60.00)(70.00,80.00)(90.00,80.00)
}
\put(110.00,-10.00){\line(-1,0){90.00}}
\put(130.00,70.00){\line(0,-1){60.00}}
\put(90.00,90.00){\line(1,0){20.00}}
\put(70.00,40.00){\line(0,1){30.00}}
\put(0.00,40.00){\line(0,-1){30.00}}
\put(-90.00,40.00){\makebox(0.00,0.00){$=$}}
\put(40,10)
{
\put(0,-45){
\put(45.00,107.00){\makebox(0.00,0.00)[t]{$\cdots$}}
\put(20.00,107.00){\makebox(0.00,0.00)[t]{$\cdot\, \cdot$}}
}
\put(30.00,30.00){\makebox(0.00,0.00){$\bullet$}}
\put(35.00,0.00){\makebox(0.00,0.00)[b]{$\cdots$}}
\put(30.00,30.00){\vector(1,1){30.00}}
\put(30.00,30.00){\vector(-2,3){20.00}}
\put(30.00,30.00){\vector(-1,1){30.00}}
\put(60.00,0.00){\vector(-1,1){30.00}}
\put(10.00,0.00){\vector(2,3){20.00}}
\put(0.00,0.00){\vector(1,1){30.00}}
\put(16,88){\makebox(0.00,0.00){$m_1 + j$}}
}
\put(-30,10)
{
\put(0,-45){
\put(45.00,52.5){\makebox(0.00,0.00)[t]{$\cdots$}}
\put(20.00,52.5){\makebox(0.00,0.00)[t]{$\cdot\, \cdot$}}
}
\put(30.00,30.00){\makebox(0.00,0.00){$\bullet$}}
\put(35.00,53.00){\makebox(0.00,0.00)[b]{$\cdots$}}
\put(30.00,30.00){\vector(1,1){30.00}}
\put(30.00,30.00){\vector(-2,3){20.00}}
\put(30.00,30.00){\vector(-1,1){30.00}}
\put(60.00,0.00){\vector(-1,1){30.00}}
\put(10.00,0.00){\vector(2,3){20.00}}
\put(0.00,0.00){\vector(1,1){30.00}}
\put(21.00,-10.00){\makebox(0.00,0.00){$i$}}
}
}
\end{picture}
}
\end{center}
\caption{\label{fig4bis}
A dioperadic composition as a horizontal compositions followed by a
contraction.
The upper vertex of the left graph has biarity $(m_1,n_1)$, the bottom
vertex biarity $(m_2,n_2)$.
}
\end{figure}
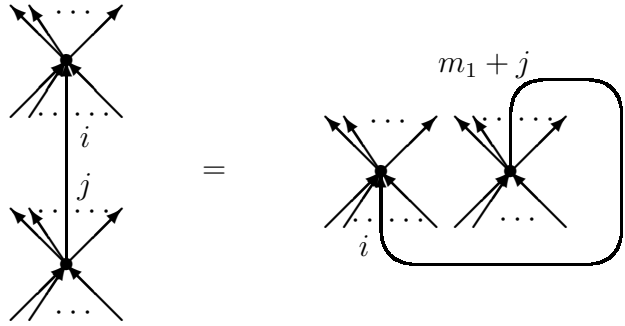

Therefore a biased definition of wheeled \PROP{s} can be given in
terms of the horizontal compositions~(\ref{hor}) that correspond to
disjoint unions of graphs, and the contractions
\begin{equation}
\label{contr}
\xi_j^i : \sfP(m,n) \to \sfP(m-1,n-1)
\end{equation}
defined for $1 \leq i \leq m$, $1 \leq j \leq n$,
that correspond to the graph in Figure~\ref{fig4}.
The `dioperadic' (see~\cite{gan} for the
terminology) compositions
\begin{equation}
\label{kocicka_Lizinka}
\circ^i_j : \sfP(m_1,n_1)  \ot  \sfP(m_2,n_2) \to  \sfP(m_1 +
m_2-1,n_1+n_2-1),
\end{equation}
where $m_1,n_2 \geq 0$, $1 \leq i \leq n_1$, $1 \leq j \leq
m_2$, corresponding to the contraction of an edge joining two
different vertices, can be expressed in terms of the above two
operations as
\[
p \circ^i_j q = \xi^{m_1 +j}_i (p \ot q),\ p \in \sfP(m_1,n_1),\
q \in \sfP(m_2,n_2),
\]
see Figure~\ref{fig4bis}.
We are, however, not going to write axioms that these operations should fulfill
here.

The obvious forgetful functor $\Box : {\tt PROP}^\circlearrowright \to {\tt PROP}$
from the category of wheeled \PROP{s} to the category of \PROP{s} is
induced by the natural inclusions $\DGr mn \hookrightarrow \DGrw
mn$, $m,n \geq 0$. It is easy to show that $\Box$ has a left adjoint
$(-)^\circlearrowright : {\tt PROP} \to {\tt PROP}^\circlearrowright$
\begin{definition}
\label{Jitka}
Given a \PROP\/ $\sfP$,
we call the wheeled \PROP\/ $\sfP^\circlearrowright$ the
{\em wheeled completion\/} of\/~$\sfP$.
\end{definition}

\begin{proposition}
\label{Winnona} For any PROP $\sfP$
there is a one-to-one correspondence between finite dimensional
representations of $\sfP$ in the category $\tt PROP$ and finite
dimensional representations of $\sfP^\circlearrowright$ in the
category of\, $\tt PROP^\circlearrowright$.
\end{proposition}

\begin{definition}
\label{suspension}

For a graded vector space $V$, let
$\susps V$ be the {\em suspension of
V\/}  defined by $(\susp V)_i := V_{i-1}$,
and let
$\desusps V$ be the {\em desuspension of
V\/}  defined by $(\susp V)_i := V_{i+1}$.

\noindent
The {\em suspension}\, wheeled properad, $\mathsf S$, is the endomorphism wheeled prop(erad) of
$\susps{\mathbf k}$ and the {\em desuspension}\, wheeled properad, $\mathsf S^{-1}$, is the endomorphism
wheeled prop(erad) of
$\desusps{\mathbf k}$. Thus $\mathsf S$ is the one dimensional $\Sigma$-bimodule
$\{\susps^{m-n}\sgn_m\ot\sgn_n\}$ and  $\mathsf S^{-1}$ is the one dimensional $\Sigma$-bimodule
$\{\desusps^{m-n}\sgn_m\ot\sgn_n\}$. Tensoring with $\mathsf S$ (or with $\mathsf S^{-1}$) defines
an endomorphism functor in the category of wheeled prop(erad)s.
\end{definition}

\subsection{Modifications and generalizations}
In~\cite{vallette:thesis,vallette:CMR04}, B.~Vallette
introduced {\em properads\/} as a suitable subcategory of the category
of (ordinary)
\PROP{s} on which (co)bar constructions and quadratic duality could be
defined. We will need, in Section~\ref{pekelna_hudba}, a wheeled
version of this notion. Recall that Vallette's
properads are algebras
over a subtriple $\rmFname$
of the free \PROP\ triple $\Gamma$ given as the colimit over {\em
connected\/} graphs, that is, for a $\Sigma$-bimodule $E$, $\rmF(E)$ is
defined by
\[
\rmF(E)(m,n) := \colim{{G \in \DGrc mn}}{{E}}(G),\
m,n \geq 0,
\]
where $\DGrc mn$ is the full subcategory of $\DGr mn$ consisting of
connected graphs. Similarly,
there is a subtriple $\WrmFname$ of the triple $\Gamma^\circlearrowright$
defined by
\begin{equation}
\label{l}
\WrmF(E)(m,n) := \colim{{G \in \DGrwc mn}}{{E}}(G),\
m,n \geq 0,
\end{equation}
where $\DGrwc mn$ is the full subcategory of $\DGrw mn$ of
connected graphs. Observe that there are precisely two connected exceptional
graphs,
\[
\uparrow\ \in \DGrwc 11\ \mbox { and } \ \hgw \in \DGrwc 00.
\]

\begin{definition}
\label{wheeled_properad}
A {\em wheeled properad\/} is an algebra over the triple  $\WrmFname$
introduced above.
\end{definition}

It is clear that, for each $\Sigma$-bimodule,
$\WrmF(E)$ is a wheeled properad, with the structure map
given by the triple multiplication $\mu_E : \WrmF({\WrmF(E)}) \to
\WrmF(E)$. It is the free wheeled properad generated by the
$\Sigma$-bimodule $E$.

Observe that each wheeled properad generates a wheeled \PROP, but not
all wheeled \PROP{s} are generated by wheeled properads. Each properad $\sfP$
has its {\em wheeled properadic completion\/} $\sfP_c^\wc$ given by an obvious
modification of Definition~\ref{Jitka}.

\begin{example}
\label{augment}
{\rm
Let $\frakk$ be the trivial properad (the initial object of the category
properads) defined~by
\[
\frakk(m,n) := \cases{\bfk}{for $(m,n) = (1,1)$, and}{0}{otherwise.}
\]
Its wheeled properadic completion $\frakk_c^\wc$ is the initial object of the
category of wheeled properads. It satisfies
\[
\frakk_c^\wc(m,n) := \cases{\bfk}{for $(m,n) = (1,1)$ or $(0,0)$, and}
{0}{otherwise.}
\]
The component $\frakk_c^\wc(0,0)$ is spanned by the contraction of the
unit $1 \in \frakk_c^\wc(1,1)$.
Of course, $\frakk_c^\wc$ equals $\WrmF(0)$, the free wheeled properad
generated by the trivial $\Sigma$-bimodule.
}\end{example}

We close this section by mentioning
an important modification of wheeled properads whose nature
resembles modular operads introduced
in~\cite{getzler-kapranov:CompM98}. Let us recall some necessary definitions.

A {\em labeled graph\/} is a connected graph $G$ together with
a map $g$ from $\vert(G)$ to
the set $\{0,1,2,\ldots\}$. The {\em genus\/}
$g(G)$
of a labeled graph $G$ is defined by
\[
g(G) :=  \dim H_1(G) +\sum_{v\in \vert(G)} g(v).
\]
Let us denote by $\Magogh gmn$ the category of
labeled directed $(m,n)$-graphs and their isomorphisms.

By
a {\em modular $\Sigma$-bimodule\/} we mean a
system $E= \{ E(g;m,n)\}_{g,m,n \geq 0}$ such that each
$E(g;m,n)$ is a left $\bfk[\Sigma_m]$- right $\bfk[\Sigma_n]$-bimodule.
For a modular $\Sigma$-bimodule
$E$ and $G \in \Magogh gmn$, let
\[
E(G) := \bigotimes_{v \in \vert(G)} E(g(v);\Out(v),\In(v)).
\]
Define finally
\[
\freeMproperad(E)(g;m,n) := \colim{{G \in \Magogh gmn}}{{E}}(G),\
g,m,n \geq 0.
\]
As before, $\freeMproperadname$ is a triple in
the category of modular $\Sigma$-bicollections.

\begin{definition}
\label{Kocicka_Micinka}
Modular wheeled properads are algebras over the triple $\freeMproperadname$.
\end{definition}

Loosely speaking, modular wheeled properads are wheeled properads
equipped with a genus grading such that the dioperadic
compositions~(\ref{kocicka_Lizinka})
preserve the genus and the contractions~(\ref{contr})
increase the genus by one. One may also say that modular wheeled
properads are {\em directed (unstable) modular operads\/}.

\begin{example}
\label{hlavne_abych_byl_zdravy}
For each  $\Sigma$-bimodule $E$, the free wheeled properad $\WrmF(E)$
is modular, with the genus grading induced by the genus of underlying
graphs. Obviously, each wheeled properad~$\sfP$ of the form $\sfP =
\WrmF(E)/I$, with the ideal $I$ generated by elements supported by
genus zero graphs, has the induced modular structure.
This means that the wheeled properadic completion of
a~dioperad~\cite{gan} is a modular wheeled properad.
\end{example}

\begin{example}
\label{i1}
If one allows the genus to be an arbitrary integer, then each wheeled
properad $\sfP$ can be turned into a modular one by assigning
$\sfP(m,n)$ the genus $1-m$.
\end{example}


\section{Master equations}
\label{Master}

\subsection{Differential $\Z_2$-graded Lie algebras.}
In the context of Batalin-Vilkovisky
formalism in quantum field theory it
is more suitable to work with the {\em odd}\, version of the usual
notion of  differential   Lie superalgebra. By definition, this is a
$\Z_2$-graded vector space, $\fg=\fg_{\tlo}\oplus \fg_{\tln}$,
equipped with two odd linear maps
$$
d: \fg \rar \fg, \ \ \ \ \ \mbox{and} \ \ \ \ \ [\ \bullet \ ]:
\fg\ot \fg \rar \fg,
$$
such that $d^2=0$, $[a\bullet b]=-(-1)^{(\tl{a}+1)(\tl{b}+1)}[b \bullet a]$,
$$
d[a\bullet b]= [da\bullet b] + (-1)^{(\tl{a}+1)}[a\bullet db],
$$
 and
 $$
           [a\bullet [b\bullet c]]= [[a\bullet b]\bullet c] +
           (-1)^{(\tl{a}+1)(\tl{b}+1)} [b\bullet [a\bullet c],
$$
for all $a,b,c\in \fg_{\tlo}\cup \fg_{\tln}$.

In many important examples, the $\Z_2$-grading in $\fg$ comes from
an underlying $\Z$-grading, i.e. $\fg=\oplus_{i\in \Z}\fg^i$,
$\fg_{\tlo}= \oplus_{i\, \mathrm even}\fg^i$, $\fg_{\tln}=
\oplus_{i\, \mathrm odd}\fg^i$, and the basic operations satisfy
$d\fg^i\subset \fg^{i+1}$, $[\fg^i\bullet \fg^j]\subset
\fg^{i+j-1}$.

Clearly, the parity change functor transforms this structure into
the ordinary
structure of differential Lie superalgebra on the vector
superspace $\Pi \fg$.

\subsection{Differential Gerstenhaber-Batalin-Vilkovisky algebras.}
Such an algebra is a quadruple
$(\fg, \circ, d, \Delta)$, where $(\fg,
\circ)$ is a unital supercommutative algebra over a field $k$, and
$(d,\Delta)$ is a pair of supercommuting  odd derivations of
$(\fg,\circ)$ of orders $1$ and, respectively, $2$  which satisfy
$d^2=\Delta^2=0$.

More explicitly, a dGBV algebra is  a differential supercommutative
algebra with unit, $(\fg,\circ, d)$, equipped an odd linear map
$\Delta: \fg
\rar \fg$ satisfying
\Bi
\item[(i)] $\Delta^2=0$,
\item[(ii)] $d\Delta + \Delta d =0$, and
\item[(iii)] for any $a,b,c\in \fg$,
\Beqrn \Delta(a\circ b \circ c) & = & \Delta(a\circ b) \circ c +
(-1)^{\tl{b}(\tl{a}+1)} b\circ \Delta(a\circ c) +
(-1)^{\tl{a}}a\circ \Delta(b\circ c)\\
&& - \Delta(a) \circ b \circ c - (-1)^{\tl{a}} a\circ \Delta(b)\circ
c - (-1)^{(\tl{a}+\tl{b})} a\circ b \circ \Delta(c).
\Eeqrn
\Ei

Note that $\Delta(1)=0$.

It is not hard to check using  identity (iii) that the linear map
$$
\Ba{rccl}
[\ \bullet \ ] & : \fg\ot \fg & \lon & \fg \\
& \hskip .6em a \ot b & \lon & [a\bullet b]:= (-1)^{\tl{a}} \Delta(a\circ b) -
(-1)^{\tl{a}} \Delta(a)\circ b - a\circ \Delta(b) \Ea
$$
makes $\fg$ into an odd Lie superalgebra.

Moreover, both triples
$(\fg,
[\ \bullet\ ], d)$ and $(\fg, [\ \bullet \ ], \Delta)$ are odd
differential Lie superalgebras, and the following odd Poisson
identity,
$$
[a\bullet (b\circ c)] = [a\bullet b]\circ c +
(-1)^{\tl{a}(\tl{b}+1)}b\circ [a\bullet b],
$$
holds for any $a,b,c\in \fg$.

\subsection{Master equation.}
An even element $S$ in a dGBV algebra
$\fg$ is called a {\em master function}\, if it satisfies
the {\em master equation}
$$
dS + \Delta S + \frac{1}{2}[S\bullet S]=0.
$$

Assume that $\fg$ is such that
the formal power series,
$$
e^S= 1 + S + \frac{1}{2!}S\circ S + \frac{1}{3!}S\circ S\circ S +\ldots,
$$
makes sense, i.e.\ gives a well-defined element of $\fg$ (often this is
achieved by introducing a formal parameter $\hbar$ and working
in $\fg[[\hbar]]$). One of the central observation in the theory
of master equations
is the following

\begin{lemma}
$S$ is a master function if and only if
$$
(d+\Delta)e^S=0.
$$
\end{lemma}

\begin{proof}
One checks that $(d+\Delta)e^S= (dS + \Delta S + \frac{1}{2}[S\bullet S])\circ
e^S$.
\end{proof}

\subsection{Master equations in geometry}
\label{MA}
Let $M$ be a smooth $n$-dimensional
manifold with the tangent sheaf denoted by $T_M$ and the
sheaf of differential forms denoted by $\Omega^\bullet_M=\bigoplus_{i=0}^n
\Omega^i_M$. It is well-known that the sheaf
of polyvector  fields,
$$
\wedge^\bullet T_M= \bigoplus_{i=0}^n \wedge^i T_M
$$
is a sheaf of supercommutative algebras with respect to the wedge
product,
$\circ=\wedge$, and also a sheaf of odd Lie superalgebras
with respect to the Schouten bracket,
$$
[\ \bullet \ ]_{Schouten}:  \wedge^iT_M \ot \wedge^jT_M
\lon \wedge^{i+j-1}T_M.
$$
Moreover, the odd Poisson identity,
$$
[a\bullet (b\wedge c)]_{Schouten} = [a\bullet b]_{Schouten}\wedge c +
(-1)^{\tl{a}(\tl{b}+1)}b\wedge [a\bullet b]_{Schouten},
$$
holds for any $a,b,c\in \wedge^\bullet T_M$.

Assume now that $M$ is
 equipped with a volume form $\nu$, that is, with a nowhere
vanishing section of $\Omega^n_M$. In particular, the associated section
$\nu^{-1}$ of $\wedge^n T_M$ is  well defined.

Define  a differential operator $\Delta: \wedge^\bullet T_M\rar
\wedge^\bullet T_M $ as the composition
$$
\Delta: \wedge^i T_M \stackrel{\lrcorner \nu}{\lon}
\Omega^{n-i}_M \stackrel{d}{\lon} \Omega^{n-i+1}_M
\stackrel{\lrcorner \nu^{-1}}{\lon} \wedge^{i-1} T_M
$$
where $d$ stands for de~Rham
differential and $\lrcorner$ for the
natural contraction of dual tensors. It is well-known (and easy to check,
see a local coordinate description below) that the data
$(\wedge^\bullet T_M, d=0, \Delta)$ is dGBV algebra with the associated
odd Lie algebra structure being exactly the Schouten structure.
The associated master equation has the form
$$
\Delta S + \frac{1}{2}[S\bullet S]_{Schouten}=0.
$$

Sections of the bundle $\wedge^\bullet T_M$ can be understood as
functions on the supermanifold $\Pi \Omega^1_M$, the total space
of the cotangent bundle with the parity of fiber
coordinates changed.
If $\{ x^\al \}_{1\leq \al\leq n}$ are local
coordinates on $M$, then the functions
 $\{ x^\al, \psi_\al:=\Pi\p/\p x^\al \}_{1\leq \al \leq n}$
form a local coordinate system on $\Pi \Omega^1_M$. The volume
form $\nu$ gets the coordinate representation
$\nu=f dx^1\wedge\ldots \wedge dx^n$ for some non-vanishing function $f$, while the divergence operator $\Delta$
gets explicitly represented as
$$
\Delta=\sum_{\al=1}^n \frac{\p^2}{\p x^\al \p \psi_\al} +
\frac{\p \ln f}{\p x^\al}\frac{\p}{\p \psi_\al}.
$$
Hence one has
are
$$
[f\bullet g]_{Schouten}= \sum_{\al=1}^n\left(
\frac{{\p} f}{\p x^\al} \frac{{\p} g}{\p \psi_\al}
+ (-1)^{\tl{g}\tl{f}}
\frac{{\p} f}{\p x^\al} \frac{{\p} g}{\p \psi_\al}
\right).
$$
for arbitrary $f,g\in \wedge^\bullet T_M$.

\begin{CY}
{\rm
There is a  variant of the above master equation for CY manifolds
with $\nu$ being the holomorphic volume form. If one defines
$\Delta_{\p}$ on holomorphic vector fields,
 $\wedge^\bullet \cT_M$, precisely as above
with the full de~Rham differential $d=\p+ \bar{\p}$ replaced
by its $(1,0)$-part
$\p$, then the sheaf
$\wedge^\bullet \cT_M\ot \Omega^{0,\bullet}_M$ is
a dGBV algebra with differential $d:=\Id\ot\bar{\p}$
and the 2nd order derivation $\Delta:=\Delta_{\p}\ot \Id$.
The associated master equation
$$
dS + \hbar\Delta S + \frac{1}{2}[S\bullet S]=0,\ \ S\in
\wedge^\bullet \cT_M\ot \Omega^{0,\bullet}_M[[\hbar]]
$$
plays a key role in the Barannikov-Kontsevich approach \cite{B,BK} to
the $B$-model side of the Mirror Symmetry. It solutions describe extended
deformations of the complex structure on a Calabi-Yau
manifold.
}
\end{CY}

\begin{MES}
In theoretical physics, one is more interested in a version of the
construction in \S\ref{MA} when the underlying space $M$ is a
supermanifold rather than a manifold. In both cases the supermanifold
$\cM:=\Pi \Omega^1 M$ has an odd symplectic structure $\omega\in \Omega^2\cM$,
 but the notions of volume forms are different --- in the supermanifold case it should be
understood as a section of the Berezinian bundle,  rather than a
differential form. Given such a section, $\nu \in \mbox{Ber}({\cM})$,
one obtains, for
an arbitrary
vector field $\zeta\in \cT_\cM$, another
section, $Lie_\zeta \nu \in \mbox{Ber}({\cM})$. As $\nu$ is a basis section,
one can write,
$$
Lie_\zeta \nu = (\mbox{div} \zeta) \nu,
$$
for some well-defined function $\mbox{div} \zeta\in \f_\cM$ called the
{\em divergence}\, of $\zeta$. As the odd 2-form $\omega$ is non-degenerate,
for any function $f\in \f_\cM$ there exists an associated Hamiltonian vector
field, $H_f\in \cT_\cM$, uniquely defined by the equation,
$$
df=\omega\lrcorner H_f.
$$
Then one can define an odd differential operator $\Delta$ on
$\f_\cM\simeq \wedge^\bullet T_M$ as follows,
$$
\Delta f:= \mbox{div}H_f,
$$
and check that $\Delta^2=0$. If there exist Darboux coordinates,
$\{ x^\al, \psi_\al \}_{1\leq \al\leq n}$, $|{\psi_\al}|
=|{x^\al}| + 1\mod 2\Z$, on $\cM$ such that the Berezinian
$\nu$ is locally constant, then the derivation $\Delta$ gets a
very simple form,
\begin{equation}\label{lapl}
\Delta=\sum_{\al=1}^n \frac{\p^2}{\p x^\al \p \psi_\al}.
\end{equation}
In this case the data $(\cM, \omega, \nu)$ is called an
 $SP$-{\em manifold} \cite{Sch},
and provides us with the most general
and down-to-earth mathematical description of key structures in
the Batalin-Vilkovisky quantization.

A formal $SP$-manifold can be canonically associated with an arbitrary finite
dimensional complex $(M_0, d_0)$ over a field $k$.
Let $\{e_\al\}$ be a homogeneous basis
of $M_0$ and let $\{x^\al\}$ be the dual basis of $M_0^*$.
The differential
 $d_0$ can be understood as a linear
odd vector field, $d_0\in T_M$, on the formal manifold $M$
associated with the vector superspace $M_0$. Indeed, if
$$
{d}_0(e_\al)= \sum_{\al,\be}{D^\be_\al}e_{\be}, \ \ \ \ \ {D^\be_\al}\in k,
$$
then the associated vector field is given by
$$
\vec{d}_0=\sum_{\al\be} (-1)^{\al}x^{\al}D_{\al}^{\be}\frac{\p}{\p x^\al}.
$$
Clearly, the equation $d_0^2=0$ is equivalent to $[\vec{d}_0, \vec{d}_0]_{Schouten}=0$.
Hence the $\f_M$-module
$\wedge^\bullet T_M$, $\f_M:=\widehat{\odot} M_0^*$ being the ring of formal
smooth functions, is naturally a {\em differential}\,
$\f_M$-module with the differential given by
$$
d f:= [d_0\bullet f]_{Schouten}
$$
for any $f\in \wedge^\bullet T_M$. The sheaf $\mbox{Ber}(\Pi \Omega^1_M)$ has a
distinguished constant section $\nu$ such that the associated odd Laplacian
$\Delta$ has the form (\ref{lapl}) with $\psi_\al=\Pi\p/\p x^\al$.
The Schouten brackets on $\f_{\Pi \Omega^1_M}=\wedge^\bullet T_M$ get the form,
$$
[f\bullet g]= \sum_{\al=1}^n\left(
\frac{\overleftarrow{\p} f}{\p x^\al} \frac{\overrightarrow{\p} g}{\p \psi_\al}
+ (-1)^{\tl{g}\tl{f}}
\frac{\overleftarrow{\p} f}{\p x^\al} \frac{\overrightarrow{\p} g}{\p \psi_\al}
\right).
$$
The resulting data $(\wedge^\bullet T_M, \Delta, d)$ is a dGBV algebra
canonically associated with the complex $(M_0,d_0)$. Hence it makes
sense to associate with the latter the following Master equation,
\begin{equation}\label{meq}
dS + \Delta S + \frac{1}{2}[S\bullet S]=0.
\end{equation}
\end{MES}

\begin{theorem}
There exists a differential wheeled PROP
$(\mbox{\sf PolyV}^\circlearrowright, \delta)$ whose
wheeled representations in a finite
dimensional complex $(M_0, d_0)$ are in one-to-one correspondence
with Master functions $S$ in the dGBV algebra
$(\wedge^\bullet T_M,  \Delta, d)$.
\end{theorem}

\sip

\begin{proof}
 Let $E=\{E(m,n):=\sgn_m\ot \id_n\}_{m,n\geq 0}$
be a collection of $\Z_2$-graded $(\Sigma_m, \Sigma_n)$ bimodules
concentrated in degree $m\mod 2\Z$.
Here $\sgn_m$ stands for the sign representation of
$\Sigma_m$
and $\id_n$ for the trivial representation of $\Sigma_n$. Every $E(m,n)$
is therefore a one-dimensional space whose basis vector we denote
graphically as a planar $(m,n)$-corolla,
$$
 
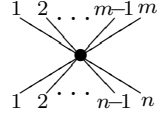
\begin{xy}
 <0mm,0mm>*{\bullet};<0mm,0mm>*{}**@{},
 <0mm,0mm>*{};<-8mm,5mm>*{}**@{-},
 <0mm,0mm>*{};<-4.5mm,5mm>*{}**@{-},
 <0mm,0mm>*{};<-1mm,5mm>*{\ldots}**@{},
 <0mm,0mm>*{};<4.5mm,5mm>*{}**@{-},
 <0mm,0mm>*{};<8mm,5mm>*{}**@{-},
   <0mm,0mm>*{};<-8.5mm,5.5mm>*{^1}**@{},
   <0mm,0mm>*{};<-5mm,5.5mm>*{^2}**@{},
   <0mm,0mm>*{};<4.5mm,5.5mm>*{^{m\hspace{-0.5mm}-\hspace{-0.5mm}1}}**@{},
   <0mm,0mm>*{};<9.0mm,5.5mm>*{^m}**@{},
 <0mm,0mm>*{};<-8mm,-5mm>*{}**@{-},
 <0mm,0mm>*{};<-4.5mm,-5mm>*{}**@{-},
 <0mm,0mm>*{};<-1mm,-5mm>*{\ldots}**@{},
 <0mm,0mm>*{};<4.5mm,-5mm>*{}**@{-},
 <0mm,0mm>*{};<8mm,-5mm>*{}**@{-},
   <0mm,0mm>*{};<-8.5mm,-6.9mm>*{^1}**@{},
   <0mm,0mm>*{};<-5mm,-6.9mm>*{^2}**@{},
   <0mm,0mm>*{};<4.5mm,-6.9mm>*{^{n\hspace{-0.5mm}-\hspace{-0.5mm}1}}**@{},
   <0mm,0mm>*{};<9.0mm,-6.9mm>*{^n}**@{},
 \end{xy}
$$
with $m$ skew-symmetric outgoing legs and $n$ symmetric ingoing legs.

\sip

Let
$\mbox{\sf PolyV}$ be the free PROP
generated by the bimodule $E$, and let
$\mbox{\sf PolyV}^\circlearrowright$ be the free wheeled PROP
generated by the same bimodule $E$.
Any derivation,
$\delta: \mbox{\sf PolyV}^\circlearrowright\rar
\mbox{\sf PolyV}^\circlearrowright$, is completely determined by its
values on the above corollas.

\bip

{\em Claim.} The derivation, $\delta$, defined by
$$
\delta \begin{xy}
 <0mm,0mm>*{\bullet};<0mm,0mm>*{}**@{},
 <0mm,0mm>*{};<-8mm,5mm>*{}**@{-},
 <0mm,0mm>*{};<-4.5mm,5mm>*{}**@{-},
 <0mm,0mm>*{};<-1mm,5mm>*{\ldots}**@{},
 <0mm,0mm>*{};<4.5mm,5mm>*{}**@{-},
 <0mm,0mm>*{};<8mm,5mm>*{}**@{-},
   <0mm,0mm>*{};<-8.5mm,5.5mm>*{^1}**@{},
   <0mm,0mm>*{};<-5mm,5.5mm>*{^2}**@{},
   <0mm,0mm>*{};<4.5mm,5.5mm>*{^{m\hspace{-0.5mm}-\hspace{-0.5mm}1}}**@{},
   <0mm,0mm>*{};<9.0mm,5.5mm>*{^m}**@{},
 <0mm,0mm>*{};<-8mm,-5mm>*{}**@{-},
 <0mm,0mm>*{};<-4.5mm,-5mm>*{}**@{-},
 <0mm,0mm>*{};<-1mm,-5mm>*{\ldots}**@{},
 <0mm,0mm>*{};<4.5mm,-5mm>*{}**@{-},
 <0mm,0mm>*{};<8mm,-5mm>*{}**@{-},
   <0mm,0mm>*{};<-8.5mm,-6.9mm>*{^1}**@{},
   <0mm,0mm>*{};<-5mm,-6.9mm>*{^2}**@{},
   <0mm,0mm>*{};<4.5mm,-6.9mm>*{^{n\hspace{-0.5mm}-\hspace{-0.5mm}1}}**@{},
   <0mm,0mm>*{};<9.0mm,-6.9mm>*{^n}**@{},
 \end{xy}
\ \ = \ \
 \sum_{I_1\sqcup I_2=(1,\ldots,m)\atop {J_1\sqcup J_2=(1,\ldots,n)\atop
}
}\hspace{0mm}
(-1)^{\sigma(I_1\sqcup I_2) + |I_1||I_2|}
 \begin{xy}
 <0mm,0mm>*{\bullet};<0mm,0mm>*{}**@{},
 <0mm,0mm>*{};<-8mm,5mm>*{}**@{-},
 <0mm,0mm>*{};<-4.5mm,5mm>*{}**@{-},
 <0mm,0mm>*{};<0mm,5mm>*{\ldots}**@{},
 <0mm,0mm>*{};<4.5mm,5mm>*{}**@{-},
 <0mm,0mm>*{};<13mm,5mm>*{}**@{-},
     <0mm,0mm>*{};<-2mm,7mm>*{\overbrace{\ \ \ \ \ \ \ \ \ \ \ \ }}**@{},
     <0mm,0mm>*{};<-2mm,9mm>*{^{I_1}}**@{},
 <0mm,0mm>*{};<-8mm,-5mm>*{}**@{-},
 <0mm,0mm>*{};<-4.5mm,-5mm>*{}**@{-},
 <0mm,0mm>*{};<-1mm,-5mm>*{\ldots}**@{},
 <0mm,0mm>*{};<4.5mm,-5mm>*{}**@{-},
 <0mm,0mm>*{};<8mm,-5mm>*{}**@{-},
      <0mm,0mm>*{};<0mm,-7mm>*{\underbrace{\ \ \ \ \ \ \ \ \ \ \ \ \ \ \
      }}**@{},
      <0mm,0mm>*{};<0mm,-10.6mm>*{_{J_1}}**@{},
 <13mm,5mm>*{};<13mm,5mm>*{\bullet}**@{},
 <13mm,5mm>*{};<5mm,10mm>*{}**@{-},
 <13mm,5mm>*{};<8.5mm,10mm>*{}**@{-},
 <13mm,5mm>*{};<13mm,10mm>*{\ldots}**@{},
 <13mm,5mm>*{};<16.5mm,10mm>*{}**@{-},
 <13mm,5mm>*{};<20mm,10mm>*{}**@{-},
      <13mm,5mm>*{};<13mm,12mm>*{\overbrace{\ \ \ \ \ \ \ \ \ \ \ \ \ \ }}**@{},
      <13mm,5mm>*{};<13mm,14mm>*{^{I_2}}**@{},
 <13mm,5mm>*{};<8mm,0mm>*{}**@{-},
 <13mm,5mm>*{};<12mm,0mm>*{\ldots}**@{},
 <13mm,5mm>*{};<16.5mm,0mm>*{}**@{-},
 <13mm,5mm>*{};<20mm,0mm>*{}**@{-},
     <13mm,5mm>*{};<14.3mm,-2mm>*{\underbrace{\ \ \ \ \ \ \ \ \ \ \ }}**@{},
     <13mm,5mm>*{};<14.3mm,-4.5mm>*{_{J_2}}**@{},
 \end{xy}
\ \ \ + \ \ \ (-1)^m
\begin{xy}
 <0mm,0mm>*{\bullet};<0mm,0mm>*{}**@{},
 <0mm,0mm>*{};<-8mm,5mm>*{}**@{-},
 <0mm,0mm>*{};<-4.5mm,5mm>*{}**@{-},
 <0mm,0mm>*{};<-1mm,5mm>*{\ldots}**@{},
 <0mm,0mm>*{};<4.5mm,5mm>*{}**@{-},
 <0mm,0mm>*{};<8mm,5mm>*{}**@{-},
   <0mm,0mm>*{};<-8.5mm,5.5mm>*{^1}**@{},
   <0mm,0mm>*{};<-5mm,5.5mm>*{^2}**@{},
   <0mm,0mm>*{};<4.5mm,5.5mm>*{^{m}}**@{},
 <0mm,0mm>*{};<-8mm,-5mm>*{}**@{-},
 <0mm,0mm>*{};<-4.5mm,-5mm>*{}**@{-},
 <0mm,0mm>*{};<-1mm,-5mm>*{\ldots}**@{},
 <0mm,0mm>*{};<4.5mm,-5mm>*{}**@{-},
 <0mm,0mm>*{};<8mm,-5mm>*{}**@{-},
   <0mm,0mm>*{};<-8.5mm,-6.9mm>*{^1}**@{},
   <0mm,0mm>*{};<-5mm,-6.9mm>*{^2}**@{},
   <0mm,0mm>*{};<4.5mm,-6.9mm>*{^{n}}**@{},
%
(7.9,5.0)*{}
   \ar@{->}@(ur,dr) (7.9,-5.0)*{}
 \end{xy}
$$
where $\sigma(I_1\sqcup I_2)$ is the sign of the shuffle
$I_1\sqcup I_2=(1,\ldots, m)$, is a differential, i.e.\ $\delta^2=0$.

{\em Proof}\, is a straightforward but tedious calculation.

\sip

For example
\[
\begin{array}{rclrcl}
\delta\,
\begin{xy}
 <0mm,0mm>*{\bullet};<0mm,0mm>*{}**@{},
\end{xy}
&=&
\begin{xy}
 <0mm,2.1mm>*{\bullet};<0mm,5mm>*{}**@{},
<0mm,-2.1mm>*{\bullet};<0mm,-5mm>*{}**@{},
<0mm,-2.5mm>*{};<0mm,2.5mm>*{}**@{-}
\end{xy}
\ \ + \ \
\begin{xy}
 <0mm,0mm>*{\bullet};<0mm,0mm>*{}**@{},
<0mm,-2.5mm>*{};<0mm,2.5mm>*{}**@{-},
(0,2.3)*{}
   \ar@{->}@(ur,dr) (0,-2.3)*{}
\end{xy},
&
\delta \hskip .15em
\begin{xy}
 <0mm,-1mm>*{\bullet};<0mm,0mm>*{}**@{},
<0mm,-1mm>*{\bullet};<0mm,3mm>*{}**@{-},
\end{xy}
&=&
\begin{xy}
 <0mm,2.1mm>*{\bullet};<0mm,0mm>*{}**@{},
<0mm,-2.1mm>*{\bullet};<0mm,0mm>*{}**@{},
<0mm,-2.5mm>*{};<0mm,6mm>*{}**@{-},
\end{xy}
\ \  + \ \
\begin{xy}
 <0mm,-2.1mm>*{\bullet};<0mm,0mm>*{}**@{},
<3mm,2.1mm>*{\bullet};<0mm,0mm>*{}**@{},
<3mm,2.1mm>*{};<0mm,-2.1mm>*{}**@{-},
<0mm,-2.1mm>*{};<-3mm,2.3mm>*{}**@{-},
\end{xy}
\ \  -  \ \
\begin{xy}
 <0mm,0mm>*{\bullet};<0mm,0mm>*{}**@{},
<0mm,0mm>*{};<0mm,-2.5mm>*{}**@{-},
<0mm,0mm>*{};<2.5mm,2.5mm>*{}**@{-},
<0mm,0mm>*{};<-2.5mm,2.5mm>*{}**@{-},
<0mm,0mm>*{};<2.5mm,2.5mm>*{}**@{-},
(2.4,2.3)*{}
   \ar@{->}@(ur,dr) (0,-2.3)*{}
\end{xy},
\\
\delta\,
\begin{xy}
 <0mm,1mm>*{\bullet};<0mm,0mm>*{}**@{},
<0mm,-2mm>*{};<0mm,1mm>*{}**@{-},
\end{xy}
&=&
\begin{xy}
 <0mm,3.5mm>*{\bullet};<0mm,0mm>*{}**@{},
<0mm,0mm>*{\bullet};<0mm,0mm>*{}**@{},
<0mm,-3.5mm>*{};<0mm,3mm>*{}**@{-},
\end{xy}
\ \  +  \ \
\begin{xy}
 <0mm,2.1mm>*{\bullet};<0mm,0mm>*{}**@{},
<-3mm,-2.1mm>*{\bullet};<0mm,0mm>*{}**@{},
<-3.5mm,-2.5mm>*{};<0mm,2.3mm>*{}**@{-},
<3.5mm,-2.5mm>*{};<0mm,2.3mm>*{}**@{-},
\end{xy}
\ \ -  \ \
\begin{xy}
 <0mm,0mm>*{\bullet};<0mm,0mm>*{}**@{},
<0mm,0mm>*{};<0mm,-2.5mm>*{}**@{-},
<0mm,0mm>*{};<2.5mm,2.5mm>*{}**@{-},
<0mm,0mm>*{};<-2.5mm,2.5mm>*{}**@{-},
<0mm,0mm>*{};<2.5mm,2.5mm>*{}**@{-},
(2.4,2.3)*{}
   \ar@{->}@(ur,dr) (0,-2.3)*{}
\end{xy},
&
\delta\,
\begin{xy}
 <0mm,0mm>*{\bullet};<0mm,0mm>*{}**@{},
<0mm,-3mm>*{};<0mm,3mm>*{}**@{-},
\end{xy}
&=&
\begin{xy}
 <0mm,2mm>*{\bullet};<0mm,0mm>*{}**@{},
<0mm,-2mm>*{\bullet};<0mm,0mm>*{}**@{},
<0mm,-5mm>*{};<0mm,5mm>*{}**@{-},
\end{xy}
\ \  +  \ \
\begin{xy}
 <0mm,2mm>*{\bullet};<0mm,0mm>*{}**@{},
<-3mm,-2mm>*{\bullet};<0mm,0mm>*{}**@{},
<-3.5mm,-2.5mm>*{};<0mm,2.3mm>*{}**@{-},
<3.5mm,-2.5mm>*{};<0mm,2.3mm>*{}**@{-},
<0mm,2.1mm>*{};<0mm,6mm>*{}**@{-},
\end{xy}
\ \ +  \ \
\begin{xy}
 <0mm,-2mm>*{\bullet};<0mm,0mm>*{}**@{},
<3mm,2mm>*{\bullet};<0mm,0mm>*{}**@{},
<3mm,2.1mm>*{};<0mm,-2.1mm>*{}**@{-},
<0mm,-2.1mm>*{};<-3mm,2.3mm>*{}**@{-},
<0mm,-2mm>*{};<0mm,-6mm>*{}**@{-},
\end{xy}
\ \  -  \ \
\begin{xy}
 <0mm,0mm>*{\bullet};<0mm,0mm>*{}**@{},
<0mm,0mm>*{};<2.5mm,-2.5mm>*{}**@{-},
<0mm,0mm>*{};<-2.5mm,-2.5mm>*{}**@{-},
<0mm,0mm>*{};<2.5mm,2.5mm>*{}**@{-},
<0mm,0mm>*{};<-2.5mm,2.5mm>*{}**@{-},
<0mm,0mm>*{};<2.5mm,2.5mm>*{}**@{-},
(2.4,2.3)*{}
   \ar@{->}@(ur,dr) (2.4,-2.3)*{}
\end{xy}.
\end{array}
\]

Let us now show that there is an one-to-one correspondence\footnote{
In fact this correspondence can be used as another proof of the claim
that $\delta$ is a differential, cf.\ \S 2.5 in \cite{me3}.} between
wheeled representations of $\mbox{\sf PolyV}^\circlearrowright$
and solutions of the Master equations (\ref{meq}).
Let
$$
R: (\mbox{\sf PolyV}^\circlearrowright, \delta) \lon (\End_{M_0}, d_0)
$$
be a wheeled representation of $\mbox{\sf PolyV}^\circlearrowright$,
in a differential superspace $(M_0, d_0)$.
If we forget about compatibility with differentials, then, according to
Proposition~9, any such a representation is uniquely determined
by a representation,
$$
R: \mbox{\sf PolyV} \lon \End_{M_0},
$$
which in turn is completely determined by its values,
$R^m_n\in \Hom({\odot^ n}M_0, \wedge^m M_0)$,
on the $(m,n)$-corollas. In the chosen basis $\{e_\al\}$ of $M_0$
and the associated dual basis $\{x^\al\}$ of $M_0^*$, every
such $R^m_n$ can be decomposed as,
$$
R_n^m=\sum R^{\al_1\ldots\al_m}_{\be_1\ldots\be_n} (x^{\be_1}
\odot\ldots\odot x^{\be_n})\ot( e_{\al_1}\wedge\ldots\wedge e_{\al_m})
$$
for some constants $R^{\al_1\ldots\al_m}_{\be_1\ldots\be_n}$.
The main idea of the proof is to assemble these constants for all $m,n\geq 0$ into
a single generating function, $S$, on the formal supermanifold $\cM=\Pi \Omega^1 M$
with coordinates $x^\al$ and $\psi_\be:=\Pi e_\be$, as follows
$$
S:= \vec{d}_0 + \sum_{m,n\geq 0}
\frac{1}{m!n!} R^{\al_1\ldots\al_m}_{\be_1\ldots\be_n}
x^{\be_1}\ldots x^{\be_n} \psi_{\al_1}\ldots \psi_{\al_m} \in \f_{\cM}.
$$
It is now straightforward to check that compatibility of the morphism
$R$ with differentials, $R\delta=d_0 R$, is equivalent to the
master equation (\ref{meq}) for $S$.
\end{proof}

Thus germs of master functions on $SP$-manifolds are nothing but representations
of the differential wheeled PROP
$(\mbox{\sf PolyV}^\circlearrowright, \delta)$. The class of master functions arising in this way
is in a sense typical. Theoretical physicists use
its infinite dimensional analogue (with badly defined ``volume
forms" and divergent contractions of ``dual tensors").

\begin{remark}
There is an important class of dGBV algebras with the divergence operator $\Delta$ originating
 from a particular graded metric or symplectic structure on the underlying supermanifold $M$
(i.e., $\Delta$ is a kind of ``odd" Laplacian). In this case tangent and cotangent bundles on $M$ are canonically
isomorphic so that the associated master equations are better described with the help
 of modular operads \cite{Ba} rather than wheeled PROPs.
\end{remark}


\section{Bar-Cobar duality for wheeled properads}
\label{pekelna_hudba}

\subsection{Augmentations}
The bar and cobar constructions are basic tools to study homological
properties of algebraic objects.
The bar construction considered in this paper will be
a functor from the category of augmented
wheeled dg properads
to the category of coaugmented wheeled dg co-properads,
and the cobar construction a functor from the category of
coaugmented wheeled dg co-properads to
the category of augmented wheeled dg properads.
Wheeled properads were introduced in
Section~\ref{Chicago}. In the definition below, $\frakk_c^\wc$ is the
wheeled properad introduced in Example~\ref{augment}.

\begin{definition}
An {\em augmented\/} wheeled properad is a wheeled properad
$\sfP$ together with
a homomorphism $\epsilon :  \sfP \to \frakk_c^\wc$. The kernel
$\overline {\sfP} := \Ker(\epsilon) \subset \sfP$ is the {\em
augmentation ideal\/} of $\sfP$.
\end{definition}

An example is the free wheeled
properad $\WrmF(E)$ with the augmentation $\epsilon :
\WrmF(E) \to \frakk_c^\wc = \WrmF(0)$ induced by the map $E \to
0$. Its augmentation ideal equals
\begin{equation}
\label{jupi}
\bWrmF(E) = \colim{{G \in \bDGrwc mn}}{{E}}(G),\
m,n \geq 0,
\end{equation}
where $\bDGrwc mn$ is the category of wheeled directed connected
non-exceptional graphs.

Wheeled co-properads are defined by dualizing
axioms of wheeled properads.
A more precise definition can be given as
follows. First, by inverting arrows
in Figures~\ref{triple} and~\ref{Talgebra}, one recovers the
classical notions of cotriple and
coalgebras over a cotriple.
Next, there exists a cotriple $\cWname :
\sigmabimod \to \sigmabimod$ whose underlying endofunctor is the same
as the underlying endofunctor of the triple $\WrmFname:
\sigmabimod \to \sigmabimod$ in~(\ref{l}), that is
\[
\cW(E)(m,n) := \colim{{G \in \DGrwc mn}}{{E}}(G),\
m,n \geq 0,
\]
Let us indicate how the cotriple structure
$\nu : \cWname \to \cWname\cWname$
can be defined.

Let $E$ be a $\Sigma$-bimodule. As in Section~\ref{Chicago}, we may
imagine elements of $\cW({\cW(E)})$ as nested $E$-decorated graphs, with
nests encompassing decorated subgraphs that represent
elements of the ``internal'' $\cW(E)$.
Then, for
an element $c \in \cW(E)$
represented by a directed graph
$G\in \DGrwcnic$ with $E$-labeled vertices,
$\nu_E(c) \in  \cW({\cW(E)})$ is the sum
\[
\nu_E(c) = \sum_n c_n
\]
taken over all nestings of the graph $G$, with $c_n$ being $c$
interpreted as an element of $\cW({\cW(E)})$ in the way determined
by the nesting $n$.
Because $G$ is connected, such
nestings are in one-to-one correspondence with markings of edges (in
the sense introduced below) and an alternative formula for $\nu_E(c)$
can be given.

\begin{definition}
A {\em marking\/} of $G$ is a map
$m : \Edg(G) \to \{\circ,\bullet\}$
from the set of internal edges of $G$ into the two-element set
$\{\circ,\bullet\}$. We call edges from
$\Edg^\circ(G) := m^{-1}({\circ})$ {\em
white edges\/} and edges from $\Edg^\bullet(G): =m^{-1}({\bullet})$ {\em
black edges\/}.
\end{definition}

Suppose we are given a marking $m$ of $G$ as above.
Let $G^\bullet$ be the graph
obtained from $G$ by cutting all white edges in the middle and let
$G^\bullet_1,\ldots, G^\bullet_s$, $s \geq 1$, be the connected
components of $G^\bullet$. The $E$-decoration of vertices of $G$
restricts to decorations of vertices of $G_i^\bullet$, and these decorated
graphs therefore determine elements  $c_i \in \cW(E)$, $1 \leq i
\leq s$. Let $G^\circ$ be the quotient $G/\Edg^\bullet(G)$ given by
contracting black edges of $G$.

Vertices
$v_1,\ldots,v_s$  of
$G^\circ$ are in one-to-one correspondence with the graphs
$G^\bullet_1,\ldots, G^\bullet_s$ and each $c_i$ induces
a $\cW(E)$-decoration of the
vertex $v_i \in \Edg(G^\circ)$, $1 \leq i \leq s$. The
graph $G^\circ$ with this
decoration then determines an element $c_m
\in  \cW(\cW(E))$. One has
\[
\nu_E(c) := \sum_m c_m,
\]
where the sum runs over all markings of the graph $G$.

\begin{definition}
\label{coprop}
{\em Wheeled co-properads\/} are coalgebras over the cotriple
$\cWname : \sigmabimod \to \sigmabimod$.
\end{definition}

Observe that, for each $\Sigma$-bimodule $E$,
$\cW(E)$ is a wheeled co-properad, with the co-properad structure given by
the cotriple map $\nu_E : \cW(E)
\to \cW({\cW(E)})$.

Important examples of wheeled co-properads can be obtained by taking
linear duals of wheeled properads satisfying a mild finiteness
assumption. We say that a graded wheeled
(co-)properad $\sfP$ is of {\em finite type\/} if the graded space
$\sfP(m,n)$ is of finite type for each $m,n \geq 0$.
The (componentwise) linear duals of
wheeled properads of finite type are wheeled
co-properads, and
the linear duals of wheeled co-properads of finite type are wheeled
properads, compare Remark~II.3.4
of~\cite{markl-shnider-stasheff:book}. We denote the linear dual of a wheeled
(co-)properad $\sfP$ by $\sfP^\#$.

The linear dual ${\frakk^\wc_c}^\#$ of the trivial
wheeled properad ${\frakk_c^\wc}$ turns out to be
the terminal object of the category of co-properads.
Observe that ${\frakk_c^\wc}^\# = \cW(0)$.

\begin{definition}
A {\em coaugmented\/} wheeled co-properad
is a wheeled co-properad $\sfC$ equipped with
a homomorphism $\eta : {\frakk_c^\wc}^\# \to
\sfC$. Its {\em coaugmentation
coideal\/} is the coimage $\bsfC:=
\Coim(\eta)$.
\end{definition}

The wheeled co-properad $\cW(E)$ is coaugmented,  with the
coaugmentation $\eta :  {\frakk^\wc}^\#_c \to \cW(E)$ induced by
the map $0 \to E$. Its
coaugmentation coideal equals
\[
\bcW(E) = \colim{{G \in \bDGrwc mn}}{{E}}(G),\
m,n \geq 0,
\]
where $\bDGrwc mn$ is as in~(\ref{jupi}).

The
{\em wheeled  suspension\/} of
a $\Sigma$-bimodule $E = \bicol E0$ is a $\Sigma$-bimodule
$$\ww E =\bicol{\ww E}0$$
with components
\[
\ww E(m,n) =  \susp^{2m-n} E(m,n) \ot \sgn_n,\
m,n \geq 0.
\]
Its {\em wheeled desuspension\/} $\ww^{-1}E$ is defined by a similar
formula:
\[
\ww^{-1} E(m,n) =  \desusp^{2m-n} E(m,n) \ot \sgn_n,\
\mbox { for }
m,n \geq 0.
\]
The origin of the above formulas will be explained in Remark~\ref{i2}
below.

\subsection{Bar and cobar}
We are ready to introduce the bar construction.
Let $\sfP = (\sfP,\pa_\sfP)$
be an augmented wheeled dg properad with a degree $1$
differential $\pa_\sfP$.
Consider the wheeled co-properad
$\cW(\ww^{-1} \bsfP)$ cogenerated by the
$\Sigma$-bicollection $\ww^{-1}\overline{\sfP}$, where $\bsfP$
is the augmentation ideal of $\sfP$.
The differential $\pa_\sfP$
of $\sfP$ induces, in the standard way, a degree $1$
coderivation $\pa_{1}$ of $\cW(\ww^{-1}\bsfP)$.
The structure
operation $\WrmF(\bsfP) \to \sfP$
of the wheeled properad $\sfP$ induces, precisely as
in the operadic case (see~\cite{markl:zebrulka}),
another degree $1$ coderivation $\pa_2$ on
$\cW(\ww^{-1}\bsfP)$. It is not difficult to show that
both $\pa_1$ and $\pa_2$ are
differentials commuting with each other.

\begin{definition}
\label{bar}
The wheeled dg co-properad
$(\cW(\ww^{-1} \bsfP),\pa_B)$ with  $\pa_B :=
\pa_1 + \pa_2$ will be called the {\em bar construction\/} of $\sfP$ and
denoted $\wB(\sfP) = (\wB(\sfP),\pa_B)$.
\end{definition}

\begin{remark}
\label{i2}
\def\frakK{{\mathfrak K}} \def\Det{\mathrm{Det}} \def\sss{{\mathbf s}}
Definition~\ref{bar} is a modification of the bar construction (called
the Feynman transform) of a modular (co)operad,
see~\cite[Definition~II.5.58]{markl-shnider-stasheff:book} or the
original source~\cite{getzler-kapranov:CompM98}. Unlike the bar
construction of an ordinary operad, the Feynman transform of a modular
operad is a `twisted' modular operad. The twisting is specified by
the {\em dualizing cocycle\/}
$\frakK$~\cite[Example~II.5.52]{markl-shnider-stasheff:book}
\[
\frakK(G) := \desusp^{|\Edg(G)|} \Det(\Edg(G)) =
\desusp^{|\Edg(G)|} \Lambda^{|\Edg(G)|}\left(\langle Edg(G)\rangle\right),
\]
the determinant of the span of the set $\Edg(G)$ of internal
edges of $G$, placed in degree $-|\Edg(G)|$.

When $G$ is directed, $\Edg(G)$  is clearly canonically isomorphic to
the union $\bigcup_{v \in \Vert(G)} {\Out}(v)$ of outgoing edges of vertices
of $G$, minus the set of outgoing legs of $G$. Consequently, on
directed graphs, $\frakK$ is a coboundary (in the sense
of~\cite[Lemma~II.5.49]{markl-shnider-stasheff:book}) and hence its action on a decorated graph $G$ is equivalent to
the tensor product of the decoration $\Sigma$-bimodule $E(m,n)$ with the
$\Sigma$-bimodule $u(m,n) :=\sgn_m\ot
\desusp^m \bfk \otimes\one_n $.

The next step is to modify the degrees in such a
way that the bar construction of the quadratic dual of a
wheeled Koszul properad will be concentrated in degree zero. This can
be achieved by taking the tensor product with the desuspension wheeled properad $\mathsf S^{-1}$ (see \S \ref{suspension}),
%
%
\[
(\sgn_m\ot
\desusp^m \bfk \otimes\one_n)\ot
(\sgn_m\ot \desusp^{m -n}{\mathbf k} \otimes \sgn_{n})
= \one_m\ot\desusp^{2m-n}{\mathbf k} \otimes \sgn_n.
\]
We recognize the formula for the wheeled desuspension $\ww^{-1}$ which we use in the bar construction.
\end{remark}

Let us look more closely at the structure of $\wB(\sfP)$. Elements of
$\wB(\sfP)(m,n)$ are represented by linear combinations of graphs
$G \in \DGrwc mn$ with vertices decorated by appropriately
desuspensed elements of $\bsfP$. The
degree of an element $x \in \wB(\sfP)(m,n)$
with the underlying  graph $G$ is
\[
\sum_{v \in \Vert(G)} \deg(p_v) - 2m+n - e(G),
\]
where $e(G)$ is the number of internal edges of $G$ and
$\deg(p_v)$ the degree of the decoration $p_v \in \bsfP$
of a vertex $v$ of $G$.
The differential $\pa_1(x)$ decomposes into the sum
\[
\pa_1(x) = \sum_{v \in \Vert(G)}\epsilon^v \pa_v(x),
\]
where $\pa_v$ replaces the decoration
$p_v$ of $v$ by $\pa_\sfP(p_v)$ and $\epsilon^v$ is an
appropriate sign.
Likewise, $\pa_2(v)$ decomposes as
\begin{equation}
\label{po}
\pa_2(x) = \sum_{e \in \Edg(G)}\epsilon^e \pa_e(x),
\end{equation}
where $\pa_e$ acts as follows.

If the edge $e$ connects two different vertices, then $\pa_e$
contracts $e$ and decorates the new
vertex obtained by contracting $e$ by the
properadic composition of the decorations of the vertices connected by
$e$. If $e$
is a directed loop starting and ending in the same vertex, then $\pa_e$
removes $e$ and decorates the vertex by the corresponding
contraction of the original decoration.
As before, $\epsilon^e$ is an appropriate sign.

The cobar construction is defined in the dual manner. For a wheeled dg
coaugmented co-properad $\sfC = (\sfC,\pa_\sfC)$, we consider the free
wheeled properad $\fW(\ww\bsfC)$ generated by the wheeled suspension of
the coaugmentation coideal of $\sfC$. The differential $\pa_\sfC$
induces a degree $1$ derivation $\pa_1$ and the structure operations
of the co-properad $\sfC$ a degree $1$ derivation $\pa_2$. As before,
both $\pa_1$ and $\pa_2$ are differentials that commute with each other.

\begin{definition}
\label{cobar}
The wheeled dg properad
$(\fW(\ww \bsfC),\pa_\Omega)$ with  $\pa_\Omega :=
\pa_1 + \pa_2$ will be called the {\em cobar construction\/} of $\sfC$ and
denoted $\wO(\sfC) = (\wO(\sfC),\pa_\Omega)$.
\end{definition}

The bar construction of Definition~\ref{bar} clearly extends to a functor
\[
\wB : \mbox {\tt dg-Proper}^\wc_+ \to \mbox {\tt dg-coProper}^\wc_+
\]
from the category of augmented wheeled dg-properads to the category of
augmented wheeled dg-co-properads, while the cobar construction of
Definition~\ref{cobar} extends to a functor
\[
\wO :   \mbox {\tt dg-coProper}^\wc_+ \to \mbox {\tt dg-Proper}^\wc_+.
\]
The following proposition shows that $\wB$ and $\wO$, restricted to
suitable subcategories, are exact functors.

\begin{proposition}
\label{za_tri_dny_odlet_z_Chicaga}
Let $\alpha : \sfP' \to \sfP''$ be a homology isomorphism
of non-positively graded augmented dg properads. Then the induced map
$\wB(\alpha) : \wB(\sfP') \to \wB(\sfP'')$ of the bar constructions
is a homology isomorphism, too.

Similarly, the map $\wO(\beta): \wO(\sfC') \to \wO(\sfC'')$
induced by a homology isomorphism $\beta: \sfC' \to \sfC''$ of
non-negatively graded dg coaugmented co-properads is also
a homology isomorphism.
\end{proposition}

\begin{proof}
The proposition follows from a simple spectral sequence argument based
on the filtration given by the number of edges of the underlying
graphs. The assumption about the non-positivity (resp.~non-negativity)
of the grading guarantees that the induced spectral sequences
converge.
\end{proof}

The following theorem is a key technical ingredient in showing that  $\wB$ and
$\wO$ restricted to suitable subcategories
are mutual homology inverses.

\begin{theorem}
\label{ll}
For each dg wheeled augmented properad $\sfP$,
the natural projection of wheeled dg properads
$\pi: \wO( \wB(\sfP)) \to \sfP$ is a quasi-isomorphism.

Dually, for a wheeled coaugmented
dg co-properad $\sfC$, the natural inclusion
$\iota : \sfC \hookrightarrow \wB (\wO (\sfC))$ of wheeled dg
co-properads is a quasi-isomorphism.
\end{theorem}

\begin{proof}
Let us first explicitly describe the epimorphism $\pi$.
As a non-dg properad, $\wO (\wB (\sfP))$ is free, generated by the
$\Sigma$-bimodule
$\ww \bcW(\ww^{-1} \bsfP)$. Therefore each wheeled properad
homomorphism $ \wO( \wB(\sfP)) \to \sfP$ is determined by a map
$\ww\bcW(\ww^{-1}\bsfP) \to \sfP$ of $\Sigma$-bimodules.
We define $\pi$ as the homomorphism corresponding to the composition
\[
\ww\bcW(\ww^{-1}\bsfP)
\stackrel{\ww( p)}{\longrightarrow} \bsfP \hookrightarrow \sfP,
\]
where $\ww(p)$ is the wheeled suspension of the projection
$p : \bcW(\ww^{-1} \bsfP) \to \ww^{-1} \bsfP$
to the space of cogenerators.
We will show that the homology of $(\wO( \wB(\sfP)),\pa_\Omega)$ is
isomorphic to the homology of $\sfP$. It will be clear from our proof
that this isomorphism is induced by the map $\pi$ constructed above.

Fix $m,n \geq 0$. As before, elements of
$(\wO(\wB(\sfP))(m,n) = \fW(\ww \bcW(\ww^{-1} \bsfP))(m,n)$ are
represented by nested graphs $G \in \DGrwc mn$ with vertices decorated by
appropriate (de)suspensions of elements of $\bsfP$.
Since $G$ is connected, its
nestings can be equivalently described by markings $m : \Edg(G)
\to \{\circ,\bullet\}$. The connected components of the graph
$G^\bullet$ obtained from $G$ by
cutting all white edges in the middle, with the
induced decoration of vertices, determine
elements of $\bcW(\ww^{-1} \bsfP)$ whose suspensions decorate
vertices of the quotient $G^\circ := G/\Edg^\bullet(G)$.
The degree of an element
with the underlying marked graph $G$ equals
\[
-
e^\bullet(G) + \sum_{v \in \Vert(G)} \deg(p_v),
\]
where $e^\bullet(G)$ is the number of black edges of $G$ and
$\deg(p_v)$ the degree of the decoration $p_v \in \bsfP$
of a vertex $v \in \Vert(G)$.

Recall that the differential $\pa_\Omega$ of $\wO( \wB(\sfP))$
decomposes as $\pa_\Omega = \pa_\sfP + \pa_1 + \pa_2$, where $\pa_\sfP$ is the differential
induced by the differential in  $\sfP$, the sum $\pa_\sfP + \pa_1$ is the differential
 induced by
the differential $\pa_B$ of $\wB(\sfP)$, and $\pa_2$ is induced by the co-properad structure of  $\wB(\sfP)$.
Obviously the part $\pa_\sfP + \pa_1$ of the total differential
does not change the
number of white edges, while $\pa_2$  increases it by one.
In fact, the differential $\pa_2$ applied to an element $x \in
\wO(\wB(\sfP))$ with the underlying marked graph $G$ equals
\[
\pa_2(x) = \sum_{e \in \Edg^\bullet(G) } \epsilon^e\pa_e(x),
\]
where $\pa_e$ changes the color of $e$ from black to white, and
$\epsilon^e$ is an appropriate sign.

Let us consider the spectral sequence $E = \{(E^r,d^r)\}_{r \geq 0}$
induced by an increasing filtration of the complex $\wO( \wB(\sfP))$
\[
 0\subset F_0\subset F_1 \subset \ldots F_k\subset F_{k+1}\subset \ldots
\]
with $F_k$ being the subspace of  $ \wO( \wB(\sfP))$ spanned by decorated graphs, $G$, with
at most $k$ total number of internal edges, i.e.\
$$
\mbox{number of black edges of}\ G \ + \
\mbox{number of white edges of}\ G \leq k.
$$
As the filtration is bounded below and exhaustive, the spectral sequence
$E = \{(E^r,d^r)\}_{r \geq 0}$ converges to the cohomology of $\wO( \wB(\sfP))$.
To prove the first part of the Theorem we have to show that this cohomology equals to the
cohomology of the wheeled properad $\sfP$.
The 0th term of this spectral sequence, $(E^0, d^0)$ has the differential
$d^0=\pa_\sfP + \pa_2$. To compute $E^1$ we notice that the complex $E^0$ sptits into a {\em direct}\, sum
\begin{equation}
\label{Jitulka1}
E^0= \bigoplus_{[G] \in [{\mathsf G}_c^\circlearrowright]}\wO(\wB(\sfP))_{[G]},
\end{equation}
where
$[{\mathsf G}_c^\circlearrowright]$ is the set of isomorphism classes of graphs in
${\mathsf G}_c^\circlearrowright$ and $\wO(\wB(\sfP))_{[G]}$ is the subspace of $\wO(\wB(\sfP))$
spanned by elements with the underlying graph isomorphic to $G$.
The differential $d^0=\pa_\sfP+\pa_2$ clearly preserves this decomposition.
Thus to compute $E^1$ it is enough to compute cohomology of the complex
$$
(E^0_G:=\wO(\wB(\sfP))_{[G]}, d_\sfP+\pa_2)
$$
for any particular fixed $(m,n)$-graph $G$. This complex is
actually a bicomplex, $E^0_G=\{E_G^{p,q}\}$ with $p$ being the total $\sfP$-degree of $\sfP$-decorate vertices,
and $q$ is the number of white vertices. Notice that $q$ ranges from zero to the total number of edges in $G$.

The filtration of the complex  $E^0_G$ by the number of white
vertices is bounded, hence the induced spectral sequence $(\cE_G^r,
\delta^r)$ converges
to the cohomology of $E^0_G$. The 0th term, $\cE^0_G$, of this sequence has the differential $\delta^0$
 equal to $d_\sfP$.
Hence $\cE^1_G$ is spanned by graphs $G$ decorated with elements of the cohomology wheeled properad $H(\sfP)$,\and
$\delta^1$ is equal to $\pa_2$, the differential in $\wO( \wB(H(\sfP)))$ coding the coproperad structure
in  $\wB(H(\sfP))$.

To compute $\cE^2_G=H(\cE^1_G, \delta^1)$
for a non-exceptional graph
$G \in \DGrwc mn$, let $\{e_1,\ldots,e_s\} = \Edg(G)$ be the
set of edges of $G$. It is a simple exercise to prove that
the exterior algebra
\begin{equation}
\label{jituska}
(\ext(e_1,\ldots,e_s),\pa_\land), \ \deg(e_i) := -1 \ \mbox { for }
\  1 \leq i \leq s,
\end{equation}
with the differential $\pa_\land
:= \sum_{1 \leq i \leq s} \pa/\pa e_i$, is {\em acyclic\/} whenever $s
\geq 1$.
Observe also that $\pa_\land$ is equivariant under the action of the
symmetric group of $\{e_1.\ldots,e_s\}$.
The acyclicity of~(\ref{jituska}) implies that also the product
\begin{equation}
\label{prod}
(\overline{H(\sfP)}(G) \ot \ext(e_1,\ldots, e_s),\pa)
\end{equation}
in which $\overline{H(\sfP)}(G)$ is as in~(\ref{pomalu}) and
$\pa(\overline x \ot u) := (-1)^{\deg(\overline x)} \overline x \ot
\pa_\land u$, for $\overline x \in \overline{H(\sfP)}(G)$ and $u \in \ext(e_1,\ldots,
e_s)$, is acyclic.

The group $\Aut(G)$ of automorphisms of $G$ acts on
$\ext(e_1,\ldots, e_s)$ by
permutations of edges, which implies that  $\pa_\land$
is $\Aut(G)$-equivariant. Therefore also the differential $\pa$ on the
product~(\ref{prod}) is $\Aut(G)$-equivariant,
hence it induces a differential $\pa_G$ on
the orbit space
\[
\left(\overline{H(\sfP)}(G) \ot
\ext(e_1,\ldots,e_s)\right)_G := \left(\overline{H(\sfP)}(G) \ot
\ext(e_1,\ldots,e_s)\right)/\Aut(G).
\]
Since $\Aut(G)$ is a finite group and the ground field has
characteristic zero, the acyclicity of~(\ref{prod}) implies the
acyclicity of $((\overline{H(\sfP)}(G) \ot
\ext(e_1,\ldots,e_s))_G, \pa)$, for each $s \geq 1$.

There is an isomorphism
\begin{equation}
\label{pp}
\Psi_G:
\left(\wO(\wB(H(\sfP)))_{[G]},\pa_2\right)
\stackrel{\cong}\to
\left(\left(\overline{H(\sfP)}(G) \ot \ext(e_1,\ldots,e_s)\right)_G,\pa_G\right)
\end{equation}
that sends the isomorphism class of an element $x \in
\wO(\wB(H(\sfP)))_{[G]}$ represented by a $\overline{H(\sfP)}$-decorated graph $G$
with marking $m$, into the orbit of
\[
\overline x \ot (e_{i_1} \land \cdots \land e_{i_t}) \in \overline{H(\sfP)}(G) \ot
\ext(e_1,\ldots,e_s),
\]
where $\overline x \in \overline{H(\sfP)}(G)$
is obtained from $x$ by forgetting the marking of $G$,
and
\begin{equation}
\label{mas}
\{e_{i_1},\ldots,e_{i_t}\} =
\Edg^\bullet(G),\  i_1 < \cdots < i_t.
\end{equation}
It is easy to see that $\Psi$ is
well-defined and that it commutes with the differentials. Its inverse
\[
\Psi^{-1}_G:
\left(\left(\overline{H(\sfP)}(G) \ot \ext(e_1,\ldots,e_s)\right)_G,\pa_G\right)
\stackrel{\cong}\to
\left(\wO(\wB(H(\sfP)))_{[G]},\pa_2\right)
\]
maps $\overline {x} \otimes
e_{i_1} \land \cdots \land e_{i_t}$, where $\overline {x} \in
\overline{H(\sfP)}(G)$ is  represented by a $\overline{H(\sfP)}$-decoration
of $G$, into $x \in \wO(\wB(H(\sfP)))_{[G]}$ given by the same
$\overline{H(\sfP)}$-decorated graph as $\overline{x}$ but with marking
defined by~(\ref{mas}).

We conclude from~(\ref{pp}) that $\wO(\wB(H(\sfP)))_{[G]}$ is $\pa_2$-acyclic
whenever $G$ has at least one internal edge. If $G$ has no internal
edge, i~e.~when $G$ is an $(m,n)$-corolla, then clearly
\[
H(\wO(\wB(H(\sfP)))_{[G]},\pa_2) \cong \overline{H(\sfP)}(G)(m,n).
\]
The two exceptional graphs $\uparrow$ and $\hgw$ contribute
to~(\ref{Jitulka1}) by $\frakk^\wc_c(1,1)$ if $(m,n) = (1,1)$
and by $\frakk^\wc_c(0,0)$ if $(m,n) = (0,0)$. As all the higher diffferentials in both of our spectral sequences
are obviously zero,
we  obtain an isomorphism,
\[
H(\wO(\wB(\sfP))) \cong \overline{H(\sfP)} \oplus \frakk^\wc_c \cong
H(\sfP)
\]
proving the first claim in the Theorem that the surjection $\pi$
is a quasi-isomorphism.

The  second claim on the map $\iota  : \sfC \to \wB (\wO (\sfC))$ can be proven by dualizing the above arguments.
\end{proof}


\section{Wheeled operads, quadratic duality and Koszulness}
\label{Tomicek}

\subsection{Wheeled operads}
As operads form a subcategory of \PROP{s}~\cite{handbook},
there is a similar
subcategory of wheeled operads in the category of wheeled \PROP{s}. In
this section we introduce these wheeled operads, define their quadratic
duals and study their Koszulness. Methods of this
section will provide a conceptual understanding of minimal models of
wheeled completions of classical quadratic operads as
$\Ass$, $\Lie$ or $\Com$.

\begin{definition}
A {\em wheeled operad\/} is a wheeled properad $\calP = \bicol \calP0$,
in the sense of
Definition~\ref{wheeled_properad}, such that
$\calP(m,n) = 0$ whenever $m \geq 2$.
\end{definition}

It is obvious from this definition
that a wheeled operad $\calP = \bicol \calP0$ consists of
\begin{itemize}
\item[(i)]
an ordinary operad $\calP_o := \{\calP(1,n)\}_{n \geq 0}$,
\item[(ii)]
a right $\calP_o$-module
$\calP_w := \{\calP(0,n)\}_{n \geq 0}$
(see~\cite[Def.~II.3.26]{markl-shnider-stasheff:book} for a definition
of right operadic modules), and
\item[(iii)]
contractions $\xi_i : \calP_o(n) \to \calP_w(n-1)$, $1 \leq i \leq
n$, that are compatible, in the appropriate sense, with the structures
in~(i) and~(ii).
\end{itemize}

\begin{definition}
The operad $\calP_o$ and the right $\calP_o$-module $\calP_w$ defined
above are called the {\em operadic\/} and {\em wheeled\/}
parts of the wheeled operad
$\calP$, respectively.
\end{definition}

Wheeled properadic completions of ordinary operads provide
examples of wheeled operads, but, as we will see immediately, not all wheeled
operads are of this form.
Recall~\cite[Definition~II.3.31]{markl-shnider-stasheff:book}
that an ordinary operad
$\calP$ is {\em quadratic\/} if it is of
the form
\begin{equation}
\label{tt}
\calP = \freeOP(E)/I,
\end{equation}
where $\freeOP(E)$ is the free operad on
a right $\Sigma_2$-module considered as a
$\Sigma$-module with
\[
E(n) := \cases{E}{if $n = 2$ and}{0}{otherwise,}
\]
and $I$ an operadic ideal generated by a subspace $R \subset
\freeOP(E)(3)$. The quadratic dual $\calP^!$ of a quadratic operad
$\calP$ is defined as
\begin{equation}
\label{99}
\calP^! := \freeOP(E^\vee)/(R^\perp),
\end{equation}
where $E^\vee := \Lin(E,\bfk) \ot \sgn_2$ is the Czech dual of $E$ and
$R^\perp$ the annihilator of $R$ in $\freeOP(E^\vee)(3) \cong
\freeOP(E)(3)^\vee$.
See~\cite[Definition~II.3.37]{markl-shnider-stasheff:book} for details.
Let us introduce wheeled versions of the above notions.

\begin{definition}
\label{ee}
{\em Quadratic wheeled operad\/} is a wheeled operad $\calP$ of the
form
\[
\calP =  \fW(E)/I
\]
where $E$ is a left $\Sigma_1$- right $\Sigma_2$-bimodule considered as a
$\Sigma$-bimodule with
\[
E(m,n) := \cases{E}{if $(m,n) = (1,2)$ and}{0}{otherwise,}
\]
and $I$ a wheeled operadic ideal generated by a subspace
\begin{equation}
\label{ow}
R_o \oplus
R_w \subset \fW(E)(1,3) \oplus \fW(E)(0,1).
\end{equation}
We will always assume that $E$ is finite dimensional.
\end{definition}

In the rest of this paper, we will often work with the wheeled properadic
completion of operads. To simplify the terminology and notation, we
will call it simply the wheeled completion and drop the subscript from
the notation, that is, write $\calP^\wc$ instead of $\calP^\wc_c$. We
believe that no confusion is possible.

\begin{example}
\label{Karmen}
A wheeled quadratic operad is the
wheeled completion of an ordinary quadratic operad if and
only if the space $R_w$ in~(\ref{ow}) is trivial. If $\calP$ is an ordinary
quadratic operad as in~(\ref{tt}) with $I$ generated by a subspace $R \subset
\freeOP(E)(3)$, then $\calP^\wc$ is wheeled quadratic with the same
space of generators and with $R_o = R$, $R_w = 0$.

In particular, wheeled  completions $\Ass^\wc$,
$\Com^\wc$ and $\Lie^\wc$ of quadratic operads $\Ass$,
$\Com$ and $\Lie$ for associative, commutative associative and Lie
algebras, respectively, are wheeled quadratic operads. The
initial wheeled properad $\frakk^\wc$ is a wheeled quadratic operad
generated by the trivial $\Sigma$-bimodule.
Therefore every wheeled quadratic operad $\calP$
is augmented, with the augmentation $\calP
\to\frakk^\wc$ induced by the map $E \to 0$ of generators.
\end{example}

\subsection{Quadratic duality and Koszulness}
We need to extend the definition of the Czech dual to $\Sigma$-bimodules.
For a graded left $\bfk[\Sigma_m]$- right
$\bfk[\Sigma_n]$-bimodule $U$ define $U^\vee$
by
\[
(U^\vee)^i := \sgn_m \ot \Lin(U_{-i},\bfk) \ot \sgn_n,
\]
where $\Lin(U_{-i},\bfk)$ is the ordinary linear dual of $U_{-i}$
and $\sgn$ the signum representation. The Czech dual of a
$\Sigma$-bimodule is then the componentwise application of
the above operation.

As in the operadic case
(see~\cite[Section~II.3.2]{markl-shnider-stasheff:book}) there is, for
each $\Sigma$-bimodule $E$, a natural biequivariant
isomorphism $\fW(E^\vee) \cong  \fW(E)^\vee$. Therefore the
annihilator of each subspace
$S \subset \fW(E)(m,n)$ can be regarded as a subspace
$S^\perp \subset \fW(E^\vee)(m,n)$, $m,n \geq 0$.

\begin{definition}
\label{vcera_jsem_videl_Drinfelda}
Let $\calP$ be a wheeled quadratic operad as in
Definition~\ref{ee}. Define its {\em wheeled quadratic dual\/} as the
wheeled quadratic operad
\[
\calP^! := \fW(E^\vee)/I^\perp
\]
with $I^\perp$ the ideal generated by
\begin{equation}
\label{virzinka_z_Prahy}
R_o^\perp \oplus
R_w^\perp \subset \fW(E^\vee)(1,3) \oplus \fW(E^\vee)(0,1).
\end{equation}
\end{definition}

The wheeled quadratic dual is clearly an involution in the category of
wheeled quadratic operads, $(\calP^!)^! = \calP$.
Let us prove
the following simple but
useful proposition.

\begin{proposition}
\label{tyden_pred_odletem_z_Chicaga}
Let $\calP$ be a wheeled quadratic operad as
in Definition~\ref{ee}.
Then the operadic part of its wheeled quadratic dual
equals the
ordinary quadratic dual of its operadic part
\begin{equation}
\label{Prasatko}
(\calP^!)_o \cong (\calP_o)^!;
\end{equation}
we denote it simply $\calP_o^!$.
The wheeled part of $\calP^!$ is the quotient of
the wheeled part of the wheeled  completion of $\calP^!_o$,
\begin{equation}
\label{aa}
(\calP^!)_w \cong ((\calP^!_o)^\wc)_w/I^\perp_w
\end{equation}
by the ideal $I^\perp_w$ generated by $R^\perp_w$.
\end{proposition}

\begin{proof}
The double coset theorem gives isomorphisms
\[
P^! = \fW(E^\vee)/(R_o^\perp,R_w^\perp)
\cong (\fW(E^\vee)/(R_o^\perp))/(R_w^\perp)
\cong (\calP^!_o)^\wc/I^\perp_w
\]
Since the ideal $I^\perp_w$ consists of
elements of biarities $(0,n)$, $n\geq 1$, it does change the
operadic part of the rightmost quotient.
Isomorphisms~(\ref{Prasatko}) and~(\ref{aa}) are then obtained by
singling out the operadic and wheeled parts of the above
display.
\end{proof}

Let us see what happens if we apply
Proposition~\ref{tyden_pred_odletem_z_Chicaga} to the
wheeled  completion of an ordinary quadratic
operad $\calP$ given by the quotient~(\ref{tt}).
As we observed in Example~\ref{Karmen},
$\calP^\wc$ is a wheeled
quadratic operad with $R_o = R$ and $R_w = 0$, therefore,
in~(\ref{virzinka_z_Prahy}), $R_o^\perp =
R^\perp$ and $R_w^\perp  =
\fW(E^\vee)(0,1)$.
By~(\ref{Prasatko}), $(\calP^\wc)^!_o \cong \calP^!$,
while~(\ref{aa}) gives
\begin{equation}
\label{aaa}
 ((\calP^\wc)^!)_w \cong ((\calP^!)^\wc)_w/I^\perp_w
\end{equation}
with
$I^\perp_w$ generated by $\fW(E^\vee)(0,1)$, that is,
by relations represented by the decorated graphs:
\begin{equation}
\label{jeste_strasim_na_departmentu}
\unitlength=1em
\thinlines
\begin{picture}(0,1.7)(5,-.5)
\put(0.05,0){\circle{2}}
\put(-1.5,-1){\vector(1,2){.5}}
\put(-.89,-.5){\vector(-1,1){0}}
\put(-1.05,0){\makebox(0,0){$\bullet$}}
\put(-1.3,0){\makebox(0,0)[rb]{$e'$}}
\put(5,0){
\put(-0.05,0){\circle{2}}
\put(1.5,-1){\vector(-1,2){.5}}
\put(.89,-.5){\vector(1,1){0}}
\put(1.05,0){\makebox(0,0){$\bullet$}}
\put(1.3,0){\makebox(0,0)[lb]{$e''$}}
}
\put(8,0){\makebox(0,0)[l]{$=0,\ \mbox { for } e',e'' \in E^\vee$.}}
\put(2.1,0){\makebox(0,0)[l]{$=$}}
\end{picture}
\end{equation}

In the following examples we analyze in detail the
wheeled quadratic duals of wheeled  completions of some
quadratic operads. As we saw above, this analysis necessarily
involves the ordinary quadratic duals of these operads for which we
refer to~\cite[Section~II.3.2]{markl-shnider-stasheff:book} or to the
original source~\cite{ginzburg-kapranov:DMJ94}.

\begin{example}
\label{co_mne_ceka_doma?}
{\rm
Let us describe, using~(\ref{aaa}),
the wheeled part of the wheeled quadratic dual
$(\Ass^\wc)^!$ of the wheeled  completion of the operad
$\Ass$.
Since $\Ass^! = \Ass$, the first
step is to understand the wheeled part $(\Ass^\wc)_w$ of the
wheeled  completion of $\Ass$.
It is easy to see that $(\Ass^\wc)_w(n)$ is,
for $n \geq 1$, spanned by graphs
\begin{eqnarray}
\label{snad_je_doma_vsechno_v_poradku}
\unitlength=1.2em
\thinlines
\begin{picture}(5,3.2)(-3,-1.5)
\thinlines
\put(-3,-3){\vector(1,1){3}}
\put(-1,-3){\vector(1,3){1}}
\put(3,-3){\vector(-1,1){3}}
\put(1,-3){\vector(-1,3){1}}
\put(2,1.54){\vector(1,0){0}}
\put(1.62,-0.1){\hkrouzek {1}}
\put(1.62,-.1){\dkrouzek {1}}
\put(0,0){\makebox(0,0){$\bullet$}}
\put(-5,-1){\makebox(0,0){$\xi^i_\sigma:=$}}
\put(-3,-3.1){\makebox(0,0)[t]{\scriptsize $\sigma(1)$}}
\put(-1.8,-2.7){\makebox(0,0)[t]{\scriptsize $\cdots$}}
\put(-1,-3.1){\makebox(0,0)[t]{\scriptsize $\sigma(i)$}}
\put(1.1,-3.1){\makebox(0,0)[t]{\scriptsize $\sigma(i\!\!+\!\!1)$}}
\put(1.8,-2.7){\makebox(0,0)[t]{\scriptsize $\cdots$}}
\put(3,-3.1){\makebox(0,0)[t]{\scriptsize $\sigma(n)$}}
\end{picture}
\\
\nonumber
\rule{0em}{1.5em}
\end{eqnarray}
where $\sigma \in \Sigma_n$ and $0 \leq i \leq n$, such that both the
`left' group of inputs (labeled $\sigma(1),\ldots,\sigma(i)$) and
the `right' group of inputs (labeled $\sigma(i+1),\ldots,\sigma(n)$)
is cyclically symmetric.

This cyclic symmetry can be seen as follows. By the
associativity of the multiplication in $\Ass$, one is allowed to
single out the
input edge labeled $\sigma(1)$, move it clockwise around the loop,
and join it at the bottom of the left group. This operation
generates a left action of the cyclic group $C_i := {\mathbb Z}/(i)$
on the left group
of labels. The cyclic symmetry of the
right group can be explained similarly.

We claim that $\xi^i_\sigma = 0$ modulo
relations~(\ref{jeste_strasim_na_departmentu})
whenever $i = n$ or $i=0$. Indeed, for $i=n$ one has, due
to the associativity of the multiplication
\[
\unitlength=1.2em
\thinlines
\begin{picture}(5,4.5)(2,-2)
\put(2,1.54){\vector(1,0){0}}
\put(1.62,-0.1){\hkrouzek {1}}
\put(1.62,-.1){\dkrouzek {1}}
\put(-2,-2){\vector(1,1){2}}
\put(-.05,-2){\vector(0,3){2}}
\put(0,0){\makebox(0,0){$\bullet$}}
\put(.7,.9){
\put(-3,-3.1){\makebox(0,0)[t]{\scriptsize $\sigma(1)$}}
\put(-1.8,-2.7){\makebox(0,0)[t]{\scriptsize $\cdots$}}
\put(-1,-3.1){\makebox(0,0)[t]{\scriptsize $\sigma(n)$}}
}
\put(10,0)
{
\put(2,1.54){\vector(1,0){0}}
\put(1.62,-0.1){\hkrouzek {1}}
\put(1.62,-.1){\dkrouzek {1}}
\put(-2,-2){\vector(1,1){2}}
\put(0,0){\makebox(0,0){$\bullet$}}
\put(1,1){
\put(-3,-3){\vector(1,1){1}}
\put(-1,-3){\vector(-1,1){1}}
\put(-2,-2){\makebox(0,0){$\bullet$}}
\put(-3,-3.1){\makebox(0,0)[t]{\scriptsize $\sigma(1)$}}
\put(-1.8,-2.7){\makebox(0,0)[t]{\scriptsize $\cdots$}}
\put(-1,-3.1){\makebox(0,0)[t]{\scriptsize $\sigma(n)$}}
}
}
\put(6,-.5){\makebox(0,0)[t]{$=$}}
\put(-4,-.5){\makebox(0,0)[t]{$\xi^n_\sigma =$}}
\end{picture}
\]
with the right element clearly belonging to the ideal generated
by~(\ref{jeste_strasim_na_departmentu}).
By similar reasons, $\xi^1_\sigma$
is also zero modulo~(\ref{jeste_strasim_na_departmentu}).
We conclude that, for $n \geq 1$,
\[
(\Ass^\wc)^!_w(n) \cong \bigoplus_{0 < i < n}
(C_i \times C_{n-i})\backslash \bfk[\Sigma_n],
\]
where $(C_i \times C_{n-i})\backslash \bfk[\Sigma_n]$
is the left quotient of $\bfk[\Sigma_n]$ by
the product of cyclic groups $C_i \times C_{n-i}$.
}
\end{example}

\begin{example}
\label{Kralicek}
{\rm
Our next task is to describe the wheeled part of the wheeled quadratic dual
$(\Com^\wc)^!$ of the wheeled  completion of the operad
$\Com$.
Since $\Com^! = \Lie$, we need to analyze first
the wheeled part $(\Lie^\wc)_w$ of the
wheeled  completion of $\Lie$. The crucial fact is
that the Jacoby identity
\[
\thinlines
\begin{picture}(170.00,34.00)(-15.00,5.00)
\unitlength=.7pt
\put(110.00,30.00){\makebox(0.00,0.00){$+$}}
\put(50.00,30.00){\makebox(0.00,0.00){$+$}}
\put(120.00,0.00){\makebox(0.00,0.00){\scriptsize $3$}}
\put(80.00,0.00){\makebox(0.00,0.00){\scriptsize $3$}}
\put(40.00,0.00){\makebox(0.00,0.00){\scriptsize $3$}}
\put(160.00,0.00){\makebox(0.00,0.00){\scriptsize $2$}}
\put(60.00,0.00){\makebox(0.00,0.00){\scriptsize $2$}}
\put(20.00,0.00){\makebox(0.00,0.00){\scriptsize $2$}}
\put(140.00,0.00){\makebox(0.00,0.00){\scriptsize $1$}}
\put(100.00,0.00){\makebox(0.00,0.00){\scriptsize $1$}}
\put(0.00,0.00){\makebox(0.00,0.00){\scriptsize $1$}}
\put(30.00,20.00){\makebox(0.00,0.00){$\bullet$}}
\put(90.00,20.00){\makebox(0.00,0.00){$\bullet$}}
\put(150.00,20.00){\makebox(0.00,0.00){$\bullet$}}
\put(140.00,30.00){\makebox(0.00,0.00){$\bullet$}}
\put(80.00,30.00){\makebox(0.00,0.00){$\bullet$}}
\put(20.00,30.00){\makebox(0.00,0.00){$\bullet$}}
\put(140.00,10.00){\qbezier(0,0)(5,5)(10,10)}
\put(140.00,50.00){\line(0,-1){20.00}}
\put(140.00,30.00){\line(1,-1){20.00}}
\put(120.00,10.00){\line(1,1){20.00}}
\put(80.00,50.00){\line(0,-1){20.00}}
\put(80.00,10.00){\qbezier(0,0)(5,5)(10,10)}
\put(80.00,30.00){\line(1,-1){20.00}}
\put(60.00,10.00){\line(1,1){20.00}}
\put(20.00,10.00){\qbezier(0,0)(5,5)(10,10)}
\put(20.00,30.00){\line(1,-1){20.00}}
\put(20.00,30.00){\line(-1,-1){20.00}}
\put(20.00,50.00){\line(0,-1){20.00}}
\put(170.00,20.00){$=0$}
\end{picture}
\]
for the Lie bracket symbolized by  $\brackett$ implies that, in
$(\Lie^\wc)_w(2)$,
\begin{eqnarray}
\label{zase_strasim_v_noci}
\unitlength=.9em
\thinlines
\begin{picture}(5,2.5)(5,0)
\put(0,0){\dotkrouzek {.9}}
\put(-4,-2){\vector(1,1){2}}
\put(-4,-2){\vector(1,1){1}}
\put(-2,-2){\vector(-1,1){1}}
\put(-1.8,-1){\vector(-1,1){0}}
\put(-3,-1){\makebox(0,0){$\bullet$}}
\put(-2,-0){\makebox(0,0){$\bullet$}}
\put(-4,-2.3){\makebox(0,0)[t]{\scriptsize $1$}}
\put(-2,-2.3){\makebox(0,0)[t]{\scriptsize $2$}}
\put(3.5,0){\makebox(0,0){$=$}}
\put(8,0){
\put(-1.8,-1){\vector(-1,1){0}}
\put(0,0){\dotkrouzek {.9}}
\put(-3,-1){\vector(1,1){1}}
\put(-1.75,-2){\vector(0,1){1}}
\put(-1.75,-1){\makebox(0,0){$\bullet$}}
\put(-2,-0){\makebox(0,0){$\bullet$}}
\put(-3,-1.3){\makebox(0,0)[t]{\scriptsize $1$}}
\put(-1.75,-2.3){\makebox(0,0)[t]{\scriptsize $2$}}
\put(3.5,0){\makebox(0,0){$-$}}
}
\put(16,0){
\put(-1.8,-1){\vector(-1,1){0}}
\put(0,0){\dotkrouzek {.9}}
\put(-3,-1){\vector(1,1){1}}
\put(-1.75,-2){\vector(0,1){1}}
\put(-1.75,-1){\makebox(0,0){$\bullet$}}
\put(-2,-0){\makebox(0,0){$\bullet$}}
\put(-3,-1.3){\makebox(0,0)[t]{\scriptsize $2$}}
\put(-1.75,-2.3){\makebox(0,0)[t]{\scriptsize $1$}}
}
\end{picture}
\\
\nonumber
\rule{0pt}{2em}
\end{eqnarray}
where the dotted oriented circle is the unique wheel of the
underlying graph (which may or may not contain other vertices). It
is not difficult to conclude from this
that $(\Lie^\wc)_w(n)$ is, for $n \geq 1$, spanned by elements
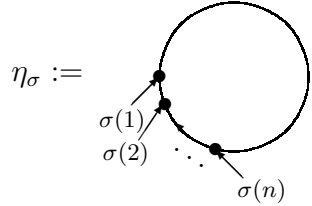
\[
\unitlength=.6em
\thinlines
\begin{picture}(5,10.5)(-4,-5.5)
\put(-3.2,-2.4){\vector(-1,1){0}}
\put(0,0){\krouzek {1.2}}
\put(-4,0){\makebox(0,0){$\bullet$}}
\put(-5.5,-1.5){\vector(1,1){1.5}}
\put(-3.7,-1.5){\makebox(0,0){$\bullet$}}
\put(-5.1,-3.4){\vector(2,3){1.3}}
\put(-1,-3.9){\makebox(0,0){$\bullet$}}
\put(.5,-5.4){\vector(-1,1){1.5}}
\put(-6,-1.7){\makebox(0,0)[t]{\scriptsize $\sigma(1)$}}
\put(-5.8,-3.5){\makebox(0,0)[t]{\scriptsize $\sigma(2)$}}
\put(1.5,-5.5){\makebox(0,0)[t]{\scriptsize $\sigma(n)$}}
\put(-2.4,-3.9){\makebox(0,0){$\ddots$}}
\put(-10,0){\makebox(0,0){$\eta_\sigma := $}}
\end{picture}
\]
with $\sigma \in \Sigma_n$.
It is also obvious that
$\eta_\sigma \in I^\perp_w$ for $n = 1$. So
\[
(\Com^\wc)_w^!(n) \cong
\cases{C_n\backslash \bfk[\Sigma_n]}{for $n \geq 2$, and}0{for $n=1$,}
\]
with the trivial action of $\Sigma_n$.
}
\end{example}

\begin{example}
\label{kasparek}
{\rm
In this example we describe the wheeled quadratic dual of the wheeled
 completion of $\Lie$. Since $\Lie^! = \Com$, we need to
start by investigating $(\Com^\wc)_w$.
The graphs that represent elements of $(\Com^\wc)_w(n)$ are the same
as in~(\ref{snad_je_doma_vsechno_v_poradku}),
but commutativity enables one to move all
inputs of $\xi^i_\sigma$ to one side (say to the left)
and arrange them into increasing order.
Therefore  $(\Com^\wc)_w(n)$ is, for $n \geq 1$, one dimensional,
spanned by the directed graph
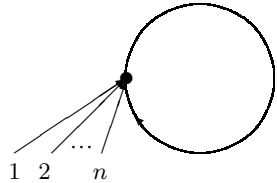
\[
\unitlength=1.2em
\thinlines
\begin{picture}(5,4.5)(-2,-2)
\put(-3,-2){\vector(3,2){3}}
\put(-.66,-2){\vector(1,3){.66}}
\put(-2,-2){\vector(1,1){2}}
\put(0.2,-1){\vector(-1,1){0}}
\put(1.97,0){\krouzek {1.2}}
\put(0,0){\makebox(0,0){$\bullet$}}
\put(0,1){
\put(-3,-3.7){\makebox(0,0)[b]{\scriptsize $1$}}
\put(-2.2,-3.7){\makebox(0,0)[b]{\scriptsize $2$}}
\put(-1.2,-2.7){\makebox(0,0)[t]{\scriptsize $\cdot \! \cdot\! \cdot$}}
\put(-.7,-3.7){\makebox(0,0)[b]{\scriptsize $n$}}
}
\end{picture}
\]
We argue as in Example~\ref{co_mne_ceka_doma?} that, by associativity,
the element
represented by the above graph belongs to $I^\perp_w$ therefore,
rather surprisingly,
\[
(\Lie^\wc)^!_w(n) = 0\ \mbox { for }\ n \geq 1.
\]
In other words, the wheeled part of $(\Lie^\wc)^!$ is {\em trivial\/}.
}
\end{example}

Since each wheeled quadratic operad $\calP$ is a particular
case of an augmented properad, it makes sense to consider its bar
construction introduced in Definition~\ref{bar}.
Let us investigate how the $0$-th homology $H_0(\wB(\calP),\pa_B)$ is
related to the wheeled quadratic dual of $\calP$.
It follows from simple combinatorics
that the cooperad $\wB(\calP) = \cW(\ww^{-1}\bcalP)$
is non-negatively graded
and that an element $x \in \wB(\calP)$ with the
underlying graph $G$ has degree $0$ if and only if all vertices of $G$
are decorated by $\calP(1,2) = E$. Therefore
\[
\wB(\calP)_0 \cong \cW(\ww^{-1}E) =\cW(E \otimes \sgn_2).
\]

It is equally obvious that the degree $1$ part
$\wB(\calP)_{1}$ is spanned by elements
whose underlying graphs have precisely one vertex decorated by
$\calP(1,3)$ and all remaining vertices decorated by $\calP(1,2)= E$.
The differential
\[
\pa_B : \cW (E  \otimes \sgn_2) \to \wB(\calP)_{1}
\]
act as in~(\ref{po}), by contracting edges of the underlying graph. Let
\[
\calP^{\mbox {\scriptsize !`}} := H_0(\wB(\calP),\pa_B)
= \Ker\left(\pa_B : \cW (E) \to \wB(\calP)_{1}\right).
\]
It is an
exercise in linear algebra to prove that the cooperad
$\calP^{\mbox {\scriptsize !`}}$ is precisely
the linear dual  the quadratic dual of $\calP$,
that is
\[
{\calP^{\mbox {\scriptsize !`}}}^\# = \calP^!.
\]

\begin{definition}
\label{wKoszul}
A quadratic wheeled operad $\calP$ is {\em wheeled Koszul\/}
if the canonical inclusion
of wheeled cooperads
\begin{equation}
\label{jsem_zvedav_jestli_mne_zamknou}
\iota : (\calP^{\mbox {\scriptsize\rm !`}},\pa =0)
\hookrightarrow (\wB(\calP),\pa_B)
\end{equation}
is a homology isomorphism.
\end{definition}

It follows from the above analysis that
$\calP$ is wheeled quadratic Koszul if and only if
\[
H_{> 0}(\wB(\calP),\pa_B) =0.
\]
Another equivalent formulation is given in the following
proposition.

\begin{proposition}
\label{Jitka_nejak_neroste}
A quadratic wheeled operad $\calP$ is wheeled Koszul if and only if the
natural projection
\begin{equation}
\label{oi}
p: (\wO(\calP^{\mbox {\scriptsize \rm !`}}),\pa_\Omega)\to  (\calP,\pa = 0)
\end{equation}
of wheeled dg operads is a homology isomorphism.
\end{proposition}

\begin{proof}
Observe that the projection~(\ref{oi}) coincides with
the composition
\[
(\wO(\calP^{\mbox {\scriptsize \rm !`}}),\pa_\Omega)
\stackrel{\wO(\iota)}{\longrightarrow}
(\wO(\wB(\calP),\pa_\Omega) \stackrel{\pi}{\longrightarrow} \calP
\]
in which $\wO(\iota)$ is induced
by the inclusion~(\ref{jsem_zvedav_jestli_mne_zamknou}) and
$\pi$ is the canonical map of Theorem~\ref{ll} which is a homology.
If $\calP$ is
wheeled Koszul then, by definition,
$\iota$ is a homology isomorphism of non-negatively graded
co-properads, hence
$\wO(\iota)$ is, by Proposition~\ref{za_tri_dny_odlet_z_Chicaga},
a homology isomorphism, too.
The projection $p$ in~(\ref{oi}) is then
a composition of homology isomorphisms so it is also
a homology isomorphism.

The opposite implication follows from a~similar analysis of the
composition
\[
(\calP^{\mbox {\scriptsize \rm !`}},\pa =0) \stackrel{c}{\to}
(\wB(\wO(\calP)),\pa_B) \stackrel{\wB(p)}{\longrightarrow}
(\wB(\calP),\pa_B).
\]
\vskip -2em
\end{proof}

The importance of Proposition~\ref{Jitka_nejak_neroste} is
that~(\ref{oi}) provides a functorial minimal resolution of wheeled
quadratic Koszul operads in the category of wheeled operads.
It is easy to verify that the left adjoint
\[
\widehat{}:\mbox {\tt Proper}^\wc \to \catWPROP
\]
to the forgetful functor
$\catWPROP \to \mbox {\tt Proper}^\wc$,  given by disjoint unions
of the underlying graphs, it is
an exact functor.  Therefore
$\hskip .3em\widehat{}\hskip .3em$ applied to~(\ref{oi})
gives a functorial minimal resolution
\[
\widehat p:
(\whO(\calP^{\mbox {\scriptsize \rm !`}}),\pa_{\widehat{\Omega}})\to
(\widehat{\calP},\pa = 0)
\]
of the wheeled \PROP\ generated by $\calP$ in the category of wheeled
\PROP{s}. It is clear that
\[
(\whO(\calP^{\mbox {\scriptsize !`}}), \pa_{\widehat{\Omega}}) =
(\freeWPROP(\ww\calP^{\mbox {\scriptsize !`}}),\pa_{\widehat{\Omega}})
\]
with $\pa_{\widehat{\Omega}}$ given on generators by the same
formula as the cobar differential $\pa_\Omega$. Notice, however, that
$(\whO(\calP^{\mbox {\scriptsize !`}}), \pa_{\widehat{\Omega}})$ is
not a cobar construction of $\calP^{\mbox {\scriptsize !`}}$
in the category of wheeled \PROP{s}.

\begin{proposition}
A wheeled quadratic operad
$\calP$ is wheeled Koszul if and only if $\calP^!$ is wheeled
Koszul. The operadic part $\calP_o$ of a
wheeled Koszul operad $\calP$ is Koszul in the
ordinary sense.
\end{proposition}

\begin{proof}
The first statement is based on the isomorphism of the
linear duals
\[
(\wB(\calP),\pa_B)^\# \cong (\wO(\calP^\#),\pa_\Omega)
\]
and on the fact that, under the above isomorphism, the linear
dual of the inclusion
\[
\iota:
(\calP^{\mbox {\scriptsize \rm !`}},\pa =0)
\hookrightarrow (\wB(\calP),\pa_B)
\]
is the projection
\[
p:
(\wO(\calP),\pa_\Omega)\to  (\calP^!,\pa = 0).
\]
Since the linear dual is an exact functor, the Koszulity of
$\calP$ implies, by Proposition~\ref{Jitka_nejak_neroste},
the Koszulity of $\calP^!$. The opposite implication follows from
the above statement applied to
$\calP^!$ and the involution property $(\calP^!)^! \cong \calP$.

To prove the second part of the proposition, one needs to observe
that the operadic part of the
map~(\ref{jsem_zvedav_jestli_mne_zamknou}) is the canonical inclusion
\[
\iota_o:
(\calP_o^{\mbox {\scriptsize \rm !`}},\pa =0)
\hookrightarrow (B(\calP_o),\pa_B),
\]
of the co-operad $\calP_o^{\mbox {\scriptsize \rm !`}}$ into its
ordinary bar construction
$(B(\calP_o),\pa_B)$,
and then invoke the definition of Koszulness in the form given
in~\cite[Definition~2.23]{getzler-jones:preprint}.
\end{proof}

\begin{example}
{\rm
\label{misa}
Wheeled operads
$\Com^\wc$, $\Ass^\wc$ and $\Lie^\wc$ are wheeled quadratic Koszul. Indeed,
wheeled Koszulness of $\Com^\wc$ is
Theorem~\ref{thm-Com-wheeled-Koszul} of Section~\ref{Terezka}. The
Koszulity of $\Ass^\wc$ follows from comparing the minimal model of
$\Ass^\wc_\infty$ described in Theorem~A to the bar construction of
$(\Ass^\wc)^!$ given in Example~\ref{co_mne_ceka_doma?}.
Finally,
wheeled Koszulness of $\Lie$ was proved in~\cite{me3}.

We do not know  whether there are (ordinary) quadratic
Koszul operads whose wheeled completion is not wheeled Koszul.
}
\end{example}

\begin{definition}
\label{snad_to_nebude_nic_vazneho}
Let $\calP$ be an ordinary quadratic Koszul operad. We say that $\calP$ is
{\em stably Koszul\/} if the completion
$(\Omega(\calP^!),\partial_\Omega)^\wc \to (\calP^\wc,\partial = 0)$
of the
canonical homology isomorphism
$(\Omega(\calP^!),\partial_\Omega) \to (\calP,\partial = 0)$
is a homology isomorphism, too.
\end{definition}

Equivalently, a quadratic Koszul operad
$\calP$ is stably Koszul if the wheeled
completion $\calP^\wc_\infty$ of its minimal model $\calP_\infty$ is a
minimal model of the wheeled completion $\calP^\wc$. The
proof of the following proposition is a simple exercise.

\begin{proposition}
\label{oslicek}
Let $\calP$ be a quadratic operad such that $\calP^\wc$ is wheeled
Koszul. Then $\calP$
is stably Koszul if and only if
$(\calP^\wc)^!_w = 0$.
\end{proposition}

Proposition~\ref{oslicek}
together with Example~\ref{misa} and
Examples~\ref{co_mne_ceka_doma?}--\ref{kasparek} imply
that the operad $\Lie$ is stably Koszul while the operads $\Ass$ and
$\Com$ are not.


\bip

\section{Wheeled resolution of  $\Ass$ and its applications}

\sip

\subsection{Wheeled completions of $\Ass$ and $\Ass_\infty$}
As the
natural functor ${\tt Proper}^\circlearrowright\rar {\tt
PROP}^\circlearrowright$ is exact, it is enough to understand
homotopy type of
wheeled
completions of $\Ass$ and $\Ass_\infty$ in the category of wheeled
properads. The generalization of all our results in this section to
wheeled PROPs is immediate.  Thus we work from now in the categories
${\tt Proper}$ and ${\tt Proper}^\circlearrowright$.

The operad, {\sf Ass},  of associative algebras is defined as the quotient,
\[
\Ass= \Gamma_{op}\langle A \rangle/(R).
\]
of the free operad, $\freeOP (A)$,
generated by an
$\bS$-module $A=\{A(n)\}$,
$$
A(n):= \left\{ \Ba{ll} k[\bS_2]=\mbox{span}\ \left(
\ \
\begin{xy}
<0mm,0mm>*{\bullet},
<0mm,5mm>*{}**@{-},
<-3mm,-5mm>*{}**@{-},
<3mm,-5mm>*{}**@{-},
<-4mm,-7mm>*{_{\sigma(1)}},
<4mm,-7mm>*{_{\sigma(2)}},
\end{xy}
\
\right)_{\sigma\in \Sigma_2} & \mbox{for}\ n=2 \vspace{3mm}\\

0 & \mbox{otherwise},
\Ea
\right.
$$
modulo the ideal  generated by
relations,
\Beq\label{assrel}
R:\ \ \
\begin{xy}
<0mm,0mm>*{\bullet},
<0mm,5mm>*{}**@{-},
<-3mm,-5mm>*{}**@{-},
<3mm,-5mm>*{}**@{-},
<-3mm,-5mm>*{\bullet};
<-6mm,-10mm>*{}**@{-},
<0mm,-10mm>*{}**@{-},
<-6mm,-12mm>*{_{\sigma(1)}},
<0mm,-12mm>*{_{\sigma(2)}},
<3mm,-7mm>*{_{\sigma(3)}},
\end{xy}
-
\begin{xy}
<0mm,0mm>*{\bullet},
<0mm,5mm>*{}**@{-},
<-3mm,-5mm>*{}**@{-},
<3mm,-5mm>*{}**@{-},
<3mm,-5mm>*{\bullet};
<6mm,-10mm>*{}**@{-},
<0mm,-10mm>*{}**@{-},
<6mm,-12mm>*{_{\sigma(3)}},
<0mm,-12mm>*{_{\sigma(2)}},
<-3mm,-7mm>*{_{\sigma(1)}},
\end{xy}
=0,  \ \ \ \ \forall \sigma\in \bS_3.
\Eeq

It is well-known \cite{St}
that the minimal resolution of $\Ass$ in the
category of ordinary
properads
is the dg free properad, ${\Ass}_\infty:={\Gamma}_{op}\langle E \rangle$,
generated by the
$\Sigma$-bimodule
$E=\{E(m,n)\}$,
$$
E(m,n)=\left\{
\Ba{cr}
\id_1\ot \downarrow^{n-2}\,
\bfk[\Sigma_n]= \mbox{span}
\left(
\begin{xy}
 <0mm,0mm>*{\bullet};<0mm,0mm>*{}**@{},
 <0mm,0mm>*{};<-8mm,-5mm>*{}**@{-},
 <0mm,0mm>*{};<-4.5mm,-5mm>*{}**@{-},
 <0mm,0mm>*{};<0mm,-4mm>*{\ldots}**@{},
 <0mm,0mm>*{};<4.5mm,-5mm>*{}**@{-},
 <0mm,0mm>*{};<8mm,-5mm>*{}**@{-},
   <0mm,0mm>*{};<-11mm,-7.9mm>*{^{\sigma(1)}}**@{},
   <0mm,0mm>*{};<-4mm,-7.9mm>*{^{\sigma(2)}}**@{},
   <0mm,0mm>*{};<10.0mm,-7.9mm>*{^{\sigma(n)}}**@{},
 <0mm,0mm>*{};<0mm,5mm>*{}**@{-},
 \end{xy}\ \
\right)_{\sigma\in \bS_n}
 & \mbox{for}\ m=1, n\geq 2 \vspace{6mm}\\
0 & \mbox{otherwise},
\Ea
\right.
$$
and equipped with the differential given on the generators as

\begin{equation} \label{ass-diff}
\p
\begin{xy}
<0mm,0mm>*{\bullet},
<0mm,5mm>*{}**@{-},
<-5mm,-5mm>*{}**@{-},
<-2mm,-5mm>*{}**@{-},
<2mm,-5mm>*{}**@{-},
<5mm,-5mm>*{}**@{-},
<0mm,-7mm>*{_{\sigma(1)\ \ \ \ldots\ \ \ \sigma(n)}},
\end{xy}
=\sum_{k=0}^{n-2}\sum_{l=2}^{n-k}
(-1)^{k+l(n-k-l)+1}
\begin{xy}
<0mm,0mm>*{\bullet},
<0mm,5mm>*{}**@{-},
<4mm,-7mm>*{^{\sigma(1)\dots\sigma(k)\qquad\sigma(k+l+1)\dots\sigma(n)}},
<-14mm,-5mm>*{}**@{-},
<-6mm,-5mm>*{}**@{-},
<20mm,-5mm>*{}**@{-},
<8mm,-5mm>*{}**@{-},
<0mm,-5mm>*{}**@{-},
<0mm,-5mm>*{\bullet};
<-5mm,-10mm>*{}**@{-},
<-2mm,-10mm>*{}**@{-},
<2mm,-10mm>*{}**@{-},
<5mm,-10mm>*{}**@{-},
<0mm,-12mm>*{_{\sigma(k+1)\dots\sigma(k+l)}},
\end{xy}.
\end{equation}

\noindent A natural morphism of dg properads,
$$
p: ({\Ass}_\infty, \p) \lon ({\Ass}, 0)
$$
given on generators by
$$
p\left(
\begin{xy}
 <0mm,0mm>*{\bullet};<0mm,0mm>*{}**@{},
 <0mm,0mm>*{};<-8mm,-5mm>*{}**@{-},
 <0mm,0mm>*{};<-4.5mm,-5mm>*{}**@{-},
 <0mm,0mm>*{};<0mm,-4mm>*{\ldots}**@{},
 <0mm,0mm>*{};<4.5mm,-5mm>*{}**@{-},
 <0mm,0mm>*{};<8mm,-5mm>*{}**@{-},
   <0mm,0mm>*{};<-11mm,-7.9mm>*{^{\sigma(1)}}**@{},
   <0mm,0mm>*{};<-4mm,-7.9mm>*{^{\sigma(2)}}**@{},
   <0mm,0mm>*{};<10.0mm,-7.9mm>*{^{\sigma(n)}}**@{},
 <0mm,0mm>*{};<0mm,5mm>*{}**@{-},
 \end{xy}\ \
\right):=
\left\{
\Ba{cr}
\Id
 & \mbox{for}\ n=2 \vspace{6mm}\\
0 & \mbox{otherwise}
\Ea
\right.
$$
is a quasi-isomorphism. The morphism $p$ induces an associated morphism of dg wheeled properads,
$$
p^\circlearrowright: \left(({\Ass}_\infty)^\circlearrowright, \p\right) \lon
\left({\Ass^\circlearrowright}, 0\right),
$$
which, however, can {\em not}\, be a
quasi-isomorphism.
Indeed \cite{me3}, while
$$
p\left(
\begin{xy}
 <0mm,-0.55mm>*{};<0mm,3mm>*{}**@{-},
 <0mm,0.5mm>*{};<4.2mm,-4mm>*{}**@{-},
 <0mm,0.6mm>*{};<0mm,-6mm>*{}**@{-},
 <-0mm,0.48mm>*{};<-4.2mm,-4mm>*{}**@{-},
 <0mm,0mm>*{\bullet};<0mm,0mm>*{}**@{},
   <0.5mm,0.5mm>*{};<4.7mm,-6.5mm>*{^2}**@{},
   <-0.48mm,0.48mm>*{};<-4.7mm,-6.5mm>*{^1}**@{},
(0,3)*{}
   \ar@{->}@(ur,dr) (0,-6)*{}
 \end{xy}
\right)=0,
 $$
the definition of $\p$ implies
\Beqrn
\p
\begin{xy}
 <0mm,-0.55mm>*{};<0mm,3mm>*{}**@{-},
 <0mm,0.5mm>*{};<4.2mm,-4mm>*{}**@{-},
 <0mm,0.6mm>*{};<0mm,-6mm>*{}**@{-},
 <-0mm,0.48mm>*{};<-4.2mm,-4mm>*{}**@{-},
 <0mm,0mm>*{\bullet};<0mm,0mm>*{}**@{},
   <0.5mm,0.5mm>*{};<4.7mm,-6.5mm>*{^2}**@{},
   <-0.48mm,0.48mm>*{};<-4.7mm,-6.5mm>*{^1}**@{},
(0,3)*{}
   \ar@{->}@(ur,dr) (0,-6)*{}
 \end{xy}
 =  -
\begin{xy}
 <0mm,3mm>*{\bullet};<0mm,0mm>*{}**@{},
 <0mm,3mm>*{};<0mm,7mm>*{}**@{-},
 <0mm,3mm>*{};<3mm,-1mm>*{}**@{-},
 <-0mm,3mm>*{};<-3mm,-1mm>*{}**@{-},
 <-3mm,-1mm>*{\bullet};<-2.3mm,2.3mm>*{}**@{},
 <-3mm,-1mm>*{};<-6mm,-5mm>*{}**@{-},
 <-3mm,-1mm>*{};<0mm,-5mm>*{}**@{-},
   <0.49mm,0.49mm>*{};<3mm,-3.5mm>*{^{2}}**@{},
   <-3mm,-1mm>*{};<-6mm,-7.5mm>*{^{1}}**@{},
   (0,7)*{}
   \ar@{->}@(ur,dr) (0,-5)*{}
 \end{xy}
 \ + \
\begin{xy}
 <0mm,3mm>*{\bullet};<0mm,0mm>*{}**@{},
 <0mm,3mm>*{};<0mm,7mm>*{}**@{-},
 <0mm,3mm>*{};<3mm,-1mm>*{}**@{-},
 <-0mm,3mm>*{};<-3mm,-1mm>*{}**@{-},
 <3mm,-1mm>*{\bullet};<-2.3mm,2.3mm>*{}**@{},
 <3mm,-1mm>*{};<0mm,-5mm>*{}**@{-},
 <3mm,-1mm>*{};<6mm,-5mm>*{}**@{-},
   <0.49mm,0.49mm>*{};<6.2mm,-7.5mm>*{^{2}}**@{},
   <-3mm,-1mm>*{};<-3mm,-3.5mm>*{^{1}}**@{},
   (0,7)*{}
   \ar@{->}@(ul,dl) (0,-5)*{}
 \end{xy}
  = -
 \begin{xy}
 <0mm,3mm>*{\bullet};<0mm,0mm>*{}**@{},
 <0mm,3mm>*{};<0mm,7mm>*{}**@{-},
 <0mm,3mm>*{};<3mm,-1mm>*{}**@{-},
 <-0mm,3mm>*{};<-3mm,-1mm>*{}**@{-},
 <-3mm,-1mm>*{\bullet};<-2.3mm,2.3mm>*{}**@{},
 <-3mm,-1mm>*{};<-6mm,-5mm>*{}**@{-},
 <-3mm,-1mm>*{};<0mm,-5mm>*{}**@{-},
   <0.49mm,0.49mm>*{};<3mm,-3.5mm>*{^{2}}**@{},
   <-3mm,-1mm>*{};<-6mm,-7.5mm>*{^{1}}**@{},
   (0,7)*{}
   \ar@{->}@(ur,dr) (0,-5)*{}
 \end{xy}
 \ + \
\begin{xy}
 <0mm,3mm>*{\bullet};<0mm,0mm>*{}**@{},
 <0mm,3mm>*{};<0mm,7mm>*{}**@{-},
 <0mm,3mm>*{};<3mm,-1mm>*{}**@{-},
 <-0mm,3mm>*{};<-3mm,-1mm>*{}**@{-},
 <-3mm,-1mm>*{\bullet};<-2.3mm,2.3mm>*{}**@{},
 <-3mm,-1mm>*{};<-6mm,-5mm>*{}**@{-},
 <-3mm,-1mm>*{};<0mm,-5mm>*{}**@{-},
   <0.49mm,0.49mm>*{};<3mm,-3.5mm>*{^{2}}**@{},
   <-3mm,-1mm>*{};<-6mm,-7.5mm>*{^{1}}**@{},
   (0,7)*{}
   \ar@{->}@(ur,dr) (0,-5)*{}
 \end{xy}
 = 0.
\Eeqrn
We therefore have non-trivial cohomology classes,
$$
\begin{xy}
<0mm,0mm>*{\bullet},
<-3mm,-5mm>*{}**@{-},
<3mm,-5mm>*{}**@{-},
<-4mm,-7mm>*{_{\sigma(1)}},
<4mm,-7mm>*{_{\sigma(2)}},
\end{xy} := \ \left[ \rule{0pt}{2.6em}\right.
\begin{xy}
 <0mm,-0mm>*{};<3mm,4mm>*{}**@{-},
 <0mm,0.5mm>*{};<4.2mm,-4mm>*{}**@{-},
 <0mm,0.6mm>*{};<0mm,-8mm>*{}**@{-},
 <-0mm,0.48mm>*{};<-4.2mm,-4mm>*{}**@{-},
 <0mm,0mm>*{\bullet};<0mm,0mm>*{}**@{},
   <0.5mm,0.5mm>*{};<4.7mm,-6.5mm>*{^{\sigma(2)}}**@{},
   <-0.48mm,0.48mm>*{};<-4.7mm,-6.5mm>*{^{\sigma(1)}}**@{},
(3,4)*{}
   \ar@{->}@(ur,dr) (0,-8)*{}
 \end{xy} \left. \rule{0pt}{2.6em} \right] \  \ \ \ \ \ \sigma\in\Sigma_2,
$$
in  $H^{-1}\left({(\Ass}_\infty)^\circlearrowright,\p\right)$
belonging to the kernel of $H^*(p^\wc)$.

\begin{theorem}
The cohomology group
$H^\bullet\left((\Ass_\infty)^\circlearrowright,\p\right)$
is concentrated in degrees $0$ and $-1$.
Moreover,
$$
H^0\left((\Ass_\infty)^\circlearrowright,\p\right)= \Ass^\circlearrowright,
$$
and, as $\bS$-bimodules,
$$
H^{-1}\left((\Ass_\infty)^\circlearrowright,\p\right)(m,n)=\left\{
\Ba{ll}
\bigoplus_{p=1}^{n-1}k[\Sigma_n]^{C_p\times C_{n-p}} & \mbox{for}\ m=0, n\geq 2\\
0 & \mbox{otherwise}\\
\Ea
\right.
$$
where $C_p\times C_{n-p}$  is the  subgroup of\, $\Sigma_n$ generated by two commuting
cyclic permutations $(12\ldots p)$ and $(p+1\ldots n)$, and
$k[\Sigma_n]^{C_p\times C_{n-p}}$ stands for coinvariants.
\end{theorem}

\begin{proof}
The space $(\Ass_\infty)^\circlearrowright$  is spanned by graphs of genus 0 and 1,
and the differential $\p$ preserves the associated genus decomposition,
$$
(\Ass_\infty)^\circlearrowright= \fA^\uparrow \oplus \fA^\circlearrowright.
$$
The subcomplex $(\fA^\uparrow,\p)$ is spanned, by definition, by graphs of genus zero and hence is isomorphic
to $(\Ass_\infty,\p)$ so that $H(\fA^\uparrow,\p)=\Ass$. Thus the main job is to compute the cohomology,
$H(\fA^\circlearrowright,\p)$,
of the subcomplex $\fA^\circlearrowright$ spanned by graphs of genus one,
i.e.\ by graphs of the form
$$
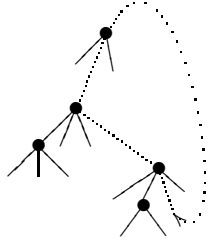
\begin{xy}
<0mm,0mm>*{\bullet},
<4mm,10mm>*{}**@{.},
<-5mm,-5mm>*{}**@{-},
<-2mm,-5mm>*{}**@{-},
<2mm,-5mm>*{}**@{-},
%
<-5mm,-5mm>*{\bullet},
<-5mm,-5mm>*{};<-1mm,-9mm>*{}**@{-},
<-5mm,-5mm>*{};<-5mm,-9mm>*{}**@{-},
<-5mm,-5mm>*{};<-9mm,-9mm>*{}**@{-},
<4mm,10mm>*{\bullet};
<4mm,10mm>*{};<5mm,5mm>*{}**@{-},
<-0mm,6mm>*{}**@{-},
<4mm,10mm>*{};<6mm,13mm>*{}**@{.},
<11mm,-8mm>*{\bullet};
<11mm,-8mm>*{};<0mm,0mm>*{}**@{.},
<11mm,-8mm>*{};<9mm,-12mm>*{}**@{-},
<11mm,-8mm>*{};<5mm,-12mm>*{}**@{-},
<11mm,-8mm>*{};<13mm,-14mm>*{}**@{.},
<11mm,-8mm>*{};<14.5mm,-11mm>*{}**@{-},
<9mm,-12.9mm>*{\bullet};
<9mm,-12.9mm>*{};<6mm,-16.9mm>*{}**@{-},
<9mm,-12.9mm>*{};<12mm,-16.9mm>*{}**@{-},
(6,13)*{}
   \ar@{.>}@(ur,dr) (13,-14)*{}
\end{xy}\ \
$$
with internal edges lying in the wheel dotted.

The following terminology will be useful: the vertices of an element
$G\in \fA^\circlearrowright$ which lie on the wheel are called {\em cyclic}.
For example, the graph shown above has three cyclic vertices and two noncyclic vertices.

It is clear that
$$
 F_p\fA^\circlearrowright:= {\rm span}\left\{ f\in \fA^\circlearrowright:\
 \mbox{total number of cyclic vertices}
   \geq p\right\},
 $$
defines a filtration in the complex $(\fA^\circlearrowright, \p)$. Let $\{\fA_p^\circlearrowright, \p_p\}_{p\geq 0}$
be the associated spectral sequence.

{\bf Step 1}. Our first target is  to compute the cohomology, $\fA_1^\circlearrowright$, of the zeroth term,
 $\{\fA_0^\circlearrowright, \p_0\}$,
of this spectral sequence. Consider the filtration,
$$
 \cF_p\fA_0^\circlearrowright:= {\rm span}\left\{ G\in \fA_0^\circlearrowright:\
 \Ba{l}
 \mbox{total number of internal edges and legs
  }\\
  \mbox{attached to cyclic vertices of}\ G\Ea
   \leq p\right\},
 $$
and let $\{\cE_r\fA_0^\circlearrowright, \delta_r\}_{r\geq 0}$ be the associated spectral
sequence. We shall show below
that  the latter degenerates at the
second term so that $\cE_2\fA_0^\circlearrowright=
H({\mathsf \fA_0^\circlearrowright}, \p_0)=\fA_1^\circlearrowright$.

 The
differential $\delta_0$ in $\cE_0\fA_0^\circlearrowright$ is given by its values on the
vertices as follows:
\Bi
\item[(i)]  on every noncyclic vertex one has $\delta_0=\p$, the differential in
$\Ass_\infty$;
\item[(ii)] on every cyclic vertex  $\delta_0=0$.
\Ei
Hence, modulo the action of finite groups, the complex
 $(\cE_0\fA_0^\circlearrowright, \delta_0)$ is isomorphic to the {direct} sum of tensor products of
copies of the complex
$( {\Ass_\infty} ,\p)$. By
 K\"{u}nneth's and Mashke's
theorems, we get,
$$
\cE_1\fA_0^\circlearrowright= V_1/h(V_2),
$$
where
\Bi
\item[-] $V_1$ is the subspace of $\fA^\circlearrowright$
consisting of all those graphs whose every noncyclic vertex is
$
\begin{xy}
 <0mm,0.66mm>*{};<0mm,3mm>*{}**@{-},
 <0.39mm,-0.39mm>*{};<2.2mm,-2.2mm>*{}**@{-},
 <-0.35mm,-0.35mm>*{};<-2.2mm,-2.2mm>*{}**@{-},
 <0mm,0mm>*{\bullet};<0mm,0mm>*{}**@{},
\end{xy} ;
$
\item[-] $V_2$ is the subspace of  $\fA^\circlearrowright$
 whose
every noncyclic vertex is either
$
\begin{xy}
 <0mm,0.66mm>*{};<0mm,3mm>*{}**@{-},
 <0.39mm,-0.39mm>*{};<2.2mm,-2.2mm>*{}**@{-},
 <-0.35mm,-0.35mm>*{};<-2.2mm,-2.2mm>*{}**@{-},
 <0mm,0mm>*{\bullet};<0mm,0mm>*{}**@{},
\end{xy}
$
 or
$
\begin{xy}
 <0mm,0.66mm>*{};<0mm,3mm>*{}**@{-},
 <0.39mm,-0.39mm>*{};<2.2mm,-2.2mm>*{}**@{-},
 <-0.35mm,-0.35mm>*{};<-2.2mm,-2.2mm>*{}**@{-},
 <0mm,0mm>*{\bullet};<0mm,0mm>*{}**@{},
 <0mm,-0.66mm>*{};<0mm,-2.2mm>*{}**@{-},
\end{xy}
$
with the number of vertices of the latter type $\geq 1$
;
\item[-] the map $h:V_2\rar V_1$ is $\Sigma$-equivariant and is given on noncyclic vertices by
$$
h\left(
\begin{xy}
 <0mm,0.66mm>*{};<0mm,3mm>*{}**@{-},
 <0.39mm,-0.39mm>*{};<2.2mm,-2.2mm>*{}**@{-},
 <-0.35mm,-0.35mm>*{};<-2.2mm,-2.2mm>*{}**@{-},
 <0mm,0mm>*{\bullet};<0mm,0mm>*{}**@{},
\end{xy}
\right) =
\begin{xy}
 <0mm,0.66mm>*{};<0mm,3mm>*{}**@{-},
 <0.39mm,-0.39mm>*{};<2.2mm,-2.2mm>*{}**@{-},
 <-0.35mm,-0.35mm>*{};<-2.2mm,-2.2mm>*{}**@{-},
 <0mm,0mm>*{\bullet};<0mm,0mm>*{}**@{},
\end{xy} \ \
, \ \ \ \ \ \ \ \
h\left(
\begin{xy}
 <0mm,0.66mm>*{};<0mm,3mm>*{}**@{-},
 <0.39mm,-0.39mm>*{};<2.2mm,-2.2mm>*{}**@{-},
 <-0.35mm,-0.35mm>*{};<-2.2mm,-2.2mm>*{}**@{-},
 <0mm,0mm>*{\bullet};<0mm,0mm>*{}**@{},
 <0mm,-0.66mm>*{};<0mm,-2.2mm>*{}**@{-},
 <0mm,0.66mm>*{};<0mm,-4mm>*{^2}**@{},
   <0.39mm,-0.39mm>*{};<2.9mm,-4mm>*{^3}**@{},
   <-0.35mm,-0.35mm>*{};<-2.8mm,-4mm>*{^1}**@{},
\end{xy}
\right) = -
 \begin{xy}
 <0mm,0mm>*{\bullet};<0mm,0mm>*{}**@{},
 <0mm,0.69mm>*{};<0mm,3.0mm>*{}**@{-},
 <0.39mm,-0.39mm>*{};<2.4mm,-2.4mm>*{}**@{-},
 <-0.35mm,-0.35mm>*{};<-1.9mm,-1.9mm>*{}**@{-},
 <-2.4mm,-2.4mm>*{\bullet};<-2.4mm,-2.4mm>*{}**@{},
 <-2.0mm,-2.8mm>*{};<0mm,-4.9mm>*{}**@{-},
 <-2.8mm,-2.9mm>*{};<-4.7mm,-4.9mm>*{}**@{-},
    <0.39mm,-0.39mm>*{};<3.3mm,-4.0mm>*{^3}**@{},
    <-2.0mm,-2.8mm>*{};<0.5mm,-6.7mm>*{^2}**@{},
    <-2.8mm,-2.9mm>*{};<-5.2mm,-6.7mm>*{^1}**@{},
 \end{xy}
\ + \
 \begin{xy}
 <0mm,0mm>*{\bullet};<0mm,0mm>*{}**@{},
 <0mm,0.69mm>*{};<0mm,3.0mm>*{}**@{-},
 <0.39mm,-0.39mm>*{};<2.4mm,-2.4mm>*{}**@{-},
 <-0.35mm,-0.35mm>*{};<-1.9mm,-1.9mm>*{}**@{-},
 <2.4mm,-2.4mm>*{\bullet};<-2.4mm,-2.4mm>*{}**@{},
 <2.4mm,-2.4mm>*{};<0mm,-4.9mm>*{}**@{-},
 <2.4mm,-2.4mm>*{};<4.7mm,-4.9mm>*{}**@{-},
    <0.39mm,-0.39mm>*{};<5.8mm,-6.9mm>*{^3}**@{},
    <-2.0mm,-2.8mm>*{};<-2.5mm,-3.9mm>*{^1}**@{},
    <-2.8mm,-2.9mm>*{};<-0.2mm,-6.9mm>*{^2}**@{},
 \end{xy}
 $$
and  on all  cyclic vertices $h$ is set to be the identity.
\Ei

A representative of a
typical element in  $\cE_1\fA^\circlearrowright$ looks as
$$
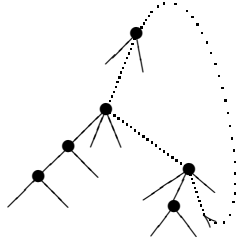
\begin{xy}
<0mm,0mm>*{\bullet},
<4mm,10mm>*{}**@{.},
<-5mm,-5mm>*{}**@{-},
<-2mm,-5mm>*{}**@{-},
<2mm,-5mm>*{}**@{-},
%
<-5mm,-5mm>*{\bullet},
<-5mm,-5mm>*{};<-1mm,-9mm>*{}**@{-},
<-5mm,-5mm>*{};<-13mm,-13mm>*{}**@{-},
<-9mm,-9mm>*{\bullet},
<-9mm,-9mm>*{};<-5mm,-13mm>*{}**@{-},
<4mm,10mm>*{\bullet};
<4mm,10mm>*{};<5mm,5mm>*{}**@{-},
<-0mm,6mm>*{}**@{-},
<4mm,10mm>*{};<6mm,13mm>*{}**@{.},
<11mm,-8mm>*{\bullet};
<11mm,-8mm>*{};<0mm,0mm>*{}**@{.},
<11mm,-8mm>*{};<9mm,-12mm>*{}**@{-},
<11mm,-8mm>*{};<5mm,-12mm>*{}**@{-},
<11mm,-8mm>*{};<13mm,-14mm>*{}**@{.},
<11mm,-8mm>*{};<14.5mm,-11mm>*{}**@{-},
<9mm,-12.9mm>*{\bullet};
<9mm,-12.9mm>*{};<6mm,-16.9mm>*{}**@{-},
<9mm,-12.9mm>*{};<12mm,-16.9mm>*{}**@{-},
(6,13)*{}
   \ar@{.>}@(ur,dr) (13,-14)*{}
\end{xy}\ \
$$
The differential
$\delta_1$ in $\cE_1\fA^\circlearrowright$ is given by its values on
vertices as
\Bi
\item[(i)]  on every noncyclic vertex one has $\delta_1=0$;
\item[(ii)] on every cyclic vertex
one has
\Beqr\label{delta1}
\delta_1
\begin{xy}
<0mm,0mm>*{\bullet},
<0mm,-6mm>*{}**@{.},
<0mm,6mm>*{}**@{.},
<-16mm,-4mm>*{}**@{-},
<-12mm,-4mm>*{}**@{-},
<-4mm,-4mm>*{}**@{-},
<-16.5mm,-6mm>*{_1},
<-12.5mm,-6mm>*{_2},
<-8.5mm,-6mm>*{_{..}},
<-4.5mm,-6mm>*{_p},
<16mm,-4mm>*{}**@{-},
<12mm,-4mm>*{}**@{-},
<4mm,-4mm>*{}**@{-},
<16.9mm,-6mm>*{_{n}},
<11.5mm,-6mm>*{_{...}},
<4.5mm,-6mm>*{_{p+1}},
\end{xy}
&=&
\sum_{i=0}^{p-2}(-1)^{i+1}
\begin{xy}
<0mm,0mm>*{\bullet},
<0mm,-6mm>*{}**@{.},
<0mm,6mm>*{}**@{.},
<-20mm,-4mm>*{}**@{-},
<-15mm,-4mm>*{}**@{-},
<-11mm,-4mm>*{}**@{-},
<-4mm,-4mm>*{}**@{-},
<-21mm,-6mm>*{_1},
<-16mm,-6mm>*{_i},
<-18.5mm,-6mm>*{_{..}},
<-8mm,-6mm>*{_{..}},
<-4.5mm,-6mm>*{_p},
<16mm,-4mm>*{}**@{-},
<12mm,-4mm>*{}**@{-},
<4mm,-4mm>*{}**@{-},
<16.9mm,-6mm>*{_{n}},
<11.5mm,-6mm>*{_{...}},
<4.5mm,-6mm>*{_{p+1}},
<-11mm,-4.9mm>*{\bullet};
<-11mm,-4.9mm>*{};<-8mm,-8.9mm>*{}**@{-},
<-11mm,-4.9mm>*{};<-14mm,-8.9mm>*{}**@{-},
<-6.5mm,-10.5mm>*{_{i+2}},
<-14.5mm,-10.5mm>*{_{i+1}},
\end{xy}  \\
&&+
\sum_{i=p+1}^{n-2}(-1)^{i+1}
\begin{xy}
<0mm,0mm>*{\bullet},
<0mm,-6mm>*{}**@{.},
<0mm,6mm>*{}**@{.},
<-16mm,-4mm>*{}**@{-},
<-12mm,-4mm>*{}**@{-},
<-4mm,-4mm>*{}**@{-},
<-16.5mm,-6mm>*{_1},
<-12.5mm,-6mm>*{_2},
<-8.5mm,-6mm>*{_{..}},
<-4.5mm,-6mm>*{_p},
<22mm,-4mm>*{}**@{-},
<10mm,-4mm>*{}**@{-},
<15.5mm,-4.5mm>*{}**@{-},
<4mm,-4mm>*{}**@{-},
<10.1mm,-5.6mm>*{_{i}},
<23mm,-6mm>*{_{n}},
<7.5mm,-6mm>*{_{..}},
<3.6mm,-6mm>*{_{p+1}},
<15mm,-5mm>*{\bullet};
<19mm,-6mm>*{_{..}},
<15mm,-5mm>*{};<12mm,-9mm>*{}**@{-},
<15mm,-5mm>*{};<18mm,-9mm>*{}**@{-},
<11.5mm,-10.5mm>*{_{i+1}},
<18.5mm,-10.5mm>*{_{i+2}},
\end{xy}
\nonumber
\Eeqr
\Ei
To compute $\cE_2\fA^\circlearrowright=H(\cE_1\fA^\circlearrowright,\delta_1)$ let us
return back to the well-known complex $\Ass_\infty$: the data
$$
F_p{\Ass_\infty} := {\rm span} \left\{G \in {\mathsf \Ass_\infty}:
\mbox{number of edges attached to the
root vertex of}\ G \leq p\right\}
$$
is clearly a filtration of the  complex $({\Ass_\infty}, \delta)$ whose spectral sequence,
$\{E_r\Ass_\infty,d_r\}_{r\geq 0}$, must converge to $\Ass$. Its first term,
$E_1\Ass_\infty=H(E_0\Ass_\infty,d_0)$ is spanned by trees
 whose root vertex may have any number of attached half-edges while all other vertices are
binary,
$
\begin{xy}
 <0mm,0.66mm>*{};<0mm,3mm>*{}**@{-},
 <0.39mm,-0.39mm>*{};<2.2mm,-2.2mm>*{}**@{-},
 <-0.35mm,-0.35mm>*{};<-2.2mm,-2.2mm>*{}**@{-},
 <0mm,0mm>*{\bullet};<0mm,0mm>*{}**@{},
\end{xy}
$. The differential $d_1$ is non-trivial only on the root vertex on
which it is given by,
$$
d_1\begin{xy}
 <0mm,0mm>*{\bullet};<0mm,0mm>*{}**@{},
 <0mm,0mm>*{};<-8mm,-5mm>*{}**@{-},
 <0mm,0mm>*{};<-4.5mm,-5mm>*{}**@{-},
 <0mm,0mm>*{};<0mm,-4mm>*{\ldots}**@{},
 <0mm,0mm>*{};<4.5mm,-5mm>*{}**@{-},
 <0mm,0mm>*{};<8mm,-5mm>*{}**@{-},
   <0mm,0mm>*{};<-11mm,-7.9mm>*{^{1}}**@{},
   <0mm,0mm>*{};<-4mm,-7.9mm>*{^{2}}**@{},
   <0mm,0mm>*{};<10.0mm,-7.9mm>*{^n}**@{},
 <0mm,0mm>*{};<0mm,5mm>*{}**@{-},
 \end{xy}=
 \sum_{i=0}^{n-2}(-1)^{i+1}
\begin{xy}
 <0mm,0mm>*{\bullet};<0mm,0mm>*{}**@{},
 <0mm,0mm>*{};<-8mm,-5mm>*{}**@{-},
 <0mm,0mm>*{};<-4mm,-5mm>*{}**@{-},
 <0mm,0mm>*{};<-5.6mm,-5mm>*{..}**@{},
 <0mm,0mm>*{};<2.6mm,-5mm>*{..}**@{},
 <0mm,0mm>*{};<4.5mm,-5mm>*{}**@{-},
 <0mm,0mm>*{};<-0.5mm,-5mm>*{}**@{-},
 <0mm,0mm>*{};<8mm,-5mm>*{}**@{-},
   <0mm,0mm>*{};<-11mm,-7.9mm>*{^1}**@{},
   <0mm,0mm>*{};<-4.5mm,-7.9mm>*{^i}**@{},
   <0mm,0mm>*{};<10.0mm,-7.9mm>*{^n}**@{},
 <0mm,0mm>*{};<0mm,5mm>*{}**@{-},
  <-0.5mm,-5.4mm>*{\bullet};
<-0.7mm,-5.4mm>*{};<-3.6mm,-9mm>*{}**@{-},
<-0.7mm,-5.4mm>*{};<3mm,-9mm>*{}**@{-},
<-0.7mm,-5.4mm>*{};<-4mm,-11mm>*{_{i+1}}**@{},
<-0.7mm,-5.4mm>*{};<4mm,-11mm>*{_{i+2}}**@{},
 \end{xy}
$$
It is clear that this spectral sequence must  degenerate at
$E_2\Ass_\infty=H(E_1\Ass_\infty,d_1)$ implying
the isomorphism $H(E_1\Ass_\infty,d_1)=\Ass$. Now let us modify
the complex $(E_1\Ass_\infty,d_1)$ by adding to the space
$E_1\Ass_\infty$ the trees whose root vertex is a degree $-1$
corolla
$\begin{xy}
 <0mm,-0.55mm>*{};<0mm,-2.5mm>*{}**@{-},
 <0mm,0mm>*{};<0mm,2.5mm>*{}**@{-},
 <0mm,0mm>*{\bullet};<0mm,0mm>*{}**@{},
 \end{xy}$
 while all other vertices are binary
$
\begin{xy}
 <0mm,0.66mm>*{};<0mm,3mm>*{}**@{-},
 <0.39mm,-0.39mm>*{};<2.2mm,-2.2mm>*{}**@{-},
 <-0.35mm,-0.35mm>*{};<-2.2mm,-2.2mm>*{}**@{-},
 <0mm,0mm>*{\bullet};<0mm,0mm>*{}**@{},
\end{xy}
$.

Denote this extension of $E_1\Ass_\infty$ by
$E_1^+\Ass_\infty$, and define a differential $d_1^+$ on $E_1^+\Ass_\infty$
be setting its values on the non-root vertices to be zero while on the root
 $(1,n)$-vertex as follows
\Beq\label{d1+}
d_1^+\begin{xy}
 <0mm,0mm>*{\bullet};<0mm,0mm>*{}**@{},
 <0mm,0mm>*{};<-8mm,-5mm>*{}**@{-},
 <0mm,0mm>*{};<-4.5mm,-5mm>*{}**@{-},
 <0mm,0mm>*{};<0mm,-4mm>*{\ldots}**@{},
 <0mm,0mm>*{};<4.5mm,-5mm>*{}**@{-},
 <0mm,0mm>*{};<8mm,-5mm>*{}**@{-},
   <0mm,0mm>*{};<-11mm,-7.9mm>*{^{1}}**@{},
   <0mm,0mm>*{};<-4mm,-7.9mm>*{^{2}}**@{},
   <0mm,0mm>*{};<10.0mm,-7.9mm>*{^n}**@{},
 <0mm,0mm>*{};<0mm,5mm>*{}**@{-},
 \end{xy}=\left\{\Ba{lr}
 \sum_{i=0}^{n-2}(-1)^{i+1}
\begin{xy}
 <0mm,0mm>*{\bullet};<0mm,0mm>*{}**@{},
 <0mm,0mm>*{};<-8mm,-5mm>*{}**@{-},
 <0mm,0mm>*{};<-4mm,-5mm>*{}**@{-},
 <0mm,0mm>*{};<-5.6mm,-5mm>*{..}**@{},
 <0mm,0mm>*{};<2.6mm,-5mm>*{..}**@{},
 <0mm,0mm>*{};<4.5mm,-5mm>*{}**@{-},
 <0mm,0mm>*{};<-0.5mm,-5mm>*{}**@{-},
 <0mm,0mm>*{};<8mm,-5mm>*{}**@{-},
   <0mm,0mm>*{};<-11mm,-7.9mm>*{^1}**@{},
   <0mm,0mm>*{};<-4.5mm,-7.9mm>*{^i}**@{},
   <0mm,0mm>*{};<10.0mm,-7.9mm>*{^n}**@{},
 <0mm,0mm>*{};<0mm,5mm>*{}**@{-},
  <-0.5mm,-5.4mm>*{\bullet};
<-0.7mm,-5.4mm>*{};<-3.6mm,-9mm>*{}**@{-},
<-0.7mm,-5.4mm>*{};<3mm,-9mm>*{}**@{-},
<-0.7mm,-5.4mm>*{};<-4mm,-11mm>*{_{i+1}}**@{},
<-0.7mm,-5.4mm>*{};<4mm,-11mm>*{_{i+2}}**@{},
 \end{xy} & \mbox{for}\ n\geq 3\\
\begin{xy}
 <0mm,0mm>*{\bullet};<0mm,0mm>*{}**@{},
 <0mm,0mm>*{};<0mm,-5mm>*{}**@{-},
 <0mm,0mm>*{};<0mm,5mm>*{}**@{-},
  <0mm,-5.4mm>*{\bullet};
<-0mm,-5.4mm>*{};<-3mm,-9mm>*{}**@{-},
<-0.mm,-5.4mm>*{};<3mm,-9mm>*{}**@{-},
<-0.7mm,-5.4mm>*{};<-4mm,-11mm>*{_{1}}**@{},
<-0.7mm,-5.4mm>*{};<4mm,-11mm>*{_{2}}**@{},
 \end{xy}
 & \mbox{for}\ n=2 \vspace{3mm}  \\
0 & \mbox{for}\ n=1.\\
 \Ea
 \right.
\Eeq

\noindent{\em Claim. The cohomology of the complex
$(E^+_1{\Ass_\infty}, d_1^+)$ is a
 one dimensional vector space spanned by $\begin{xy}
 <0mm,-0.55mm>*{};<0mm,-2.5mm>*{}**@{-},
 <0mm,0mm>*{};<0mm,2.5mm>*{}**@{-},
 <0mm,0mm>*{\bullet};<0mm,0mm>*{}**@{},
 \end{xy}$\, .}

\noindent
{\em Proof of the claim}. Consider a 2-step filtration, $F_0\subset F_1$ of the complex
$(E_1^+\Ass_\infty, d_1^+)$ by the number of $\begin{xy}
 <0mm,-0.55mm>*{};<0mm,-2.5mm>*{}**@{-},
 <0mm,0mm>*{};<0mm,2.5mm>*{}**@{-},
 <0mm,0mm>*{\bullet};<0mm,0mm>*{}**@{},
 \end{xy}$ . The zero-th term of the associated spectral sequence is isomorphic
to the direct sum of the complexes,
$$
(E_1\Ass_\infty, d_1)\oplus (\desusp E_1\Ass_\infty, d_1)
\oplus (\mbox{span}\langle \begin{xy}
 <0mm,-0.55mm>*{};<0mm,-2.5mm>*{}**@{-},
 <0mm,0mm>*{};<0mm,2.5mm>*{}**@{-},
 <0mm,0mm>*{\bullet};<0mm,0mm>*{}**@{},
 \end{xy}\rangle, 0)
$$
so that the next term of the spectral sequence is
$$
{\Ass}\ \oplus \desusp{\Ass} \oplus \langle \begin{xy}
 <0mm,-0.55mm>*{};<0mm,-2.5mm>*{}**@{-},
 <0mm,0mm>*{};<0mm,2.5mm>*{}**@{-},
 <0mm,0mm>*{\bullet};<0mm,0mm>*{}**@{},
 \end{xy}\rangle
$$
with the differential being zero on $\desusp\Ass \oplus \langle \begin{xy}
 <0mm,-0.55mm>*{};<0mm,-2.5mm>*{}**@{-},
 <0mm,0mm>*{};<0mm,2.5mm>*{}**@{-},
 <0mm,0mm>*{\bullet};<0mm,0mm>*{}**@{},
 \end{xy}\rangle$ and the
natural isomorphism,
$$
\desusp\ : {\Ass}\lon \desusp \Ass
$$
on the remaining summand. Hence the claim follows.

Comparing differentials (\ref{delta1}) and (\ref{d1+}), we see that,
 modulo actions of finite groups, the complex  $(\cE_1\fA^\circlearrowright_0,\delta_1)$
is isomorphic to the tensor product of a trivial complex (i.e.\ one with zero differential)
with the tensor powers
of the complex, $(E_1^+\Ass_\infty,d_1^+)$. Then the above Claim implies that
$\cE_2\fA^\circlearrowright_0=H(\cE_1\fA^\circlearrowright_0,\delta_1)$ is spanned
by the wheeled graphs whose every vertex is cyclic and is either binary or ternary,
i.e.\ by graphs of the form,
\Beq\label{www}
\cE_2\fA^\circlearrowright_0=\mbox{ span}\ \ \
\xy
*\xycircle(10,10){.};
<10mm,0mm>*{\bullet};
<10mm,0mm>*{};<13mm,3mm>*{}**@{-},
<8mm,6mm>*{\bullet};
<8mm,6mm>*{};<4mm,6mm>*{}**@{-},
<3mm,9.5mm>*{\bullet};
<3mm,9.5mm>*{};<0mm,12mm>*{}**@{-},
<-4mm,9mm>*{\bullet};
<-4mm,9mm>*{};<-8mm,9.5mm>*{}**@{-},
<-4mm,9mm>*{};<-5.5mm,5.5mm>*{}**@{-},
<-9mm,3mm>*{\bullet};
<-9mm,3mm>*{};<-7mm,0mm>*{}**@{-},
<-9mm,-4mm>*{\bullet};
<-9mm,-4mm>*{};<-10.7mm,-7mm>*{}**@{-},
<-9mm,-4mm>*{};<-5mm,-5.7mm>*{}**@{-},
<-2mm,-10mm>*{\bullet};
<-2mm,-10mm>*{};<2mm,-12mm>*{}**@{-},
<6mm,-8mm>*{\bullet};
<6mm,-8mm>*{};<9mm,-10mm>*{}**@{-},
\endxy
\Eeq
with legs numbered by integers (not shown).
The induced differential, $\delta_2$, on such graphs is obviously zero,
so that the spectral sequence $(\cE_p\fA^\circlearrowright_0,\delta_p)$ degenerates
giving us an isomorphism
$\fA^\circlearrowright_1=\lim_{p\rar\infty}
(\cE_p\fA^\circlearrowright_0,\delta_p)\simeq \cE_2\fA^\circlearrowright_0$.

\sip

{\bf Step 2.} We have shown above that the first term,
$\fA^\circlearrowright_1$, of the spectral sequence
$\{\fA^\circlearrowright_p, \p_p\}_{p\geq 0}$
 can be identified with the vector space
spanned by wheeled graphs of the form (\ref{www}). The induced differential
$\p_1$ on $\fA^\circlearrowright_1$ can then be described on generators as follows:
$$
\p_1
\begin{xy}
 <0mm,0mm>*{\bullet};
 <0mm,0mm>*{};<0mm,4mm>*{}**@{.},
 <0mm,0mm>*{};<0mm,-4mm>*{}**@{.},
 <0mm,0mm>*{};<-3mm,-4mm>*{}**@{-},
 \end{xy} = 0, \ \ \ \ \
\p_1
\begin{xy}
 <0mm,0mm>*{\bullet};
 <0mm,0mm>*{};<0mm,4mm>*{}**@{.},
 <0mm,0mm>*{};<0mm,-4mm>*{}**@{.},
 <0mm,0mm>*{};<3mm,-4mm>*{}**@{-},
 \end{xy} = 0, \ \ \ \ \
 \p_1
\begin{xy}
 <0mm,0mm>*{\bullet};
 <0mm,0mm>*{};<0mm,4mm>*{}**@{.},
 <0mm,0mm>*{};<0mm,-4mm>*{}**@{.},
 <0mm,0mm>*{};<3mm,-4mm>*{}**@{-},
 <0mm,0mm>*{};<-3mm,-4mm>*{}**@{-},
 \end{xy}=
\begin{xy}
 <0mm,3mm>*{\bullet};
 <0mm,3mm>*{};<0mm,7mm>*{}**@{.},
 <0mm,3mm>*{};<0mm,-3mm>*{}**@{.},
 <0mm,3mm>*{};<-3mm,-1mm>*{}**@{-},
<0mm,-3mm>*{\bullet};
 <0mm,-3mm>*{};<0mm,-7mm>*{}**@{.},
 <0mm,-3mm>*{};<3mm,-7mm>*{}**@{-},
 \end{xy}
 \ - \
\begin{xy}
 <0mm,3mm>*{\bullet};
 <0mm,3mm>*{};<0mm,7mm>*{}**@{.},
 <0mm,3mm>*{};<0mm,-3mm>*{}**@{.},
 <0mm,3mm>*{};<3mm,-1mm>*{}**@{-},
<0mm,-3mm>*{\bullet};
 <0mm,-3mm>*{};<0mm,-7mm>*{}**@{.},
 <0mm,-3mm>*{};<-3mm,-7mm>*{}**@{-},
 \end{xy}
$$
Now the theorem follows immediately from the following

\noindent{\em Claim. The spectral sequence $\{\fA^\circlearrowright_p, \p_p\}_{p\geq 0}$
degenerates at $p=2$ with an isomorphism of $\bS$-bimodules,
$$
\fA^\circlearrowright_2(0,n)=H(\fA^\circlearrowright_1, \p_1)(0,n)=
\left\{\Ba{ll}
\Ass^\circlearrowright(0,1) & \mbox{for}\ n=1 \\
\Ass^\circlearrowright(0,n) \
\oplus \
{\bigoplus_{p=1}^{n-1}\desusp k[\Sigma_n]^{C_p\times C_{n-p}}}
& \mbox{for}\ n\geq 2.
\Ea
\right.
$$

\noindent
Proof of the Claim.}
Let $(\UPQ rs,\pa_1)$ be, for $r +s = n$, the subcomplex of
$(\fA^\circlearrowright_1, \p_1)(0,n)$ spanned by graphs whose legs
entering the cycle from the left have labels $1,\ldots,r$ and legs
merging from the right labels $r+1,\ldots,n$.
Clearly each $\UPQ rs$ is a $(\Sigma_r \times \Sigma_s)$-module and the
$\Sigma_n$-space
$\fA^\circlearrowright_1$ is the sum of induced representations
\begin{equation}
\label{prvni_den_v_Bonnu}
\fA^\circlearrowright_1 \cong
\bigoplus_{r +s = n} {\rm Ind}^{\Sigma_n}_{\Sigma_r \times
\Sigma_s}\ \UPQ rs.
\end{equation}
Our proof of the claim will be based on showing that the cohomology of
$(\UPQ rs,\pa_1)$ is isomorphic
to the (regraded) cohomology
of the moduli space of
configurations
$\Cof rs$ defined in the following paragraph.

First, let $\tM rs$ be the space of configurations of two types of labeled, not
necessarily distinct,
points on the unit cycle $S^1$ -- ``left'' points labeled by $1,\ldots,r$
and ``right'' points labeled $\hbox{$r+1$},\ldots,n$. Let $\tCof rs$ be the
open subspace of configurations such
that points of the same type
do not collide, that is, only ``left-right'' collisions are
allowed in $\tCof rs$. We denote by
$\MM rs$ and  $\Cof rs$ the corresponding moduli spaces,
$\MM rs := \tM rs/S^1$ and $\Cof rs := \tCof rs/S^1$.
Let us prove that
\begin{equation}
\label{treti_den_v_Bonnu}
H^{\bullet}(\UPQ rs,\pa_1) \cong H^{-\bullet}(\Cof rs).
\end{equation}

To this end, consider the complement $\NN rs : = \MM rs \setminus \Cof rs$.
It is clear that $\MM rs$ is a compact orientable $(r+s -1)$-dimensional
manifold and $\NN rs$ its closed subspace. By definition, $\NN rs$
consists of equivalence classes of configurations such that two (or
more) points of the same type coincide. It therefore
looks locally as an intersection of hyperplanes, thus it is a strong
deformation retract of some open neighborhood $\OO \supset \NN
rs$. Denote finally $\KK := \MM rs \setminus \OO$.
By~\cite[Proposition~3.46]{hatcher:AT},
\[
H^{-\bullet}(\KK) \cong H_{r+s-1 +\bullet}(\MM rs,\OO),
\]
where, since $\NN rs$ is a deformation retract of $\OO$,
$H_\bullet(\MM rs,\OO) \cong H_\bullet(\MM rs,\NN rs)$ and, similarly,
$H^\bullet(\KK)
\cong H^\bullet(\Cof rs)$. We see that
\begin{equation}
\label{co_doma}
H^{-\bullet}(\Cof rs) \cong H_{r+s-1 +\bullet}(\MM rs,\NN rs).
\end{equation}

To describe the right hand side of~(\ref{co_doma}), notice that
$\MM rs$ has a cell structure, with
codimension $d$ cells corresponding to types of configurations with exactly
$d$ collisions. The closed subspace $\NN rs \subset \MM rs$
is a cell subcomplex and
codimension $d$ cells
of the relative cellular chain complex $C(\MM rs,\NN rs)$
correspond to types of configurations with exactly $d$ left-right
collisions.
Obviously, these types are parametrized by
graphs~(\ref{www}) with $d$
triple points. One easily sees
that $(C(\MM rs,\NN rs),\partial)$ is isomorphic
to $(\UPQ rs,\partial_1)$ with the opposite grading shifted by
$r+s-1$, giving rise to the isomorphism
\[
H_{r+s-1+\bullet}(\MM rs,\NN rs) \cong H_{r+s-1+\bullet}(C(\MM rs,\NN
rs),\partial) \cong H^\bullet (\UPQ rs,\partial_1)
\]
which, combined with~(\ref{co_doma}), gives~(\ref{treti_den_v_Bonnu}).

On the other hand, the homotopy type of $\tCof rs$ is easy to
describe: there exists an equivariant deformation retraction that
distributes left points evenly around the cycle so that two adjacent
points are precisely $2\pi/r$ apart, leaving the point
labeled by $1$ unchanged, and similarly distributes right
points leaving the one labeled $r+1$ fixed.
The configurations obtained in this way are
parametrized by the position of points labeled $1$ and $p+1$, plus
the cyclic orders of left and right points, that is
\[
\tCof rs \sim
(C_r \times C_s) \backslash(\Sigma_r \times \Sigma_s)\times
S^1_r \times S^1_s,
\]
where $S_i^1 := S^1$ if $i \geq 1$ while $S_0^1 := \mbox {\rm the
point}$. By the definition of $\Cof rs$,
\[
\Cof rs \sim (C_r \times C_s) \backslash(\Sigma_r \times
\Sigma_s) \times \Srs,
\]
where $\Srs := S^1$ if $r,s \geq 1$ and
$\Srs := \mbox {the
point}$ if $(r,s) \in \{(1,0),(0,1)\}$.

Equation~(\ref{treti_den_v_Bonnu}) implies
\[
H^\bullet(\UPQ rs,\pa_1) \cong H^{-\bullet}(\Cof rs)
\cong
 \bfk[\Sigma_r \times \Sigma_s]^{C_r \times C_s}
\ot  H^{-\bullet}(\Srs),
\]
which, combined with~(\ref{prvni_den_v_Bonnu}), leads to
\Beqrn
H^\bullet(\fA^\circlearrowright_1,\pa_1)(0,n) &\cong & \bigoplus_{r +s = n}
\bfk[\Sigma_n]^{C_r \times C_s} \otimes H^{-\bullet}(\Srs)\\
&=&
\left\{\Ba{ll}
\Ass^\circlearrowright(0,1) & \mbox{for}\ n=1 \\
\Ass^\circlearrowright(0,n) \
\oplus \
{\bigoplus_{p=1}^{n-1}\desusp k[\Sigma_n]^{C_p\times C_{n-p}}}
& \mbox{for}\ n\geq 2.
\Ea
\right.
.
\Eeqrn
proving the claim.
\end{proof}

\begin{proposition}
There is an isomorphism of graded wheeled properads,
$$
H^\bullet\left((\Ass_\infty)^\circlearrowright,\p\right)= {\mathsf S}^\circlearrowright,
$$
where ${\mathsf S}^\circlearrowright$ is the wheeled completion of the
quotient properad,
$$
{\mathsf S}:= \freePROP(\hat{A})/(R),
$$
with the
$\Sigma$-bimodule
$\hat{A}=\{\hat{A}(m,n)\}_{m,n\geq 0}$ given by
$$
\hat{A}(m,n):= \left\{ \Ba{ll} \id_1\ot k[\bS_2]=\mbox{span}\ \left(
\ \
\begin{xy}
<0mm,0mm>*{\bullet},
<0mm,5mm>*{}**@{-},
<-3mm,-5mm>*{}**@{-},
<3mm,-5mm>*{}**@{-},
<-4mm,-7mm>*{_{\sigma(1)}},
<4mm,-7mm>*{_{\sigma(2)}},
\end{xy}
\
\right)_{\sigma\in \Sigma_2} & \mbox{for}\ m=1,n=2 \vspace{3mm}\\
\desusp k[\bS_2]=\mbox{span}\ \left(\ \ \
\begin{xy}
<0mm,0mm>*{\bullet},
<-3mm,-5mm>*{}**@{-},
<3mm,-5mm>*{}**@{-},
<-4mm,-7mm>*{_{\sigma(1)}},
<4mm,-7mm>*{_{\sigma(2)}},
\end{xy} \ \ \
\right)_{\sigma\in \Sigma_2} & \mbox{for}\ m=0,n=2 \vspace{3mm}\\
0 & \mbox{otherwise},
\Ea
\right.
$$
and relations given by (\ref{assrel}) and,

$$
\begin{xy}
<0mm,0mm>*{\bullet},
<-3mm,-5mm>*{}**@{-},
<3mm,-5mm>*{}**@{-},
<-3mm,-5mm>*{\bullet};
<-6mm,-10mm>*{}**@{-},
<0mm,-10mm>*{}**@{-},
<-6mm,-12mm>*{_{\sigma(1)}},
<0mm,-12mm>*{_{\sigma(2)}},
<3mm,-7mm>*{_{\sigma(3)}},
\end{xy}
-
\begin{xy}
<0mm,0mm>*{\bullet},
<-3mm,-5mm>*{}**@{-},
<3mm,-5mm>*{}**@{-},
<-3mm,-5mm>*{\bullet};
<-6mm,-10mm>*{}**@{-},
<0mm,-10mm>*{}**@{-},
<-6mm,-12mm>*{_{\sigma(2)}},
<0mm,-12mm>*{_{\sigma(1)}},
<3mm,-7mm>*{_{\sigma(3)}},
\end{xy}
=0,  \ \ \ \ \
\begin{xy}
<0mm,0mm>*{\bullet},
<-3mm,-5mm>*{}**@{-},
<3mm,-5mm>*{}**@{-},
<3mm,-5mm>*{\bullet};
<6mm,-10mm>*{}**@{-},
<0mm,-10mm>*{}**@{-},
<6mm,-12mm>*{_{\sigma(3)}},
<0mm,-12mm>*{_{\sigma(2)}},
<-3mm,-7mm>*{_{\sigma(1)}},
\end{xy}
-
\begin{xy}
<0mm,0mm>*{\bullet},
<-3mm,-5mm>*{}**@{-},
<3mm,-5mm>*{}**@{-},
<3mm,-5mm>*{\bullet};
<6mm,-10mm>*{}**@{-},
<0mm,-10mm>*{}**@{-},
<6mm,-12mm>*{_{\sigma(2)}},
<0mm,-12mm>*{_{\sigma(2)}},
<-3mm,-7mm>*{_{\sigma(1)}},
\end{xy}
=0,  \ \ \ \ \forall \sigma\in \bS_3.
$$
\end{proposition}

\begin{proof}
It was shown in the proof of Theorem~6.1.1 that every generator
$e\in k[\bS_n]^{C_p\times C_{n-p}}\subset
H^{-1}((\Ass_\infty)^\circlearrowright)$ can be canonically
identified with the $C_p\times C_{n-p}$-orbit,
$$
e=\bigoplus_{i=0}^{p-1}\zeta^i\bigoplus_{j=1}^{n-p-1}\xi^j
\ \ \
\begin{xy}
 <0mm,20mm>*{\bullet};
 <0mm,20mm>*{};<0mm,22mm>*{}**@{.},
 <0mm,20mm>*{};<0mm,-20mm>*{}**@{.},
 <0mm,20mm>*{};<3mm,16mm>*{}**@{-},
<-7mm,16mm>*{_{\zeta(1)}},
<11mm,16mm>*{_{\xi(p+1)}},
 <0mm,20mm>*{};<-3mm,16mm>*{}**@{-},
<0mm,15mm>*{\bullet};
<0mm,15mm>*{};<-3mm,11mm>*{}**@{-},
<-7mm,11mm>*{_{\zeta(2)}},
<-3mm,9mm>*{.};
<-3mm,8mm>*{.};
<-3mm,7mm>*{.};
<0mm,6mm>*{\bullet};
<0mm,6mm>*{};<-3mm,2mm>*{}**@{-},
<-7mm,1mm>*{_{\zeta(p)}},
<0mm,1mm>*{\bullet};
<0mm,1mm>*{};<3mm,-3mm>*{}**@{-},
<11mm,-4mm>*{_{\zeta(p+2)}},
<0mm,-4mm>*{\bullet};
<0mm,-4mm>*{};<3mm,-8mm>*{}**@{-},
<11mm,-9mm>*{_{\zeta(p+3)}},
<3mm,-9mm>*{.};
<3mm,-10mm>*{.};
<3mm,-11mm>*{.};
<0mm,-12mm>*{\bullet};
<0mm,-12mm>*{};<3mm,-16mm>*{}**@{-},
<10mm,-16mm>*{_{\zeta(n)}},
(0,22)*{}
   \ar@{.>}@(ur,dr) (0,-20)*{}
 \end{xy}
$$
of a planar wheeled graph  which has precisely one ternary
vertex and
$n-2$ binary vertices of which $p-1$ vertices have the non-cyclic input leg
pointing ``outside''
the planar wheel. To prove the Proposition it is enough to show
that every such a
linear combination of graphs
is homologically equivalent in the complex
$((\Ass_\infty)^\circlearrowright, \p)$ to a uniquely defined element in
${\mathsf S}^\circlearrowright$.
The latter can be easily established by
induction
on the number of vertices starting with the following initial step,
$$
\begin{xy}
 <0mm,0mm>*{\bullet};
 <0mm,2mm>*{};<0mm,-9mm>*{}**@{.},
 <0mm,0mm>*{};<3mm,-4mm>*{}**@{-},
<-5mm,-4mm>*{_{1}},
<5mm,-4mm>*{_{3}},
 <0mm,0mm>*{};<-3mm,-4mm>*{}**@{-},
<0mm,-5mm>*{\bullet};
<0mm,-5mm>*{};<-3mm,-9mm>*{}**@{-},
<-5mm,-9mm>*{_{2}},
(0,2)*{}
   \ar@{.>}@(ur,dr) (0,-9)*{}
\end{xy}
\ +
\
\begin{xy}
 <0mm,0mm>*{\bullet};
 <0mm,2mm>*{};<0mm,-9mm>*{}**@{.},
 <0mm,0mm>*{};<3mm,-4mm>*{}**@{-},
<-5mm,-4mm>*{_{2}},
<5mm,-4mm>*{_{3}},
 <0mm,0mm>*{};<-3mm,-4mm>*{}**@{-},
<0mm,-5mm>*{\bullet};
<0mm,-5mm>*{};<-3mm,-9mm>*{}**@{-},
<-5mm,-9mm>*{_{1}},
(0,2)*{}
   \ar@{.>}@(ur,dr) (0,-9)*{}
\end{xy}
\ - \
\begin{xy}
 <0mm,0mm>*{\bullet};
 <0mm,3mm>*{};<0mm,-4mm>*{}**@{.},
 <0mm,0mm>*{};<3mm,-4mm>*{}**@{-},
<0.6mm,-10mm>*{_{2}},
<5mm,-4mm>*{_{3}},
 <0mm,0mm>*{};<-3mm,-4mm>*{}**@{-},
<-3mm,-4mm>*{\bullet};
<-3mm,-4mm>*{};<-6mm,-8mm>*{}**@{-},
<-3mm,-4mm>*{};<0mm,-8mm>*{}**@{-},
<-6.6mm,-10mm>*{_{1}},
(0,3)*{}
   \ar@{.>}@(ur,dr) (0,-4)*{}
\end{xy}
\ =\
\p
\begin{xy}
 <0mm,0mm>*{\bullet};
 <0mm,3mm>*{};<0mm,-4mm>*{}**@{.},
 <0mm,0mm>*{};<3mm,-4mm>*{}**@{-},
<-3mm,-6mm>*{_{2}},
<5mm,-4mm>*{_{3}},
 <0mm,0mm>*{};<-3mm,-4mm>*{}**@{-},
<0mm,0mm>*{};<-6mm,-4mm>*{}**@{-},
<-6.7mm,-6mm>*{_{1}},
(0,3)*{}
   \ar@{.>}@(ur,dr) (0,-4)*{}
\end{xy}
$$
\vskip -2em
\end{proof}

\begin{remark}
It is worth emphasizing  a correspondence,
$$
\begin{xy}
<0mm,0mm>*{\bullet},
<-3mm,-5mm>*{}**@{-},
<3mm,-5mm>*{}**@{-},
<-4mm,-7mm>*{_{\sigma(1)}},
<4mm,-7mm>*{_{\sigma(2)}},
\end{xy}
\leftrightsquigarrow
\begin{xy}
 <0mm,-0mm>*{};<3mm,4mm>*{}**@{-},
 <0mm,0.5mm>*{};<4.2mm,-4mm>*{}**@{-},
 <0mm,0.6mm>*{};<0mm,-8mm>*{}**@{-},
 <-0mm,0.48mm>*{};<-4.2mm,-4mm>*{}**@{-},
 <0mm,0mm>*{\bullet};<0mm,0mm>*{}**@{},
   <0.5mm,0.5mm>*{};<4.7mm,-6.5mm>*{^{\sigma(2)}}**@{},
   <-0.48mm,0.48mm>*{};<-4.7mm,-6.5mm>*{^{\sigma(1)}}**@{},
(3,4)*{}
   \ar@{->}@(ur,dr) (0,-8)*{}
 \end{xy}
$$
between the $(0,2)$ generators of the operad $\sf S$ and the wheeled elements
in $(\Ass_\infty)^\circlearrowright$.
\end{remark}

\subsection{Directed oriented ribbon graphs}
The  authors of \cite{MV} were able to compute cohomology of the
directed, $(\fG^\uparrow_g(m,n),\delta)$, version (without wheels
though) of Kontsevich's ribbon graph complex \cite{konts} consisting
of oriented ribbon graphs of genus $g$
 with $n$ input and $m$ output legs,  and show that it is acyclic almost
 everywhere. {\em Wheeled}\, directed ribbon graphs,
$(\fG_g^\circlearrowright(m,n), \delta)$, provide us with a finer
approximation to the original Kontsevich's complex than
$(\fG_g^\uparrow(m,n),\delta)$, and in this case the cohomology
groups $H^\bullet(\fG_g^\circlearrowright(m,n),\delta)$ turn out to
be non-trivial in many degrees.

Let us now give a precise description of all key actors. Consider
first the vector space $\widehat{\fG}^\circlearrowright_g(m,n)$
spanned  by directed {\em oriented ribbon}\,
$(m,n)$-graphs\footnote{The set, $\sG^\circlearrowright(m,n)$, of
directed $(m,n)$-graphs are defined in \S 2. } which are triples,
($G, or(G), \mbox{ribbon structure})$, where

$\bullet$ $G$ is a graph from $\sG^\circlearrowright(m,n)$ satisfying
the conditions: (i) every vertex of $G$ has valence at least $3$
(with at least one ingoing and at least one outgoing edges),
 (ii) vertices of any particular closed path in $G$  are purely
``operadic", i.e.\ they all have either precisely one incoming edge
or they all have precisely one outgoing edge, and (iii) the
associated geometric realization $|G|$ has genus $g$;

$\bullet$ an {\em orientation}\, on a directed $(m,n)$-graph $G$ is,
by definition, an orientation on the vector space,
$\R^{vert(G)}\oplus \R^m\oplus \R^n$ which is in fact the same as an
orientation on $\R^{e(G)}\oplus H_\bullet(|G|,\R)$, where $e(G)$ is
the cardinality of the set, $\edge(G)$, of internal edges of $G$,
and $H_\bullet(|G|,\R)$ is the homology of $|G|$ (see, e.g.
\cite{CV}). Thus we can understand  an orientation of $G$ as an
element
$$
or(G):=or_1(G)\ot or_2(G)\in \det \R^{e(G)}\ot
\det H_\bullet(|G|,\R).
$$

$\bullet$ a {\em ribbon structure}\,  on a directed $(m,n)$-graph
$G$ is,
by definition, an ordering of the set, $\In(v)$, of incoming edges
and the set, $\Out(v)$,
of outgoing edges for each vertex $v\in \vert(G)$;
it can be equivalently understood  as a cyclic ordering of the set
$\In(v)\cup \Out(v)$.

We often abbreviate a triple from $\hat{\fG}^\circlearrowright_g$ by
$(G, or(G))$ or even by $G$.

\begin{definition}
The quotient of the vector space
$\hat{\fG}^\circlearrowright(m,n)$ by
the equivalence relation
$$
(G, -or(G))= -(G, or(G))
$$
is denoted by $\fG_g(m,n)$. It is naturally an
 $\mathbb N$-graded vector space,
$$
{\fG}^\circlearrowright_g(m,n)= \bigoplus_{n\geq 1}{\fG}_g^n(m,n),
$$
with respect to the number of vertices of its elements,
i.e.\ ${\fG}_g^n(m,n):=\{G\in{\fG}^\circlearrowright_g(m,n): \, |vert(G)|=n\}$.
\end{definition}

In fact, it is only the part $\fG_g:= \fG_g^\circlearrowright(0,0)$
which can be regarded as a directed version of Kontsevich's ribbon
graph complex; by contrast to $\fG^\uparrow_g(0,0)=0$, the part
$\fG_g^\circlearrowright(0,0)$ is highly non-trivial.

\begin{fact}[\cite{konts, MV}]
The graded vector space $\fG_g^\circlearrowright(m,n)$ can be made into a cochain complex
by setting
$$
\p(G, or(G)):=\sum_{G'| G'/e=G}(G', or(G'))
$$
where the summation goes over all connected directed oriented ribbon
$(m,n)$-graphs $G'$ such that

$\bullet$
$G$ can be obtained from $G'$ by contracting an internal edge
  $e\in e(G')$ which is not a loop;

$\bullet$
the cyclic ordering of edges at the vertex $v\in vert(G)$ into
which $e$ contracts agrees in the obvious sense with the one induced from the
contraction (see, e.g., \S 2.2.2 in \cite{CV} for full details and pictures);

$\bullet$
$or(G')=(or_1(G)\wedge e)\ot or_2(G)$.

Set $\fG_g^\circlearrowright:=\{ \fG_g^\circlearrowright(m,n)\}_{m,n\geq 0}$.
It is naturally an $\bS$-bimodule.
\end{fact}

\begin{theorem}\label{hg}
$H^\bullet(\fG_g^\circlearrowright)$ is isomorphic as an
$\bS$-bimodule to the ordinary PROP,\footnote{For a PROP
$P=\{P(m,n)\}$ we denote by $P^\dag=\{P^\dag(m,n)\}$ the associated
PROP with ``reversed flow", i.e. $P^\dag(m,n):= P(n,m)$.}
$$
\frac{{\mathsf S}^\circlearrowright *
({\mathsf S}^\circlearrowright)^\dag}{I_0}
$$
where ${\mathsf S}^\circlearrowright$ is defined in \S 6.1.2,
$*$ stands
for the free product of ordinary PROPs, and the ideal $I_0$ is generated by
the elements
$$
\begin{xy}
<0mm,0mm>*{\bullet},
<0mm,4mm>*{}**@{-},
<-3mm,-4mm>*{}**@{-},
<3mm,-4mm>*{}**@{-},
<-4mm,-7mm>*{_{\sigma(1)}},
<4mm,-7mm>*{_{\sigma(2)}},
<0mm,4mm>*{\bullet},
<0mm,4mm>*{};<-3mm,8mm>*{}**@{-},
<-3mm,8mm>*{}**@{-},
<3mm,8mm>*{}**@{-},
<-4mm,11mm>*{_{\tau(1)}},
<4mm,11mm>*{_{\tau(2)}},
\end{xy}\ \ \ , \ \ \ \ \forall \sigma, \tau\in \bS_2.
$$
\end{theorem}

\begin{proof}
Exactly the same argument as in the final section in \cite{MV}
establishes an isomorphism of the complex $(\fG_g^\circlearrowright,
\p)$ with the operadic wheeled completion\footnote{see \S 3.11 in
\cite{me3} for the definition of the operadic wheelification
functor.}, $({ \sf IB}_\infty^\looparrowright, \delta)$, of the
minimal dg PROP resolution, $({\sf IB}_\infty, \delta)$, of the
PROP, $\sf IB$, of infinitesimal associative bialgebras, which is,
by definition, the quotient,
$$
{\sf IB}:=
\frac{{\Ass} *
{\Ass}^\dag}{I},
$$
by the ideal $I$ generated by
the elements
$$
\begin{xy}
<0mm,0mm>*{\bullet},
<0mm,4mm>*{}**@{-},
<-3mm,-4mm>*{}**@{-},
<3mm,-4mm>*{}**@{-},
<-4mm,-7mm>*{_{\sigma(1)}},
<4mm,-7mm>*{_{\sigma(2)}},
<0mm,4mm>*{\bullet},
<0mm,4mm>*{};<-3mm,8mm>*{}**@{-},
<-3mm,8mm>*{}**@{-},
<3mm,8mm>*{}**@{-},
<-4mm,11mm>*{_{\tau(1)}},
<4mm,11mm>*{_{\tau(2)}},
\end{xy}
\  - \
\begin{xy}
<0mm,0mm>*{\bullet};<0mm,0.8mm>*{}**@{},
 <3mm,4mm>*{\bullet};<2.4mm,2.4mm>*{}**@{},
 <0mm,0mm>*{};<3mm,4mm>*{}**@{-},
 <0mm,0mm>*{};<-3mm,4mm>*{}**@{-},
 <0mm,0mm>*{};<0mm,-4mm>*{}**@{-},
 <3mm,4mm>*{};<6mm,0mm>*{}**@{-},
 <3mm,4mm>*{};<3mm,8mm>*{}**@{-},
 <-5mm,5.6mm>*{^{\tau(1)}}**@{},
<3mm,10mm>*{_{\tau(2)}}**@{},
<0mm,-7mm>*{_{\sigma(1)}}**@{},
<7.8mm,-2mm>*{_{\sigma(2)}}**@{},
    \end{xy}
   \ \ \  - \ \
\begin{xy}
<0mm,0mm>*{\bullet};<0mm,0.8mm>*{}**@{},
 <-3mm,4mm>*{\bullet};<2.4mm,2.4mm>*{}**@{},
 <0mm,0mm>*{};<3mm,4mm>*{}**@{-},
 <0mm,0mm>*{};<-3mm,4mm>*{}**@{-},
 <0mm,0mm>*{};<0mm,-4mm>*{}**@{-},
 <-3mm,4mm>*{};<-6mm,0mm>*{}**@{-},
 <-3mm,4mm>*{};<-3mm,8mm>*{}**@{-},
 <-7mm,-2.5mm>*{^{\sigma(1)}}**@{},
<-3mm,10mm>*{_{\tau(1)}}**@{},
<0mm,-7mm>*{_{\sigma(2)}}**@{},
<7.8mm,5mm>*{_{\tau(2)}}**@{},
    \end{xy}
\ \ \ , \ \ \ \ \forall \sigma, \tau\in \bS_2.
$$
As the ideal $I$ is related to the ideal $I_0$ via distributive law
(see \cite{gan,MV}), one can apply Theorem~3.11.1(ii) from
\cite{me3} to compute the cohomology of the complex $({\sf
IB}_\infty^\looparrowright, \delta)$ as an $\bS$-bimodule,
$$
H^\bullet({\sf IB}_\infty^\looparrowright, \delta)\simeq
\frac{{H^\bullet((\Ass_\infty)^\circlearrowright}) *
H^\bullet((\Ass_\infty)^\circlearrowright)^\dag}{I_0},
$$
Finally,
Proposition~6.1.2 completes the proof.
\end{proof}

\subsection{Proof of Theorem C (see \S 1)}
 Since $H^\bullet(\fG_g)$ is,
by definition of $\fG_g$, equal to the $H^\bullet(\fG_g^\circlearrowright)(0,0)$ part
of the $\bS$-bimodule  $H^\bullet(\fG_g^\circlearrowright)=\{
 H^\bullet(\fG_g^\circlearrowright)(m,n)\}_{m,n\geq 0}$, we have an
 isomorphism
$$
H^\bullet(\fG_g)\simeq \frac{H^\bullet((\Ass_\infty)^\circlearrowright) *
H^\bullet((\Ass_\infty)^\circlearrowright)^\dag}{I_0}(0,0).
$$
As to the r.h.s. contribute only elements of type $(0,m)$ from
$H^\bullet((\Ass_\infty)^\circlearrowright)$ and elements of type
$(n,0)$ from $H^\bullet((\Ass_\infty)^\circlearrowright)^\dag$, the required result
follows from Theorem~6.1.1.
\hfill $\Box$

\subsection{Proof of Corollary D (see \S 1)}
By Theorems C and 6.1.1, $H^\bullet(\fG_g)$, is isomorphic to the vector space
spanned by all possible directed $(0,0)$-graphs which can be  obtained by
gluing input legs of a disjoint union of $(0,n)$-graphs of the types,
$$
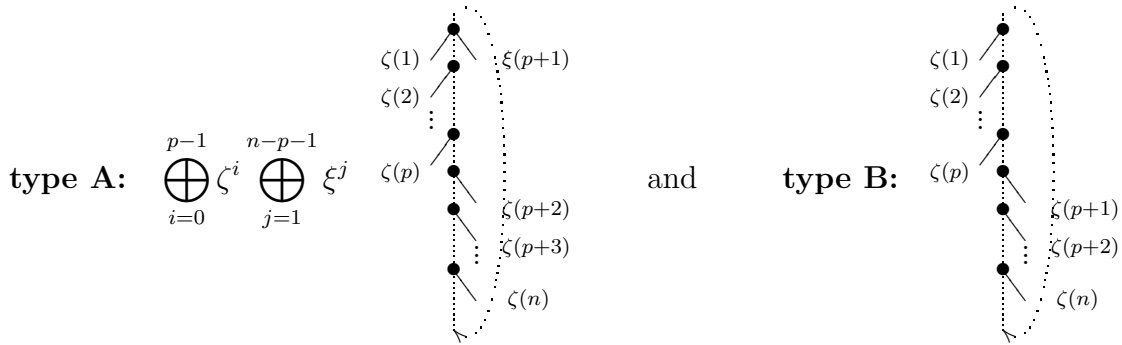
\mbox{\bf type A:}\ \ \ \bigoplus_{i=0}^{p-1}\zeta^i\bigoplus_{j=1}^{n-p-1}\xi^j
\ \ \
\begin{xy}
 <0mm,20mm>*{\bullet};
 <0mm,20mm>*{};<0mm,22mm>*{}**@{.},
 <0mm,20mm>*{};<0mm,-20mm>*{}**@{.},
 <0mm,20mm>*{};<3mm,16mm>*{}**@{-},
<-7mm,16mm>*{_{\zeta(1)}},
<11mm,16mm>*{_{\xi(p+1)}},
 <0mm,20mm>*{};<-3mm,16mm>*{}**@{-},
<0mm,15mm>*{\bullet};
<0mm,15mm>*{};<-3mm,11mm>*{}**@{-},
<-7mm,11mm>*{_{\zeta(2)}},
<-3mm,9mm>*{.};
<-3mm,8mm>*{.};
<-3mm,7mm>*{.};
<0mm,6mm>*{\bullet};
<0mm,6mm>*{};<-3mm,2mm>*{}**@{-},
<-7mm,1mm>*{_{\zeta(p)}},
<0mm,1mm>*{\bullet};
<0mm,1mm>*{};<3mm,-3mm>*{}**@{-},
<11mm,-4mm>*{_{\zeta(p+2)}},
<0mm,-4mm>*{\bullet};
<0mm,-4mm>*{};<3mm,-8mm>*{}**@{-},
<11mm,-9mm>*{_{\zeta(p+3)}},
<3mm,-9mm>*{.};
<3mm,-10mm>*{.};
<3mm,-11mm>*{.};
<0mm,-12mm>*{\bullet};
<0mm,-12mm>*{};<3mm,-16mm>*{}**@{-},
<10mm,-16mm>*{_{\zeta(n)}},
(0,22)*{}
   \ar@{.>}@(ur,dr) (0,-20)*{}
 \end{xy} \hspace{10mm} \mbox{and}\ \hspace{10mm} \mbox{\bf type B:}\ \ \
 \begin{xy}
 <0mm,20mm>*{\bullet};
 <0mm,20mm>*{};<0mm,22mm>*{}**@{.},
 <0mm,20mm>*{};<0mm,-20mm>*{}**@{.},
<-7mm,16mm>*{_{\zeta(1)}},
 <0mm,20mm>*{};<-3mm,16mm>*{}**@{-},
<0mm,15mm>*{\bullet};
<0mm,15mm>*{};<-3mm,11mm>*{}**@{-},
<-7mm,11mm>*{_{\zeta(2)}},
<-3mm,9mm>*{.};
<-3mm,8mm>*{.};
<-3mm,7mm>*{.};
<0mm,6mm>*{\bullet};
<0mm,6mm>*{};<-3mm,2mm>*{}**@{-},
<-7mm,1mm>*{_{\zeta(p)}},
<0mm,1mm>*{\bullet};
<0mm,1mm>*{};<3mm,-3mm>*{}**@{-},
<11mm,-4mm>*{_{\zeta(p+1)}},
<0mm,-4mm>*{\bullet};
<0mm,-4mm>*{};<3mm,-8mm>*{}**@{-},
<11mm,-9mm>*{_{\zeta(p+2)}},
<3mm,-9mm>*{.};
<3mm,-10mm>*{.};
<3mm,-11mm>*{.};
<0mm,-12mm>*{\bullet};
<0mm,-12mm>*{};<3mm,-16mm>*{}**@{-},
<10mm,-16mm>*{_{\zeta(n)}},
(0,22)*{}
   \ar@{.>}@(ur,dr) (0,-20)*{}
 \end{xy}
$$
with output legs of a disjoint union of
 similar $(m,0)$-graphs of the corresponding types $\bar{A}$ and $\bar{B}$.

Thus with every nonzero element $e$ of $H^\bullet(\fG_g)$ we can associate
\Bi
\item a collection of $p$ graphs of type $A$ each having $a_i$ vertices
(and hence $a_i+1$ input legs), $1\leq i \leq p$;
\item a collection of $q$ graphs of type $B$ each having $b_j$ vertices
(and hence $b_j$ input legs), $1\leq j \leq q$;
\item a collection of $\bar{p}$ graphs of type $\bar{A}$ each having $\bar{a}_i$ vertices
(and hence $\bar{a}_{\bar{i}}+1$ output legs), $1\leq \bar{i} \leq \bar{p}$;
\item a collection of $\bar{q}$ graphs of type $\bar{B}$ each having
$\bar{b}_{\bar{j}}$ vertices
(and hence $\bar{b}_{\bar{i}}$ output legs), $1\leq \bar{j} \leq \bar{q}$;
\Ei
such that total number, say $N$, of input legs of graphs of types $A$ and $B$
is equal to the total number of output legs of graphs of types $\bar{A}$ and
$\bar{B}$, i.e.
\begin{equation}\label{degree}
N= \sum_{i=1}^p a_i + \sum_{j=1}^q b_j +p =  \sum_{\bar{i}=1}^{\bar{p}}
\bar{a}_{\bar{i}} +
\sum_{\bar{j}=1}^{\bar{q}} \bar{b}_{\bar{j}} +\bar{p}.
\end{equation}
The total number of internal edges of a graph $G$ obtained from the above data
by gluing $N$ output legs with $N$ input legs is equal to
$3N  -p - \bar{p}$ while the total number of vertices is $2N-p-\bar{p}$.
By the Euler formula, the genus of $G$ is
$$
g= 1+ (3N-p-\bar{p}) - (2N - p - \bar{p})=N+1.
$$
Thus our non-zero element $e$ belongs to $H^{2g-2-p-\bar{p}}(\fG_g)$ and we immediately
conclude that $H^n(\fG_g)=0$ for all $n> 2g-2$.

\sip

It is clear from (\ref{degree}) that, for fixed $N=g-1$, the maximal
possible values of parameters $p$ and $\bar{p}$ are equal to $[N/2]$.
Hence $H^n(\fG_g)$ is indeed non-zero  only for $n$ in the range
 $g-\frac{1}{2}(1-(-1)^g)\leq n \leq 2(g-1)$.
$\Box$

\begin{example}
$H^2(\fG_3)$ is two dimensional and is spanned by the following graphs,
$$
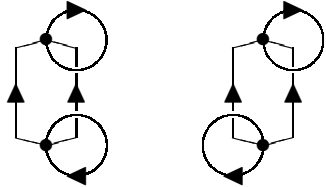
\xy
*\xycircle(4,4){-};
<0mm,4mm>*{\blacktriangleright};
<-4mm,0mm>*{\bullet};
<-4mm,0mm>*{};<0mm,-1mm>*{}**@{-},
<0mm,-1mm>*{};<0mm,-3.5mm>*{}**@{-},
<0mm,-4.5mm>*{};<0mm,-9.4mm>*{}**@{-},
<0mm,-6.5mm>*{};<0mm,-9.4mm>*{}**@{-},
<0mm,-7mm>*{\blacktriangle};
<-4mm,0mm>*{};<-8mm,-1mm>*{}**@{-},
<-8mm,-1mm>*{};<-8mm,-13mm>*{}**@{-},
<-8mm,-13mm>*{};<-4mm,-14mm>*{}**@{-},
<-8mm,-7mm>*{\blacktriangle};
<0mm,-14mm>*{}*\xycircle<11.5pt>{-};
<0mm,-18mm>*{\blacktriangleleft};
<-4mm,-14mm>*{\bullet};
<-4mm,-14mm>*{};<0mm,-13mm>*{}**@{-},
<0mm,-13mm>*{};<0mm,-10.5mm>*{}**@{-},
\endxy
\ \hspace{10mm} \
\xy
*\xycircle(4,4){-};
<0mm,4mm>*{\blacktriangleright};
<-4mm,0mm>*{\bullet};
<-4mm,0mm>*{};<0mm,-1mm>*{}**@{-},
<0mm,-1mm>*{};<0mm,-3.5mm>*{}**@{-},
<0mm,-4.5mm>*{};<0mm,-9.4mm>*{}**@{-},
<0mm,-6.5mm>*{};<0mm,-13mm>*{}**@{-},
<0mm,-7mm>*{\blacktriangle};
<-4mm,0mm>*{};<-8mm,-1mm>*{}**@{-},
<-8mm,-1mm>*{};<-8mm,-9.4mm>*{}**@{-},
<-8mm,-13mm>*{};<-4mm,-14mm>*{}**@{-},
<-8mm,-7mm>*{\blacktriangle};
<-8mm,-14mm>*{}*\xycircle<11.5pt>{-};
<-8mm,-18mm>*{\blacktriangleleft};
<-4mm,-14mm>*{\bullet};
<-4mm,-14mm>*{};<0mm,-13mm>*{}**@{-},
<-8mm,-13mm>*{};<-8mm,-10.5mm>*{}**@{-},
\endxy
$$
\end{example}

\subsection{Proof of Theorem A (see \S 1)}
We have to show that the natural morphism of dg wheeled
properads,
$$
p: \left({\Ass}_\infty^\circlearrowright, \p\right) \lon
\left(\Ass^\circlearrowright, 0\right),
$$
given on generators by
\Beqrn
p\left(
\begin{xy}
 <0mm,0mm>*{\bullet};<0mm,0mm>*{}**@{},
 <0mm,0mm>*{};<-8mm,-5mm>*{}**@{-},
 <0mm,0mm>*{};<-4.5mm,-5mm>*{}**@{-},
 <0mm,0mm>*{};<0mm,-4mm>*{\ldots}**@{},
 <0mm,0mm>*{};<4.5mm,-5mm>*{}**@{-},
 <0mm,0mm>*{};<8mm,-5mm>*{}**@{-},
   <0mm,0mm>*{};<-9mm,-7.9mm>*{^1}**@{},
   <0mm,0mm>*{};<-4.6mm,-7.9mm>*{^2}**@{},
   <0mm,0mm>*{};<9.4mm,-7.9mm>*{^n}**@{},
 <0mm,0mm>*{};<0mm,5mm>*{}**@{-},
 \end{xy}\ \
\right) &:=&
\left\{
\Ba{cr}
\Id
 & \mbox{for}\ n=2 \vspace{6mm}\\
0 & \mbox{otherwise}
\Ea
\right.
\\
p\left(
\begin{xy}
 <0mm,-0.5mm>*{\blacktriangledown};<0mm,0mm>*{}**@{},
 <0mm,0mm>*{};<-16mm,-5mm>*{}**@{-},
 <0mm,0mm>*{};<-12mm,-5mm>*{}**@{-},
 <0mm,0mm>*{};<-3.5mm,-5mm>*{}**@{-},
 <0mm,0mm>*{};<-6mm,-5mm>*{...}**@{},
   <0mm,0mm>*{};<-19mm,-7.9mm>*{^{1}}**@{},
   <0mm,0mm>*{};<-13mm,-7.9mm>*{^{2}}**@{},
   <0mm,0mm>*{};<-4mm,-7.9mm>*{^{m}}**@{},
 <0mm,0mm>*{};<16mm,-5mm>*{}**@{-},
 <0mm,0mm>*{};<12mm,-5mm>*{}**@{-},
 <0mm,0mm>*{};<3.5mm,-5mm>*{}**@{-},
 <0mm,0mm>*{};<6.6mm,-5mm>*{...}**@{},
   <0mm,0mm>*{};<19mm,-7.9mm>*{^{m\hspace{-0.5mm}+\hspace{-0.5mm}n}}**@{},
   <0mm,0mm>*{};<5mm,-7.9mm>*{^{m\hspace{-0.5mm}+\hspace{-0.5mm}1}}**@{},
 \end{xy}
\right)&:=& 0 \ \ \forall m,n\geq 1,
\Eeqrn
is a quasi-isomorphism.

Consider an increasing filtration, $F_0\subset F_1\subset\ldots \subset F_p\subset\ldots$ of the complex
$(\Ass_\infty^\circlearrowright, \p)$ by the number $p$ of $\bullet$-vertices, and let
$\{E_r,\delta_r\}_{r\geq 0}$ be the associated spectral sequence. The differential $\delta_0$ is non-zero
only on $\blacktriangledown$-vertices, and is given by
\begin{eqnarray*}
\delta_0
\hspace{-2mm}
\begin{xy}
 <0mm,-0.5mm>*{\blacktriangledown};<0mm,0mm>*{}**@{},
 <0mm,0mm>*{};<-16mm,-5mm>*{}**@{-},
 <0mm,0mm>*{};<-12mm,-5mm>*{}**@{-},
 <0mm,0mm>*{};<-3.5mm,-5mm>*{}**@{-},
 <0mm,0mm>*{};<-6mm,-5mm>*{...}**@{},
   <0mm,0mm>*{};<-19mm,-7.9mm>*{^{1}}**@{},
   <0mm,0mm>*{};<-13mm,-7.9mm>*{^{2}}**@{},
   <0mm,0mm>*{};<-4mm,-7.9mm>*{^{m}}**@{},
 <0mm,0mm>*{};<16mm,-5mm>*{}**@{-},
 <0mm,0mm>*{};<12mm,-5mm>*{}**@{-},
 <0mm,0mm>*{};<3.5mm,-5mm>*{}**@{-},
 <0mm,0mm>*{};<6.6mm,-5mm>*{...}**@{},
   <0mm,0mm>*{};<19mm,-7.9mm>*{^{m\hspace{-0.5mm}+\hspace{-0.5mm}n}}**@{},
   <0mm,0mm>*{};<5mm,-7.9mm>*{^{m\hspace{-0.5mm}+\hspace{-0.5mm}1}}**@{},
 \end{xy} \hskip -.5em
& = &
\sum_{i=0}^{m-1}\left((-1)^{m+1}\zeta\right)^i\sum_{j=1}^{n-1} \left(
(-1)^{n+1}\xi\right)^j\left(
  \sum_{k=2}^{m}
(-1)^{k+m}
\hspace{-7mm}
\begin{xy}
 <0mm,-0.5mm>*{\blacktriangledown};<0mm,0mm>*{}**@{},
 <0mm,0mm>*{};<-16mm,-5mm>*{}**@{-},
 <0mm,0mm>*{};<-11mm,-5mm>*{}**@{-},
 <0mm,0mm>*{};<-3.5mm,-5mm>*{}**@{-},
 <0mm,0mm>*{};<-6mm,-5mm>*{...}**@{},
   <0mm,0mm>*{};<-10mm,-7.9mm>*{^{k\hspace{-0.5mm}+\hspace{-0.5mm}1}}**@{},
   <0mm,0mm>*{};<-4mm,-7.9mm>*{^{m}}**@{},
 <0mm,0mm>*{};<16mm,-5mm>*{}**@{-},
 <0mm,0mm>*{};<12mm,-5mm>*{}**@{-},
 <0mm,0mm>*{};<3.5mm,-5mm>*{}**@{-},
 <0mm,0mm>*{};<6.6mm,-5mm>*{...}**@{},
   <0mm,0mm>*{};<17.3mm,-7.9mm>*{^{m\hspace{-0.5mm}+\hspace{-0.5mm}n}}**@{},
   <0mm,0mm>*{};<5mm,-7.9mm>*{^{m\hspace{-0.5mm}+\hspace{-0.5mm}1}}**@{},
 <-16mm,-5.5mm>*{\bullet};<0mm,0mm>*{}**@{},
 <-16mm,-5.5mm>*{};<-20mm,-11mm>*{}**@{-},
 <-16mm,-5.5mm>*{};<-12mm,-11mm>*{}**@{-},
 <-16mm,-5.5mm>*{};<-18mm,-11mm>*{}**@{-},
 <-16mm,-5.5mm>*{};<-14mm,-11mm>*{}**@{-},
 <-16mm,-13mm>*{_{1\ \, \dots\ \ k}},
 \end{xy}\right.\\
\nonumber\\
&&\left. + \sum_{k=2}^{n-2}
(-1)^{k+m+n}
\begin{xy}
<0mm,-0.5mm>*{\blacktriangledown};<0mm,0mm>*{}**@{},
 <0mm,0mm>*{};<-16mm,-5mm>*{}**@{-},
 <0mm,0mm>*{};<-12mm,-5mm>*{}**@{-},
 <0mm,0mm>*{};<-3.5mm,-5mm>*{}**@{-},
 <0mm,0mm>*{};<-6mm,-5mm>*{...}**@{},
   <0mm,0mm>*{};<-18mm,-7.9mm>*{^{1}}**@{},
   <0mm,0mm>*{};<-13mm,-7.9mm>*{^{2}}**@{},
   <0mm,0mm>*{};<-4mm,-7.9mm>*{^{m}}**@{},
 <0mm,0mm>*{};<17mm,-5mm>*{}**@{-},
 <0mm,0mm>*{};<7mm,-5mm>*{}**@{-},
 <0mm,0mm>*{};<3.5mm,-5mm>*{}**@{-},
 <0mm,0mm>*{};<11.6mm,-5mm>*{...}**@{},
   <0mm,0mm>*{};<21mm,-7.9mm>*{^{m\hspace{-0.5mm}+\hspace{-0.5mm}n}}**@{},
   <0mm,0mm>*{};<11mm,-7.9mm>*{^{m\hspace{-0.5mm}+\hspace{-0.5mm}k+\hspace{-0.5mm}1}}**@{},
 <3.5mm,-5.5mm>*{\bullet};<0mm,0mm>*{}**@{},
 <3.5mm,-5.5mm>*{};<-0.5mm,-11mm>*{}**@{-},
 <3.5mm,-5.5mm>*{};<1.5mm,-11mm>*{}**@{-},
 <3.5mm,-5.5mm>*{};<5.5mm,-11mm>*{}**@{-},
 <3.5mm,-5.5mm>*{};<7.5mm,-11mm>*{}**@{-},
 <3.5mm,-13mm>*{_{m+1\, \dots\ m+k}},
 \end{xy}\right)
\nonumber
\end{eqnarray*}

Let $(CC_{-*}(W), \p_W)$ be the (negatively graded)
cellular chain complex of
the cyclohedron introduced and studied in
\cite{cyclo}. A comparison of the differential $\delta_0$ with the
boundary operator $\p_W$ (given explicitly in Proposition 2.14 of
\cite{cyclo}) immediately implies that, modulo action of finite
groups, the complex $(E_0,\delta_0)$ is isomorphic to the tensor
product of a trivial complex (i.e.\ one with vanishing differential)
with the complex $\R[\bS_{*+\bullet +2}]\ot \desusp CC_{-*}(W)\ot
\hbox{$\desusp CC_{-\bullet}(W)$}$.
As cyclohedron consists of contractible polytopes, its
homology is concentrated in degree $0$ and is equal to $\R$.  Thus
$(E_1,\delta_1)$ is the free
wheeled properad generated by
$\bullet$-corollas and extra binary corollas in degree $-2$,
$$
\begin{xy}
 <0mm,-0.5mm>*{\blacktriangledown};<0mm,0mm>*{}**@{},
 <0mm,0mm>*{};<-3.5mm,-5mm>*{}**@{-},
   <0mm,0mm>*{};<-4mm,-7.9mm>*{^{\sigma(1)}}**@{},
 <0mm,0mm>*{};<3.5mm,-5mm>*{}**@{-},
   <0mm,0mm>*{};<5mm,-7.9mm>*{^{\sigma(2)}}**@{},
 \end{xy}\  \ \ \ \sigma\in \bS_2.
$$
standing for the natural basis of
$\R[\bS_2]\ot H^0(\desusp CC_{-*}(W)\ot \desusp CC_{-*}(W))=
\downarrow^2 \R[\bS_2]$.

The induced differential $\delta_1$ is given on $\bullet$-corollas by the usual $A_\infty$ formula,
while on the remaining binary $\blacktriangledown$-corollas as

$$
\delta_1
\begin{xy}
 <0mm,-0.5mm>*{\blacktriangledown};<0mm,0mm>*{}**@{},
 <0mm,0mm>*{};<-3.5mm,-5mm>*{}**@{-},
   <0mm,0mm>*{};<-4mm,-7.9mm>*{^{\sigma(1)}}**@{},
 <0mm,0mm>*{};<3.5mm,-5mm>*{}**@{-},
   <0mm,0mm>*{};<5mm,-7.9mm>*{^{\sigma(2)}}**@{},
 \end{xy} =
\begin{xy}
 <0mm,-0mm>*{};<3mm,4mm>*{}**@{-},
 <0mm,0.5mm>*{};<4.2mm,-4mm>*{}**@{-},
 <0mm,0.6mm>*{};<0mm,-8mm>*{}**@{-},
 <-0mm,0.48mm>*{};<-4.2mm,-4mm>*{}**@{-},
 <0mm,0mm>*{\bullet};<0mm,0mm>*{}**@{},
   <0.5mm,0.5mm>*{};<4.7mm,-6.5mm>*{^{\sigma(2)}}**@{},
   <-0.48mm,0.48mm>*{};<-4.7mm,-6.5mm>*{^{\sigma(1)}}**@{},
(3,4)*{}
   \ar@{->}@(ur,dr) (0,-8)*{}
 \end{xy}.
$$
By Theorem~6.1.1 and Proposition~6.1.2, $E_2$ is precisely $\Ass^\circlearrowright$. Hence the induced differential
$\delta_2$ vanishes so that converging to $H(\Ass_\infty^\circlearrowright,\p)$  spectral sequence
$\{E_r,d_r\}$ degenerates at this second term. This completes the proof.
\hfill $\Box$

\subsection{Wheeled Massey operations}
Let $(V, \cdot, d)$ be a differential graded associative algebra over a field
$\bfk$, and let
$(W, d)$ be a complex together
with morphisms of complexes,
$i: W\rar V$ and $p: V\rar W$ such that
$$
p\circ i= \Id + Qd + dQ
$$
for some degree $-1$ linear operator $Q: V\rar V$. Homotopy theory general nonsense
says that $W$ must have an induced $\Ass_\infty$-structure,
$\{\mu_n: W^{\ot n}\rar\ \susp^{n-2} W\}_{n\geq 1}$,
and the following explicit formulae for this structure have been found in \cite{Me0}:
\Beqrn
\mu_1&=&d, \\
\mu_n &=& p\circ \lambda_n\circ i^{\ot n} \ \ \ \mbox{for}\ n\geq 2,
\Eeqrn
where $\lambda_n: V^{\ot n}\rar\ \susp^{n-2} V$
are defined inductively as follows,
$$
\lambda_n(v_1, \ldots, v_n):= \sum_{k+l=n+1\atop k,l\geq 1}
(-1)^{k+(l-1)(|v_1|+\ldots+|v_n|)}Q\lambda_{k}(v_1,\ldots, v_k)\, \cdot \,
Q\lambda_{n-k}(v_{k+1},\ldots, v_n),
$$
starting with formal equality $\lambda_1:=-Q^{-1}$.

If $V$ is finite dimensional, then the datum $(V, \cdot, d)$ is
canonically an
$\Ass^\circlearrowright$-algebra so that again we may expect that
$W$ (which may not be finite-dimensional!) has a naturally induced $\Ass_\infty^\circlearrowright$-structure,
$$
\left\{\mu_n: W^{\ot n} \rar\ \susp^{n-2}W\right\}_{n\geq 1}, \ \ \ \
\left\{\mu_{p,q}: W^{\ot p+q} \rar\ \susp^{p+q} \bfk \right\}_{p,q\geq 1},
$$
which satisfy quadratic equations mimicking formulae (\ref{ass1}) and
(\ref{ass2}) respectively. Straightforward but very tedious calculations show
that this is indeed so with $\mu_n$ given by the formulae above
and the new set of operations $\mu_{p,q}$ given as follows:
if ${\it Tr}_i(Q\circ\lambda_n): V^{\ot n-1} \rar\ \susp^{n-1} \bfk$
stands for the trace of
the linear map $Q\circ \lambda_n: V^{\ot n}\rar\ \susp^{n-1} V$
with respect to the $i$-th input, $2\leq i\leq n-1$, $n\geq 3$, then one has
$$
\mu_{i-1,n-i}\left(w_1,\ldots, w_{i-1}, w_{i+1}, \ldots w_{n}\right)=
\hspace{80mm}
$$
$$
\hspace{10mm}
=\frac{1}{2}
\sum_{l=0}^{i-2}\sum_{j=1}^{n-i-1}(-1)^{li+1+(n-i+1)j}
Tr_i(Q\circ\lambda_n)\left(i(w_{\zeta^l(1)}), \ldots, i( w_{\zeta^l(i-1)}),
i(w_{\xi^j(i+1)}), \ldots,  i(w_{\xi^j(n)})\right)
$$
where
$$
\zeta:=(12\ldots i-1)\in \bS_{i-1}, \ \ \  \ \xi:=((i+1)(i+2)\ldots n)\in
\bS_{n-i}
$$
and $w_\bullet$ are arbitrary elements of $W$.

In particular, the cohomology of any finite dimensional dg associative algebra
is naturally an $\Ass_\infty^\circlearrowright$-algebra and hence admits,
in general, lots of new Massey type cyclically (skew)symmetric
operations corresponding to compositions $\mu_{p,q}$.

\subsection{Cyclic characteristic class of an $\Ass_\infty$-algebra}
Any finite dimensional representation
of the operad $\Ass$ extends naturally to a representation
of the wheeled properad $\Ass^\circlearrowright$. It is {\em not}\, true, however, that any
 finite dimensional representation of the dg operad $\Ass_\infty$ can be extended to a representation
of the dg wheeled properad $\Ass_\infty^\circlearrowright$ --- there exist a cohomological obstruction,
and the main purpose of this subsection is to describe it, that is,
to associate with an arbitrary  finite dimensional
$\Ass_\infty$-algebra
a cohomology class whose vanishing provides a
necessary and sufficient condition
for  existence of its wheeled extension.

It is well known that with any graded vector space $V$ one can associate a graded Lie algebra
$$
C^\bullet(V,V)=\bigoplus_{n}C^n(V,V), \ \ \ C^n(V,V):=\bigoplus_{j\geq 1} \Hom_{n+1}(\otexp{(\susp V)}j,V)
$$
equipped with the Gerstenhaber Lie bracket, $[\ \hskip .15em  ,\ ]_G$.
Here ${\it Hom}_{i}(V^{\ot j}, V)$ stands for the space
of homogeneous maps $V^{\ot j}\rar V$ of degree $i$.

Maurer-Cartan elements in this Lie algebra, that is, elements $\Gamma\in C^1(V,V)$ such that
$[\Gamma,\Gamma]_G=0$, are in 1-1 correspondence with $\Ass_\infty$-structures in $V$.
Such elements make also $(C^\bullet(V,V), [\ \hskip .15em ,\ ]_G)$
into a {\em differential}\, graded Lie algebra, $(C^\bullet(V,V), [\
  \hskip .15em,\ ]_G, D_\Gamma)$,
with $D_\Gamma:=[\Gamma,\ ]_G$.

Define next a graded vector space,
$$
{\it Cyc}^\bullet(V,\bfk):=\bigoplus_n {\it Cyc}^n(V,\bfk),\ \ \
{\it Cyc}^n(V,\bfk):=
 \bigoplus_{
 p,q\geq 1}{\it Hom}_n\left(((\downarrow V)^{\ot p})_{C_p} \ot
((\downarrow V)^{\ot q})_{C_q}, \bfk\right).
$$
It is naturally a module over the graded Lie algebra $(C^\bullet(V,V),
[\ \hskip .15em ,\ ]_G)$ so that the direct
sum,
$$
\fg:= C^\bullet(V,V)  \bigoplus {\it Cyc}^\bullet(V,\bfk),
$$
has structure of a graded Lie algebra with the Lie bracket defined by
$$
[a,b]:=[a,b]_G,\ \ \ [a,x]:= a\circ x, \ \ \ [x,y]=0
$$
for arbitrary $a,b\in C^\bullet(V,V)$, $x,y\in {\it Cyc}^\bullet(V,\bfk)$. Here $\circ$ stands for the action of
$C^\bullet(V,V)$ on ${\it Cyc}^\bullet(V,\bfk)$.

Any Maurer-Cartan element
$\Gamma\in C^1(V,V)$
defines a Maurer-Cartan element, $\Gamma\oplus 0$, in $\fg$ and
hence makes the latter and, in particular, the subspace ${\it Cyc}^\bullet(V)$  into a complex with the differential
$$
\delta_\Gamma:=[\Gamma\oplus 0, \ \ \ ].
$$

\begin{remark}
{\rm
In fact, general
Maurer-Cartan elements in $\fg$ describe strongly
homotopy $\Ass^\omega$-algebras, where $\Ass^\omega$-structures consist of an
associative multiplication $\mu : V \ot V\to V$ and an operation
$\omega : V \ot V \to \bfk$ such that, besides the associativity of
$\mu$, $\omega(\mu(a,b),c)$ is symmetric in $a,b$ and
$\omega(a,\mu(b,c))$ is symmetric in $b,c$.
}
\end{remark}

If $V$ is a finite dimensional vector space, then we can associate with the Maurer-Cartan element,
$\Gamma:=\{\mu_n: V^{\ot n}\rar\ \susp^{n-2} V\}_{n\geq 1}\in C^1(V,V)$,
an element
$$
\Gamma^\circlearrowright:=\left\{\Gamma_{p,q}^\circlearrowright\in
{\it Hom}_1\left(((\downarrow V)^{\ot p})_{C_p} \ot
((\downarrow V)^{\ot q})_{C_q}, \bfk\right)
\right\}_{p,q\geq 1}\in
{\it Cyc}^1(V,\bfk),
$$
as follows,
\begin{eqnarray}
\label{cyclization}
\lefteqn{
\Gamma_{p,q}\left(v_1,\ldots, v_p, v_{p+1}, \ldots v_{p+q}\right):=}
\\
\nonumber
&=&
\sum_{i=0}^{p-1}\sum_{j=1}^{q-1}(-1)^{i(p+1) + j(q+1)}
{\it Tr}_{p+1}(\mu_{p+q+1})\left(v_{\zeta^l(1)}, \ldots, v_{\zeta^l(p)},
v_{\xi^j(p+1)}, \ldots,  v_{\xi^j(p+q)}\right)
\end{eqnarray}
where
$\zeta:=(12\ldots p)\in \bS_{p}$, $\xi:=((p+1)(p+2)\ldots (p+q))\in
\bS_{q}$,
${\it Tr}_{p+1}(\mu_{p+q+1}): V^{\ot p+q}\rar\ \susp^{n-2}\ \bfk$ is the trace
of the linear map $\mu_{p+q+1}: V^{p+q+1}\rar\ \susp^{p+q-2}\ V$
 with respect
to the $(p+1)$-st input,
and $v_1,\ldots,v_{p+q}$ are arbitrary elements of $V$.

\begin{prop-def}
For any finite-dimensional $\Ass_\infty$-algebra $(V, \Gamma)$ we have
$$
\delta_\Gamma \Gamma^\circlearrowright=0.
$$
The associated cohomology class $[\Gamma^\circlearrowright]\in
H^1({\it Cyc}^\bullet(V,\bfk))$ is called the {\em cyclic
characteristic class\/} of the $\Ass_\infty$-algebra $(V, \Gamma)$.
\end{prop-def}

\noindent
{\bf Sketch of a proof.}
 It follows from  definition (\ref{ass1}) of the differential in the operad $\Ass_\infty$ that
\Beqrn
\delta
\begin{xy}
 <0mm,0mm>*{\bullet};<0mm,0mm>*{}**@{},
 <0mm,0mm>*{};<-11mm,-5mm>*{}**@{-},
 <0mm,0mm>*{};<-8.5mm,-5mm>*{}**@{-},
 <0mm,0mm>*{};<-5mm,-5mm>*{...}**@{},
 <0mm,0mm>*{};<-2.5mm,-5mm>*{}**@{-},
   <0mm,0mm>*{};<-13mm,-7.9mm>*{^1}**@{},
   <0mm,0mm>*{};<-9mm,-7.9mm>*{^2}**@{},
   <0mm,0mm>*{};<-3mm,-7.9mm>*{^p}**@{},
<0mm,0mm>*{};<11mm,-5mm>*{}**@{-},
 <0mm,0mm>*{};<8.5mm,-5mm>*{}**@{-},
 <0mm,0mm>*{};<5.5mm,-5mm>*{...}**@{},
 <0mm,0mm>*{};<2.5mm,-5mm>*{}**@{-},
   <0mm,0mm>*{};<14mm,-7.9mm>*{^{p\hspace{-0.5mm}+\hspace{-0.5mm}q}}**@{},
   <0mm,0mm>*{};<3mm,-7.9mm>*{^{p\hspace{-0.5mm}+\hspace{-0.5mm}1}}**@{},
 <0mm,-10mm>*{};<0mm,10mm>*{}**@{-},
(0,10)*{}
   \ar@{->}@(ur,dr) (0,-10)*{}
 \end{xy}\ &=&
 \sum_{k=0}^{p-2}\sum_{l=2}^{p-k}
(-1)^{k+l(p-k-l)+1}
\begin{xy}
 <0mm,0mm>*{\bullet};<0mm,0mm>*{}**@{},
 <0mm,0mm>*{};<-15mm,-5mm>*{}**@{-},
 <0mm,0mm>*{};<-12mm,-5mm>*{...}**@{},
 <0mm,0mm>*{};<-10mm,-5mm>*{}**@{-},
 <0mm,0mm>*{};<-7mm,-5mm>*{}**@{-},
 <0mm,0mm>*{};<-4.5mm,-5mm>*{...}**@{},
 <0mm,0mm>*{};<-2.5mm,-5mm>*{}**@{-},
   <0mm,0mm>*{};<-17mm,-7.9mm>*{^1}**@{},
   <0mm,0mm>*{};<-11mm,-7.9mm>*{^k}**@{},
   <0mm,0mm>*{};<-3mm,-7.9mm>*{^p}**@{},
 <-7mm,-5.5mm>*{\bullet};
<-7mm,-5.5mm>*{};<-10mm,-11mm>*{}**@{-},
<-7mm,-5.5mm>*{};<-4mm,-11mm>*{}**@{-},
<-7mm,-5.5mm>*{};<-8mm,-11mm>*{}**@{-},
<-7mm,-13mm>*{_{k+1\, \dots\, k+l}},
<0mm,0mm>*{};<11mm,-5mm>*{}**@{-},
 <0mm,0mm>*{};<8.5mm,-5mm>*{}**@{-},
 <0mm,0mm>*{};<5.5mm,-5mm>*{...}**@{},
 <0mm,0mm>*{};<2.5mm,-5mm>*{}**@{-},
   <0mm,0mm>*{};<14mm,-7.9mm>*{^{p\hspace{-0.5mm}+\hspace{-0.5mm}q}}**@{},
   <0mm,0mm>*{};<3mm,-7.9mm>*{^{p\hspace{-0.5mm}+\hspace{-0.5mm}1}}**@{},
 <0mm,-10mm>*{};<0mm,10mm>*{}**@{-},
(0,10)*{}
   \ar@{->}@(ur,dr) (0,-10)*{}
 \end{xy}
 \\
&&
\hskip -10em
+\sum_{k=0}^{p}\sum_{l=p-1}^{p+q+1-k}
(-1)^{k+l(p+q-k)+1} \hskip -.5em
\begin{xy}
<0mm,0mm>*{\bullet},
<0mm,10mm>*{}**@{-},
<4mm,-7mm>*{^{1\  \dots \  k\qquad\ \ k+l+1\ \ \dots\ \  {p\hspace{-0.5mm}+\hspace{-0.5mm}q} }},
<-14mm,-5mm>*{}**@{-},
<-6mm,-5mm>*{}**@{-},
<20mm,-5mm>*{}**@{-},
<8mm,-5mm>*{}**@{-},
<0mm,-15mm>*{}**@{-},
<0mm,-5mm>*{\bullet};
<-5mm,-10mm>*{}**@{-},
<-2mm,-10mm>*{}**@{-},
<2mm,-10mm>*{}**@{-},
<5mm,-10mm>*{}**@{-},
<0.5mm,-12mm>*{_{k+1...\, p \hspace{-0.1mm} \ p+1... k+l}},
(0,10)*{}
   \ar@{->}@(ur,dr) (0,-15)*{}
\end{xy}
+\sum_{k=0}^{n-2}\sum_{l=2}^{n-k}
(-1)^{k+l(n-k-l)+1} \hskip -.5em
\begin{xy}
 <0mm,0mm>*{\bullet};<0mm,0mm>*{}**@{},
 <0mm,0mm>*{};<-11mm,-5mm>*{}**@{-},
 <0mm,0mm>*{};<-8.5mm,-5mm>*{}**@{-},
 <0mm,0mm>*{};<-5mm,-5mm>*{...}**@{},
 <0mm,0mm>*{};<-2.5mm,-5mm>*{}**@{-},
   <0mm,0mm>*{};<-13mm,-7.9mm>*{^1}**@{},
   <0mm,0mm>*{};<-9mm,-7.9mm>*{^2}**@{},
   <0mm,0mm>*{};<-3mm,-7.9mm>*{^p}**@{},
<0mm,0mm>*{};<10mm,-5mm>*{}**@{-},
 <0mm,0mm>*{};<16.5mm,-5mm>*{}**@{-},
 <0mm,0mm>*{};<6.6mm,-5mm>*{}**@{-},
 <0mm,0mm>*{};<4.6mm,-5mm>*{..}**@{},
 <0mm,0mm>*{};<13mm,-5mm>*{..}**@{},
 <0mm,0mm>*{};<2.5mm,-5mm>*{}**@{-},
 <0mm,0mm>*{};<7.5mm,-7.9mm>*{^{k}}**@{},
   <0mm,0mm>*{};<17mm,-7.9mm>*{^{p\hspace{-0.5mm}+\hspace{-0.5mm}q}}**@{},
   <0mm,0mm>*{};<3mm,-7.9mm>*{^{p\hspace{-0.5mm}+\hspace{-0.5mm}1}}**@{},
 <10.5mm,-5.5mm>*{\bullet};
 <10.5mm,-5.5mm>*{};<7.5mm,-10.5mm>*{}**@{-},
  <10.5mm,-5.5mm>*{};<13.5mm,-10.5mm>*{}**@{-},
   <10.5mm,-5.5mm>*{};<9.5mm,-10.5mm>*{}**@{-},
   <10.5mm,-13mm>*{_{k+1\, \dots\, k+l}},
 <0mm,-10mm>*{};<0mm,10mm>*{}**@{-},
(0,10)*{}
   \ar@{->}@(ur,dr) (0,-10)*{}
 \end{xy}.
\Eeqrn
Straightforward calculations show that upon ``(skew)cyclization"
 of both sides of the above equality as in (\ref{cyclization}) the middle sum
 (consisting of graphs with {\em two}\, cyclic vertices) vanishes. Then the remaining terms,
 when represented in a vector space $V$, assemble precisely into the required
 equation $\delta_\Gamma\Gamma^\circlearrowright=0$.
 \hfill $\Box$

\begin{theorem}
An\, $\Ass_\infty$-structure,
$\Gamma:=\{\mu_n: V^{\ot n}\rar\ \susp^{n-2}V\}_{n\geq 1}$,
 in a finite-dimensional vector
space $V$ can be extended to an $\Ass_\infty^\circlearrowright$-structure if and only if the cyclic
characteristic class $[\Gamma^\circlearrowright]$ vanishes.
\end{theorem}

\begin{proof}
If the cyclic characteristic class vanishes,
then $\Gamma^\circlearrowright=\delta_\Gamma \Theta$
for some $\Theta=\{ \mu_{p,q}: V^{\ot p+q}\rar\
\susp^{p+q}V\}_{p,q\geq 1}\in \Cyc  ^0(V,\bfk)$. Hence the datum
$\{\mu_\bullet, \mu_{\bullet,\bullet}\}$ makes
$V$ into an $\Ass_\infty^\circlearrowright$-algebra.

On the other hand, if $\{\mu_\bullet, \mu_{\bullet,\bullet}\}$ is an
$\Ass_\infty^\circlearrowright$-structure on $V$, then $\Theta:=\{
\mu_{p,q}: V^{\ot p+q}\rar\ \hbox{$\susp^{p+q}V$}\}_{p,q\geq 1}\in \Cyc  ^0(V,\bfk)$
satisfies the equation $\delta_\Gamma
\Theta=\Gamma^\circlearrowright$.
\end{proof}

A unusual feature of the situation described in Theorem~6.7.2. is
that there exists {\em precisely one\/} obstruction to extension of
$\Ass_\infty$-algebras into $\Ass_\infty^\wc$-algebras while in
obstruction theory one usually has to deal with an infinite
series of obstructions. The explanation is, however, simple -- the
wheeled operad $\Ass_\infty^\wc$ has only elements of genus $0$ or
$1$.

\begin{example}
{\rm
There is an interesting class of $\Ass_\infty$-algebras which
always extend into $\Ass_\infty^\circlearrowright$-algebras,
consisting of structures $\Gamma:=\{\mu_n: V^{\ot
  n}\rar\ \susp^{n-2}V\}_{n\geq 1}$ such that the antisymmetrization
of the trace ${\it Tr}_{p+1}(\mu_{p+q+1}): V^{\ot p+q}\rar\
\susp^{n-2}\ \bfk$ in the first $p$ and the
last $q$ variables vanishes for each $p,q \geq 1$, $n = p+q$.
We suggest to call these
algebras {\em traceless $\Ass_\infty$-algebras\/}.
}
\end{example}

\section{Wheeled $\mathsf{Com}$}
\label{Terezka}

In this section we prove Theorem~B of the introduction and show
thereby that the operad
$\Com$ for commutative associative algebras is wheeled Koszul.
We also demonstrate directly  that $\Com$
is not stably Koszul in the sense of
Definition~\ref{snad_to_nebude_nic_vazneho} by calculating, in
Theorem~\ref{thm-comm-stable},
the homology of the
wheeled completion $(\Com_\infty)^\wc$.

\subsection{Resolution of wheeled $\Com$}
The operad  $\Com$ is the quotient
$$
\Com = \freeOP(E)/(R) ,
$$
of the free operad $\freeOP(E)$ on the $\bS$-module $E$ with
\begin{align*}
E(2) & := \id_2 \qquad \mbox{spanned by }
\begin{xy}
<0mm,0mm>*{\bullet},
<0mm,5mm>*{}**@{-},
<-3mm,-5mm>*{}**@{-},
<3mm,-5mm>*{}**@{-},
<-3mm,-7mm>*{_1},
<3mm,-7mm>*{_2},
\end{xy}
=
\begin{xy}
<0mm,0mm>*{\bullet},
<0mm,5mm>*{}**@{-},
<-3mm,-5mm>*{}**@{-},
<3mm,-5mm>*{}**@{-},
<-3mm,-7mm>*{_2},
<3mm,-7mm>*{_1},
\end{xy}
\\
E(n)& :=0, \qquad n \not= 1
\end{align*}
modulo the ideal generated by the relations
$$
R:
\begin{xy}
<0mm,0mm>*{\bullet},
<0mm,5mm>*{}**@{-},
<-3mm,-5mm>*{}**@{-},
<3mm,-5mm>*{}**@{-},
<-3mm,-5mm>*{\bullet};
<-6mm,-10mm>*{}**@{-},
<0mm,-10mm>*{}**@{-},
<-6mm,-12mm>*{_{\sigma(1)}},
<0mm,-12mm>*{_{\sigma(2)}},
<3mm,-7mm>*{_{\sigma(3)}},
\end{xy}
-
\begin{xy}
<0mm,0mm>*{\bullet},
<0mm,5mm>*{}**@{-},
<-3mm,-5mm>*{}**@{-},
<3mm,-5mm>*{}**@{-},
<-3mm,-5mm>*{\bullet};
<-6mm,-10mm>*{}**@{-},
<0mm,-10mm>*{}**@{-},
<-6mm,-12mm>*{_{\sigma(1)}},
<0mm,-12mm>*{_{\sigma(3)}},
<3mm,-7mm>*{_{\sigma(2)}},
\end{xy}
=0, \qquad \forall \sigma\in \bS_3.
$$

Recall that the minimal model
$\mathsf{Com}_\infty$ of the operad $\mathsf{Com}$ is generated
by standard associative $(1,n)$-corollas
in degree $2-n$, $n \geq 2$,
modulo the shuffle relations as in Theorem~B(i),
with the differential given by~(\ref{comm-diff}).

\begin{theorem}\label{thm-comm-stable}
The dimension of the homology
of $\left(\mathsf{Com}_\infty\right)^\circlearrowright(0,n)$
in degree $- k\leq 0$ equals~$\binom{n-1}{k}$.
\end{theorem}

In the above theorem which we prove in
Subsection~\ref{jsem_zvedav_kdy_a_jestli_vubec_napise}
we used the convention that $\binom nk = 0$ for $k > n$.
Theorem~\ref{thm-comm-stable} shows that,
in order to obtain a~minimal model of
$\mathsf{Com}^\circlearrowright$, one must add new generators to
$(\mathsf{Com}_\infty)^\wc$ that kill homology classes in
$H^{<0}((\mathsf{Com}_\infty)^\wc, \pa)$, then add generators that
kill relations among these new generators, etc. The result is
described in Theorem~B which we prove later in this section. We also prove:

\begin{theorem}
\label{thm-Com-wheeled-Koszul}
The operad $\mathsf{Com}$ is wheeled Koszul.
\end{theorem}

Theorem~\ref{thm-Com-wheeled-Koszul} together with the uniqueness, up
to isomorphism, of minimal models imply that the minimal model of
$\mathsf{Com}^\circlearrowright$ described in Theorem~B  coincides with the
wheeled cobar of the wheeled quadratic dual $(\Com^\wc)^!$
constructed in Example~\ref{Kralicek}.

Theorem~B
and Theorem~\ref{thm-Com-wheeled-Koszul} are proved in
Subsection~\ref{zacinaji_Zbraslavice}. Their proofs are based on an analysis
of $(\Omega(B(\mathsf{Com})))^\circlearrowright$ (the wheeled
 completion of the double bar complex for the ordinary
$\mathsf{Com}$) and
$\Omega^\circlearrowright(B^\circlearrowright(\mathsf{Com}^\circlearrowright))$
(the double bar complex of $\mathsf{Com}^\circlearrowright$ in the
category of wheeled operads).

\subsection{Double bar complex of $\mathsf{Com}$}
\label{ordinary-double-bar-Com}

Let us show how the description of $\mathsf{Com}_\infty$ presented
above can be read off from the ordinary double bar complex
$\Omega(B(\mathsf{Com}))$. Recall that $\Omega(B(\mathsf{Com}))(n)$ is
generated by trees with one output leg and $n$ input legs,
with vertices the standard commutative corollas with one
outgoing edge and $\geq 2$ symmetric incoming edges. The set of
edges $Edg(G)$ of each tree $G$ is two-colored (in pictures, black
edges are ordinary ones, and white edges are doubled ones).

In the construction below, we use the same conventions as in the proof of
Theorem~\ref{ll}, in particular, the bigrading~(\ref{Jituska})
is the one used in the proof of Theorem~\ref{ll}.
The space $\Omega(B(\mathsf{Com}))(n)$ is bigraded,\label{whops}
with the bidegree
given by
\begin{equation}
\label{Jituska}
{\rm bideg\/}(G) = (\deg_1(G), \deg_2(G)) := (-e(G),e^0(G)),
\end{equation}
where $e(G)$ is the number of internal edges of the underlying
graph $G$ and $e^\circ
(G)$ the number of white edges of $G$.
There are two differentials, $(1,0)$-and $(0,1)$-ones,
denoted by $\partial_1$ and $\partial_2$, respectively.
The first one contracts black edges, and the second one changes the
color of edges from black to white. For example,
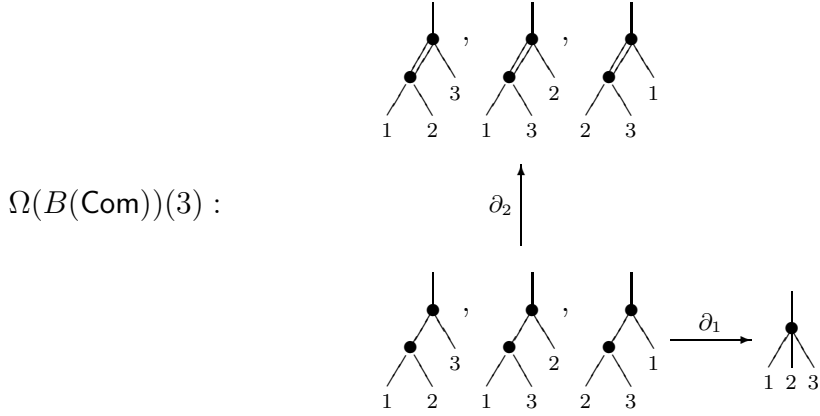
\begin{center}
\setlength{\unitlength}{3.6cm}
\begin{picture}(2,1.6)(-2,-0.3)
\put(0,0){\makebox(0,0){
\begin{xy}
<0mm,0mm>*{\bullet},
<0mm,5mm>*{}**@{-},
<-3mm,-5mm>*{}**@{-},
<0mm,-5mm>*{}**@{-},
<3mm,-5mm>*{}**@{-},
<-3mm,-7mm>*{_1},
<0mm,-7mm>*{_2},
<3mm,-7mm>*{_3},
\end{xy}
}}
\put(-1,0){\makebox(0,0){
\begin{xy}
<0mm,0mm>*{\bullet},
<0mm,5mm>*{}**@{-},
<-3mm,-5mm>*{}**@{-},
<3mm,-5mm>*{}**@{-},
<-3mm,-5mm>*{\bullet};
<-6mm,-10mm>*{}**@{-},
<0mm,-10mm>*{}**@{-},
<-6mm,-12mm>*{_{1}},
<0mm,-12mm>*{_{2}},
<3mm,-7mm>*{_{3}},
\end{xy},
\begin{xy}
<0mm,0mm>*{\bullet},
<0mm,5mm>*{}**@{-},
<-3mm,-5mm>*{}**@{-},
<3mm,-5mm>*{}**@{-},
<-3mm,-5mm>*{\bullet};
<-6mm,-10mm>*{}**@{-},
<0mm,-10mm>*{}**@{-},
<-6mm,-12mm>*{_{1}},
<0mm,-12mm>*{_{3}},
<3mm,-7mm>*{_{2}},
\end{xy},
\begin{xy}
<0mm,0mm>*{\bullet},
<0mm,5mm>*{}**@{-},
<-3mm,-5mm>*{}**@{-},
<3mm,-5mm>*{}**@{-},
<-3mm,-5mm>*{\bullet};
<-6mm,-10mm>*{}**@{-},
<0mm,-10mm>*{}**@{-},
<-6mm,-12mm>*{_{2}},
<0mm,-12mm>*{_{3}},
<3mm,-7mm>*{_{1}},
\end{xy}
}}
\put(-1,1){\makebox(0,0){
\begin{xy}
<0mm,0mm>*{\bullet},
<0mm,5mm>*{}**@{-},
<-3mm,-5mm>*{}**@{=},
<3mm,-5mm>*{}**@{-},
<-3mm,-5mm>*{\bullet};
<-6mm,-10mm>*{}**@{-},
<0mm,-10mm>*{}**@{-},
<-6mm,-12mm>*{_{1}},
<0mm,-12mm>*{_{2}},
<3mm,-7mm>*{_{3}},
\end{xy},
\begin{xy}
<0mm,0mm>*{\bullet},
<0mm,5mm>*{}**@{-},
<-3mm,-5mm>*{}**@{=},
<3mm,-5mm>*{}**@{-},
<-3mm,-5mm>*{\bullet};
<-6mm,-10mm>*{}**@{-},
<0mm,-10mm>*{}**@{-},
<-6mm,-12mm>*{_{1}},
<0mm,-12mm>*{_{3}},
<3mm,-7mm>*{_{2}},
\end{xy},
\begin{xy}
<0mm,0mm>*{\bullet},
<0mm,5mm>*{}**@{-},
<-3mm,-5mm>*{}**@{=},
<3mm,-5mm>*{}**@{-},
<-3mm,-5mm>*{\bullet};
<-6mm,-10mm>*{}**@{-},
<0mm,-10mm>*{}**@{-},
<-6mm,-12mm>*{_{2}},
<0mm,-12mm>*{_{3}},
<3mm,-7mm>*{_{1}},
\end{xy}
}}
\put(-1,0.35){\vector(0,1){.3}}
\put(-1.03,0.5){\makebox(0,0)[r]{\scriptsize $\partial_2$}}
\put(.05,0){
\put(-.5,0){\vector(1,0){.3}}
\put(-.35,0.03){\makebox(0,0)[b]{\scriptsize $\partial_1$}}}
\put(-2.5,0.5){\makebox(0,0){$\Omega(B(\mathsf{Com}))(3):$}}
\end{picture}
\end{center}

Since $\mathsf{Com}$ is Koszul, its $\partial_1$-cohomology of
$\Omega(B(\mathsf{Com}))$ is concentrated in bidegrees $(2-n,*)$ and
is isomorphic to $\mathsf{Com}_\infty$ (the differential of
$\mathsf{Com}_\infty$ is induced by $\partial_2$). The description of
$\mathsf{Com}_\infty$ given above uses associative corollas satisfying
shuffle relations. They are represented in $\Omega(B(\mathsf{Com}))$
in the following way. Consider the sum of all planar binary trees with
only black edges with $1$ output leg and $n$ input legs marked by
$\sigma(1),\dots,\sigma(n)$ from the left to the right.  There is,
modulo an overall sign, a unique
choice of signs of the graphs such that this
sum is a $\partial_1$-closed element in $\Omega(B(\mathsf{Com})(n)$.
It represents the associative corolla with $n$ inputs marked by
$\sigma(1),\dots,\sigma(n)$ from the left to the right. One can check that
these corollas satisfy the shuffle relations, and that $\partial$ is
induced by $\partial_2$.

Here are few examples (we omit the signs; associative corollas are on
the left, commutative corollas are on the right):
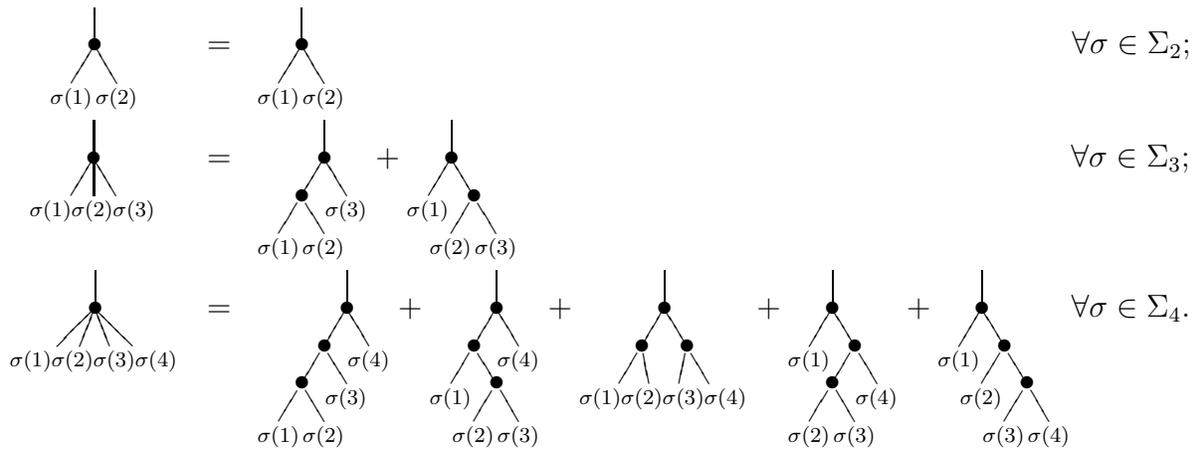
\begin{eqnarray*}
\begin{xy}
<0mm,0mm>*{\bullet},
<0mm,5mm>*{}**@{-},
<-3mm,-5mm>*{}**@{-},
<3mm,-5mm>*{}**@{-},
<-3mm,-7mm>*{_{\sigma(1)}},
<3mm,-7mm>*{_{\sigma(2)}},
\end{xy}
\hskip 5.5mm& = &
\begin{xy}
<0mm,0mm>*{\bullet},
<0mm,5mm>*{}**@{-},
<-3mm,-5mm>*{}**@{-},
<3mm,-5mm>*{}**@{-},
<-3mm,-7mm>*{_{\sigma(1)}},
<3mm,-7mm>*{_{\sigma(2)}},
\end{xy}
\hskip 9.65cm \forall\sigma\in\bS_2;
\\
\begin{xy}
<0mm,0mm>*{\bullet},
<0mm,5mm>*{}**@{-},
<-3mm,-5mm>*{}**@{-},
<0mm,-5mm>*{}**@{-},
<3mm,-5mm>*{}**@{-},
<0mm,-7mm>*{_{\sigma(1)\sigma(2)\sigma(3)}},
\end{xy}
\hskip 3mm & = &
\begin{xy}
<0mm,0mm>*{\bullet},
<0mm,5mm>*{}**@{-},
<-3mm,-5mm>*{}**@{-},
<3mm,-5mm>*{}**@{-},
<-3mm,-5mm>*{\bullet};
<-6mm,-10mm>*{}**@{-},
<0mm,-10mm>*{}**@{-},
<-6mm,-12mm>*{_{\sigma(1)}},
<0mm,-12mm>*{_{\sigma(2)}},
<3mm,-7mm>*{_{\sigma(3)}},
\end{xy}
+
\begin{xy}
<0mm,0mm>*{\bullet},
<0mm,5mm>*{}**@{-},
<-3mm,-5mm>*{}**@{-},
<3mm,-5mm>*{}**@{-},
<3mm,-5mm>*{\bullet};
<6mm,-10mm>*{}**@{-},
<0mm,-10mm>*{}**@{-},
<6mm,-12mm>*{_{\sigma(3)}},
<0mm,-12mm>*{_{\sigma(2)}},
<-3mm,-7mm>*{_{\sigma(1)}},
\end{xy}
\hskip 7.35cm \forall\sigma\in\bS_3;
\\
\begin{xy}
<0mm,0mm>*{\bullet},
<0mm,5mm>*{}**@{-},
<-5mm,-5mm>*{}**@{-},
<-2mm,-5mm>*{}**@{-},
<2mm,-5mm>*{}**@{-},
<5mm,-5mm>*{}**@{-},
<0mm,-7mm>*{_{\sigma(1)\sigma(2)\sigma(3)\sigma(4)}},
\end{xy}
& = &
\begin{xy}
<0mm,0mm>*{\bullet},
<0mm,5mm>*{}**@{-},
<-3mm,-5mm>*{}**@{-},
<3mm,-5mm>*{}**@{-},
<-3mm,-5mm>*{\bullet};
<-6mm,-10mm>*{}**@{-},
<0mm,-10mm>*{}**@{-},
<-6mm,-10mm>*{\bullet};
<-9mm,-15mm>*{}**@{-},
<-3mm,-15mm>*{}**@{-},
<-9mm,-17mm>*{_{\sigma(1)}},
<-3mm,-17mm>*{_{\sigma(2)}},
<0mm,-12mm>*{_{\sigma(3)}},
<3mm,-7mm>*{_{\sigma(4)}},
\end{xy}
+
\begin{xy}
<0mm,0mm>*{\bullet},
<0mm,5mm>*{}**@{-},
<-3mm,-5mm>*{}**@{-},
<3mm,-5mm>*{}**@{-},
<-3mm,-5mm>*{\bullet};
<-6mm,-10mm>*{}**@{-},
<0mm,-10mm>*{}**@{-},
<0mm,-10mm>*{\bullet};
<-3mm,-15mm>*{}**@{-},
<3mm,-15mm>*{}**@{-},
<-6mm,-12mm>*{_{\sigma(1)}},
<-3mm,-17mm>*{_{\sigma(2)}},
<3mm,-17mm>*{_{\sigma(3)}},
<3mm,-7mm>*{_{\sigma(4)}},
\end{xy}
+
\begin{xy}
<0mm,0mm>*{\bullet},
<0mm,5mm>*{}**@{-},
<-3mm,-5mm>*{}**@{-},
<3mm,-5mm>*{}**@{-},
<-3mm,-5mm>*{\bullet};
<-6mm,-10mm>*{}**@{-},
<-2mm,-10mm>*{}**@{-},
<3mm,-5mm>*{\bullet};
<6mm,-10mm>*{}**@{-},
<2mm,-10mm>*{}**@{-},
<0mm,-12mm>*{_{\sigma(1)\sigma(2)\sigma(3)\sigma(4)}},
\end{xy}+
\begin{xy}
<0mm,0mm>*{\bullet},
<0mm,5mm>*{}**@{-},
<-3mm,-5mm>*{}**@{-},
<3mm,-5mm>*{}**@{-},
<3mm,-5mm>*{\bullet};
<6mm,-10mm>*{}**@{-},
<0mm,-10mm>*{}**@{-},
<0mm,-10mm>*{\bullet};
<-3mm,-15mm>*{}**@{-},
<3mm,-15mm>*{}**@{-},
<6mm,-12mm>*{_{\sigma(4)}},
<3mm,-17mm>*{_{\sigma(3)}},
<-3mm,-17mm>*{_{\sigma(2)}},
<-3mm,-7mm>*{_{\sigma(1)}},
\end{xy}
+
\begin{xy}
<0mm,0mm>*{\bullet},
<0mm,5mm>*{}**@{-},
<-3mm,-5mm>*{}**@{-},
<3mm,-5mm>*{}**@{-},
<3mm,-5mm>*{\bullet};
<6mm,-10mm>*{}**@{-},
<0mm,-10mm>*{}**@{-},
<6mm,-10mm>*{\bullet};
<9mm,-15mm>*{}**@{-},
<3mm,-15mm>*{}**@{-},
<9mm,-17mm>*{_{\sigma(4)}},
<3mm,-17mm>*{_{\sigma(3)}},
<0mm,-12mm>*{_{\sigma(2)}},
<-3mm,-7mm>*{_{\sigma(1)}},
\end{xy}
 \forall\sigma\in\bS_4.
\end{eqnarray*}

\subsection{Wheeled double bar complexes}\label{wheeled-db-com}

We consider three wheeled dg-operads whose
operadic parts coincide with $\Omega(B(\mathsf{Com}))$.
The first one is its wheeled
completion $(\Omega(B(\mathsf{Com}))^\circlearrowright$.
The space $(\Omega(B(\mathsf{Com}))^\circlearrowright(0,n)$ is
spanned by all graphs with one oriented wheel and $n$ input
legs. There are black and white edges, as before, but we require
that there is at least one white edge in the wheel. Bigrading and
differentials are extended
from $\Omega(B(\mathsf{Com}))$.

The second one is the extended wheeled completion denoted by
$(\Omega(B(\mathsf{Com}))^\circlearrowright_x$ defined as
the extension of
$(\Omega(B(\mathsf{Com}))^\circlearrowright$ by graphs that have all
edges in the wheel black.  Notice that
$(\Omega(B(\mathsf{Com}))^\circlearrowright_x$ contains also graphs
with black loops (one-edge wheels).  Applying $\partial_1$ to such
graphs would require to contract such loops. But this is impossible in
$(\Omega(B(\mathsf{Com}))^\circlearrowright_x$ and we extend the
differential by postulating that contracting such a loop gives zero.

We can consider one more wheeled dg-operad where this
contraction is defined. We add graphs with one-edge wheel considered
as a ``virtual'' one (indicated by dots in pictures). Obviously,
what we get (with naturally extended gradings and differentials) is
the wheeled double bar complex
$\Omega^\circlearrowright(B^\circlearrowright(\mathsf{Com}^\circlearrowright))$.
\begin{figure}[t]
\begin{center}
\setlength{\unitlength}{3.2cm}
\begin{picture}(2,2.6)(-2.5,-0.3)
\put(0,0){\makebox(0,0){
\begin{xy}
<0mm,0mm>*{\bullet},
<-3mm,-5mm>*{}**@{-},
<3mm,-5mm>*{}**@{-},
<-3mm,-7mm>*{_1},
<3mm,-7mm>*{_2},
(0,0)*{} \ar@{.}@(u,r) (0,0)*{}
\end{xy}
}}
\put(-1,0){\makebox(0,0){
\begin{xy}
<0mm,0mm>*{\bullet},
<-3mm,-5mm>*{}**@{-},
<3mm,-5mm>*{}**@{-},
<-3mm,-7mm>*{_1},
<3mm,-7mm>*{_2},
(0,0)*{} \ar@{-}@(u,r) (0,0)*{}
\end{xy},
\begin{xy}
<0mm,0mm>*{\bullet},
<-3mm,-5mm>*{}**@{-},
<-3mm,-5mm>*{\bullet};
<-6mm,-10mm>*{}**@{-},
<0mm,-10mm>*{}**@{-},
<-6mm,-12mm>*{_{1}},
<0mm,-12mm>*{_{2}},
(0,0)*{} \ar@{.}@(u,r) (0,0)*{}
\end{xy}
}}
\put(-2.5,0){\makebox(0,0){
\begin{xy}
<0mm,0mm>*{\bullet},
<-3mm,-5mm>*{}**@{-},
<-3mm,-5mm>*{\bullet};
<-6mm,-10mm>*{}**@{-},
<0mm,-10mm>*{}**@{-},
<-6mm,-12mm>*{_{1}},
<0mm,-12mm>*{_{2}},
(0,0)*{} \ar@{-}@(u,r) (0,0)*{}
\end{xy},
\begin{xy}
<0mm,0mm>*{\bullet},
<-3mm,-5mm>*{}**@{-},
<3mm,-5mm>*{}**@{-},
<3mm,-5mm>*{\bullet};
<0mm,-10mm>*{}**@{-},
<0mm,-12mm>*{_2},
<-3mm,-7mm>*{_1},
(0,0)*{} \ar@{-}@(u,r) (3,-5)*{}
\end{xy}
}}
\put(-.9,1){\makebox(0,0){
\begin{xy}
<0mm,0mm>*{\bullet},
<-3mm,-5mm>*{}**@{-},
<3mm,-5mm>*{}**@{-},
<-3mm,-7mm>*{_1},
<3mm,-7mm>*{_2},
(0,0)*{} \ar@{=}@(u,r) (0,0)*{}
\end{xy},
\begin{xy}
<0mm,0mm>*{\bullet},
<-3mm,-5mm>*{}**@{=},
<-3mm,-5mm>*{\bullet};
<-6mm,-10mm>*{}**@{-},
<0mm,-10mm>*{}**@{-},
<-6mm,-12mm>*{_{1}},
<0mm,-12mm>*{_{2}},
(0,0)*{} \ar@{.}@(u,r) (0,0)*{}
\end{xy}
}}
\put(-2.6,1){\makebox(0,0){
\begin{xy}
<0mm,0mm>*{\bullet},
<-3mm,-5mm>*{}**@{=},
<-3mm,-5mm>*{\bullet};
<-6mm,-10mm>*{}**@{-},
<0mm,-10mm>*{}**@{-},
<-6mm,-12mm>*{_{1}},
<0mm,-12mm>*{_{2}},
(0,0)*{} \ar@{}@(u,r) (0,0)*{}
\end{xy},
\begin{xy}
<0mm,0mm>*{\bullet},
<-3mm,-5mm>*{}**@{-},
<-3mm,-5mm>*{\bullet};
<-6mm,-10mm>*{}**@{-},
<0mm,-10mm>*{}**@{-},
<-6mm,-12mm>*{_{1}},
<0mm,-12mm>*{_{2}},
(0,0)*{} \ar@{=}@(u,r) (0,0)*{}
\end{xy},
\begin{xy}
<0mm,0mm>*{\bullet},
<-3mm,-5mm>*{}**@{-},
<3mm,-5mm>*{}**@{=},
<3mm,-5mm>*{\bullet};
<0mm,-10mm>*{}**@{-},
<0mm,-12mm>*{_2},
<-3mm,-7mm>*{_1},
(0,0)*{} \ar@{}@(u,r) (3,-5)*{}
\end{xy},
\begin{xy}
<0mm,0mm>*{\bullet},
<-3mm,-5mm>*{}**@{-},
<3mm,-5mm>*{}**@{=},
<3mm,-5mm>*{\bullet};
<0mm,-10mm>*{}**@{-},
<0mm,-12mm>*{_1},
<-3mm,-7mm>*{_2},
(0,0)*{} \ar@{}@(u,r) (3,-5)*{}
\end{xy}
}}
\put(-2.5,2){\makebox(0,0){
\begin{xy}
<0mm,0mm>*{\bullet},
<-3mm,-5mm>*{}**@{=},
<-3mm,-5mm>*{\bullet};
<-6mm,-10mm>*{}**@{-},
<0mm,-10mm>*{}**@{-},
<-6mm,-12mm>*{_{1}},
<0mm,-12mm>*{_{2}},
(0,0)*{} \ar@{=}@(u,r) (0,0)*{}
\end{xy},
\begin{xy}
<0mm,0mm>*{\bullet},
<-3mm,-5mm>*{}**@{-},
<3mm,-5mm>*{}**@{=},
<3mm,-5mm>*{\bullet};
<0mm,-10mm>*{}**@{-},
<0mm,-12mm>*{_2},
<-3mm,-7mm>*{_1},
(0,0)*{} \ar@{=}@(u,r) (3,-5)*{}
\end{xy}
}}
\put(-.1,0){
\put(-.5,0){\vector(1,0){.3}}
\put(-.35,0.03){\makebox(0,0)[b]{\scriptsize $\partial_1$}}}
\put(-1.35,0){
\put(-0.5,0){\vector(1,0){.3}}
\put(-.35,0.03){\makebox(0,0)[b]{\scriptsize $\partial_1$}}}
\put(-1.1,1.05){
\put(-.5,0){\vector(1,0){.3}}
\put(-.35,0.03){\makebox(0,0)[b]{\scriptsize $\partial_1$}}}
\put(.05,0){
\put(-1,0.35){\vector(0,1){.3}}
\put(-1.03,0.5){\makebox(0,0)[r]{\scriptsize $\partial_2$}}}
\put(-1.45,0){
\put(-1,0.35){\vector(0,1){.3}}
\put(-1.03,0.5){\makebox(0,0)[r]{\scriptsize $\partial_2$}}}
\put(-1.45,1){
\put(-1,0.35){\vector(0,1){.3}}
\put(-1.03,0.5){\makebox(0,0)[r]{\scriptsize $\partial_2$}}}
\end{picture}
\end{center}
\caption{\label{ze_uz_bych_koncil?}%
The bicomplex $\Omega^\circlearrowright(B^\circlearrowright
(\mathsf{Com}^\circlearrowright))(0,2)$.}
\end{figure}

The example of
$\Omega^\circlearrowright(B^\circlearrowright
(\mathsf{Com}^\circlearrowright))(0,2)$ is given in
Figure~\ref{ze_uz_bych_koncil?}.
The bicomplex $(\Omega(B(\mathsf{Com}))^\circlearrowright_x(0,2)$ is
a quotient
of the above one consisting of graphs with no
virtual loops;
$(\Omega(B(\mathsf{Com}))^\circlearrowright(0,2)$ is a
sub-bicomplex without virtual loops and fully black wheels.

\begin{proposition}\label{x-acyclic}
The complex $(\Omega(B(\mathsf{Com}))^\circlearrowright_x(0,n)$ is,
for each $n\geq 1$, acyclic.
\end{proposition}

\begin{proof}
We use the same reasoning as in the proof of Theorem~\ref{ll}.
As before, $(\Omega(B(\mathsf{Com}))^\circlearrowright_x(0,n)$ with
differential $\partial_2$ splits into the direct sum over isomorphism
classes of graphs of complexes of exterior algebras~(\ref{prod}) on edges
of each type of graph.  It is acyclic if there is at least one edge in
the underlying graph. But here each graph has at least one
edge. Therefore, $(\Omega(B(\mathsf{Com}))^\circlearrowright_x(0,n)$
is acyclic.
\end{proof}

\subsection{Proof of Theorem~\protect\ref{thm-comm-stable}}
\label{jsem_zvedav_kdy_a_jestli_vubec_napise}

Notice that the homology of
$(\mathsf{Com}_\infty)^\circlearrowright$ coincides with the
homology of $(\Omega(B(\mathsf{Com})))^\circlearrowright$.  This
follows from the fact that the free wheeled operad functor
is a polynomial functor on the category of $\Sigma$-modules and that
the wheeled  completion of a free operad is a free wheeled
operad, so the wheeled completions of free operads generated by
quasi-isomorphic dg-$\Sigma$-modules are quasi-isomorphic.

Since $(\Omega(B(\mathsf{Com})))^\circlearrowright$ is a
subcomplex of the acyclic complex
$(\Omega(B(\mathsf{Com})))^\circlearrowright_x$, it is enough to study
their quotient,
$$
\mathcal{F}_{x/c}:=\frac{(\Omega(B(\mathsf{Com})))^\circlearrowright_x}
{(\Omega(B(\mathsf{Com})))^\circlearrowright}.
$$ Note that $\mathcal{F}_{x/c}$ is generated by graphs with fully
black wheels. The first differential, $\partial_1$, is the same as
before and contracts black edges. The second differential,
$\partial_2$, changes the color of edges \emph{outside the wheel} from
black to white.

The $\partial_2$-cohomology of $\mathcal{F}_{x/c}$ splits into the
direct sum of complexes of exterior algebras on edges, as in the
proofs of Theorem~\ref{ll} and Proposition~\ref{x-acyclic}. But here we
consider only the edges outside the wheel, therefore  generators of the
$\partial_2$-cohomology are graphs with no edges outside the wheel.

It follows from a short exact sequence argument and the identification
above that the cohomology of
$(\mathsf{Com}_\infty)^\circlearrowright(0,n)$ shifted by $1$
equals to the cohomology of the graph complex $\mathcal{W}_n$ spanned
by wheels with $n$ numbered legs, with no edge outside the wheel, and
with the differential $\partial_1$ contracting edges. For example,
$\mathcal{W}_2$:
\[
0\longrightarrow
\begin{xy}
<0mm,0mm>*{\bullet},
<-3mm,-5mm>*{}**@{-},
<3mm,-5mm>*{}**@{-},
<3mm,-5mm>*{\bullet};
<0mm,-10mm>*{}**@{-},
<0mm,-12mm>*{_2},
<-3mm,-7mm>*{_1},
(0,0)*{} \ar@{-}@(u,r) (3,-5)*{}
\end{xy}
\stackrel{\partial_1}\longrightarrow
\begin{xy}
<0mm,0mm>*{\bullet},
<-3mm,-5mm>*{}**@{-},
<3mm,-5mm>*{}**@{-},
<-3mm,-7mm>*{_1},
<3mm,-7mm>*{_2},
(0,0)*{} \ar@{-}@(u,r) (0,0)*{}
\end{xy}
\longrightarrow 0.
\]
The grading, given by the number of edges, comes from the double
bar complex. Observe that $\mathcal{W}_n$ is concentrated in degrees
$-1,\ldots,-n$.  To describe the differential $\partial_1$, it is
convenient to consider graphs $G$ in $\mathcal{W}_n$ twisted by the
determinant space $\mathrm{Det}(Edg(G)):
=\Lambda^{|Edg(G)|}\left(\langle Edg(G)\rangle\right)$. Then
\[
\partial_1\colon \mathrm{Det}(Edg(G))\otimes G \mapsto \sum_{e\in
Edg(G)} \frac{\partial}{\partial e} \mathrm{Det}(Edg(G)) \otimes
\partial_e(G),
\]
where $\partial_e(G)$ is, as before, the graph $G/e$ obtained by
contraction of the edge $e$.  If $e$ is a loop (i.~e., $G$ is one-edge
graph in degree $-1$), then $\partial_e(G)=0$.  In particular,
$\partial_1\equiv 0$ in $\mathcal{W}_1$.

Theorem~\ref{thm-comm-stable} now follows from the following
proposition whose proof is an exercise on cyclic cohomology:

\begin{proposition}
The dimension of  $H^{-k}(\mathcal{W}_n)$ equals $\binom{n-1}{k-1}$,
$-k \leq -1$.
\end{proposition}

The relation to cyclic cohomology is the following.  The complex
defining cyclic cohomology of the polynomial algebra in variables
$x_1,\dots, x_n$ contains a subcomplex generated by cyclic expressions
where each $x_i$ appears exactly once.  This subcomplex is isomorphic
to $\mathcal{W}_n$.

\subsection{Proofs of Theorem~B and Theorem~\ref{thm-Com-wheeled-Koszul}}
\label{zacinaji_Zbraslavice}

We prove Theorem~\ref{thm-Com-wheeled-Koszul} first. Let us consider
$B^\wc(\mathsf{Com}^\circlearrowright)(0,n)$. It is the complex of
graphs with $n$ input legs and no output legs whose vertices are
decorated by elements in $\mathsf{Com}^\circlearrowright$, graded by
the number of edges. The differential $\partial$ contracts the edges,
introducing appropriate signs. For example,
\[
B^\wc(\mathsf{Com}^\circlearrowright)(0,2):\quad
\begin{array}{ccccc}
\begin{xy}
<0mm,0mm>*{\bullet},
<-3mm,-5mm>*{}**@{-},
<-3mm,-5mm>*{\bullet};
<-6mm,-10mm>*{}**@{-},
<0mm,-10mm>*{}**@{-},
<-6mm,-12mm>*{_{1}},
<0mm,-12mm>*{_{2}},
(0,0)*{} \ar@{-}@(u,r) (0,0)*{}
\end{xy},
\begin{xy}
<0mm,0mm>*{\bullet},
<-3mm,-5mm>*{}**@{-},
<3mm,-5mm>*{}**@{-},
<3mm,-5mm>*{\bullet};
<0mm,-10mm>*{}**@{-},
<0mm,-12mm>*{_2},
<-3mm,-7mm>*{_1},
(0,0)*{} \ar@{-}@(u,r) (3,-5)*{}
\end{xy}
& \stackrel{\partial}\longrightarrow &
\begin{xy}
<0mm,0mm>*{\bullet},
<-3mm,-5mm>*{}**@{-},
<3mm,-5mm>*{}**@{-},
<-3mm,-7mm>*{_1},
<3mm,-7mm>*{_2},
(0,0)*{} \ar@{-}@(u,r) (0,0)*{}
\end{xy},
\begin{xy}
<0mm,0mm>*{\bullet},
<-3mm,-5mm>*{}**@{-},
<-3mm,-5mm>*{\bullet};
<-6mm,-10mm>*{}**@{-},
<0mm,-10mm>*{}**@{-},
<-6mm,-12mm>*{_{1}},
<0mm,-12mm>*{_{2}},
\end{xy}
& \stackrel{\partial}\longrightarrow &
\begin{xy}
<0mm,0mm>*{\bullet},
<-3mm,-5mm>*{}**@{-},
<3mm,-5mm>*{}**@{-},
<-3mm,-7mm>*{_1},
<3mm,-7mm>*{_2},
\end{xy}
\end{array}
\]
By definition of wheeled Koszulness,
Theorem~\ref{thm-Com-wheeled-Koszul} is equivalent to the statement
that the cohomology of $B^\wc(\mathsf{Com}^\circlearrowright)(0,n)$
is zero for $n=1$ and concentrated in degree $-n$ for $n\geq 2$, i.e., that it is generated by trivalent
graphs of genus $1$.

For $n=1$ it follows from a direct inspection that
$B^\wc(\mathsf{Com}^\circlearrowright)(0,1)$ is acyclic. If $n\geq
2$, we introduce a filtration on
$B^\wc(\mathsf{Com}^\circlearrowright)(0,n)$:
\[
0=F_{-1}\subset
F_0\subset \dots\subset
F_n=B^\wc(\mathsf{Com}^\circlearrowright)(0,n),
\]
in which $F_i$ is the space of graphs with $\leq i$ edges in the
wheel.  For example, $F_0$ is the space of graphs with no wheel (i.e.,
there is one vertex with no output leg). Let us show that
$H^*(F_i/F_{i-1},\partial)=0$ for all $i < n$.  Since $F_n/F_{n-1}$ is
spanned by graphs in degree $-n$, this will, by an elementary spectral
sequence argument, imply the statement.

Indeed, $F_i/F_{i-1}$ is spanned by graphs with exactly $i$ edges in
the wheel. The differential induced by $\partial$ acts by contractions
of edges outside the wheel. So, given a partition of the set
$\{1,\dots,n\}$ into $i$ subsets, $I_1\sqcup\dots \sqcup
I_i=\{1,\dots,n\}$, we can consider a subcomplex in $F_i/F_{i-1}$
spanned by graphs such that, if we cut all edges in the wheel, we
obtain $i$ graphs whose legs are marked by $I_1,\dots,I_i$.
Obviously, $F_i/F_{i-1}$ splits into the direct sum of such
subcomplexes.

The graphs with legs marked by $I_j$ form a complex
$\mathcal{C}(I_j)$. The subcomplex of $F_i/F_{i-1}$ corresponding to
the partition $I_1\sqcup\dots \sqcup I_i=\{1,\dots,n\}$ is the direct
sum of $(i-1)!$ copies of $\bigotimes_{j=1}^i\mathcal{C}(I_j)$.

So, if we prove that at least one complex $\mathcal{C}(I_j)$ for
$j=1,\dots,i$ is acyclic, we immediately conclude that $F_i/F_{i-1}$
is acyclic. We will see that $\mathcal{C}(I_j)$ is acyclic if
$|I_j|\geq 2$. This explains why $F_i/F_{i-1}$ is acyclic for $i<n$
(there must be at least one set $I_j$ with at least 2 elements; we
assume $n\geq 2$) but it is not acyclic for $i=n$.

Consider $\mathcal{C}(I_j)$ assuming $|I_j|\geq 2$. It consists of
trees with $1$ output leg and $|I_j|+1$ input legs. We use the
convention that the last input leg is the wheel input,
attached to the same vertex (the special vertex) as the output leg
(to handle $F_0$ in the same way, we would need to add a virtual loop at the
vertex with no output leg). We introduce a two-term filtration
$\mathcal{C}'(I_j)\subset \mathcal{C}(I_j)$ such that
$\mathcal{C}'(I_j)$ consists of graphs whose the special vertex has
biarity $(1,k)$, $k>2$. It is easy to see that $\mathcal{C}'(I_j)$ and
$\mathcal{C}(I_j)/\mathcal{C}'(I_j)$ are canonically isomorphic.
Indeed, slightly abusing notations, we can say that
$\mathcal{C}(I_j)/\mathcal{C}'(I_j)$ is generated by graphs with the
special vertex of biarity (1,2). The isomorphism is then given by the
contraction of the unique edge attached to the special vertex of
a~graph in $\mathcal{C}(I_j)/\mathcal{C}'(I_j)$. The first level of the
associated spectral sequence has two identical rows and the
differential is an isomorphism of these rows.  So, the cohomology vanishes
at the next level of this spectral sequence.  Therefore $F_i/F_{i-1}$
is acyclic for $i<n$ and  Theorem~\ref{thm-Com-wheeled-Koszul} is
proved.%
\qed

Now we can proceed with the proof of
Theorem~B.  Wheeled Koszulness implies that
the minimal resolution of $\mathsf{Com}^\circlearrowright$ is
generated by the standard $(1,n)$-generators of $\mathsf{Com}_\infty$
and by some new $(0,n)$-generators which form a basis of $H^{-n} \hskip
-.2em\left(B^\wc(\mathsf{Com}^\circlearrowright)(0,n)\right)$,
$n\geq 2$.

Our goal is to identify a basis of
$H^{-n}\hskip -.2em\left(B^\wc(\mathsf{Com}^\circlearrowright)(0,n)\right)$,
$n\geq 2$, with the generators introduced
in Theorem~B(ii), and to show that the differential
of the wheeled cobar construction on $(\mathsf{Com}^\circlearrowright)^!$
coincides with~\eqref{wheeled-delta}. To do so, we
study the wheeled double bar complex of
$\mathsf{Com}^\circlearrowright$ in the same way as we analyzed the
ordinary double bar complex of $\mathsf{Com}$ in
Subsection~\ref{ordinary-double-bar-Com}.

We define the bigrading of
$\Omega^\circlearrowright(B^\circlearrowright(\mathsf{Com}^\circlearrowright))(0,n)$
by $\deg_1(G)=-e(G)$ and $\deg_2(G)= e^\circ(G)$, where $e(G)$ denotes
the number of all edges of the underlying graph $G$ and $e^\circ(G)$
the number of white edges of $G$ (it is the same bigrading as in
the proof of Theorem~\ref{ll}).  The nontrivial part of this bicomplex
lies in
the triangle bounded by $\deg_1=-n$, $\deg_1=-\deg_2$, and $\deg_2=0$
(see example for $n=2$ in Subsection~\ref{wheeled-db-com}).

Consider the $\deg_2=0$ part of
$\Omega^\circlearrowright(B^\circlearrowright(\mathsf{Com}^\circlearrowright))(0,n)$
spanned by graphs with no white edges. Denote this subcomplex by
$\mathcal{D}_n$.  It follows from our constructions that
$\left(\mathcal{D}_n,\partial_1\right)$ is isomorphic to
$\left(B^\wc(\mathsf{Com}^\circlearrowright)(0,n),\partial\right)$.

To identify the generators of $H^{-n}(\mathcal{D}_n)$, we return to the
argument with filtration $F_*$ from the proof of
Theorem~\ref{thm-Com-wheeled-Koszul} (slightly abusing notations, we
consider it on $(\mathcal{D}_n,\partial_1)$). Note that $\partial_1$
restricted to $F_n/F_{n-1}$ is equal to zero. So, there is a basis of
$H^{-n}(\mathcal{D}_n)$, $n\geq 2$, consisting of $(n-1)!$ elements (we
need to fix the cyclic order of legs markings), which are represented
by wheels with no edges outside it, plus some correction terms arising
from the spectral sequence associated to $F_*$. As in the case
of the ordinary double bar complex, it is possible to give an explicit
description of these correction terms.

Let us fix a cyclic order of the set $\{1,\dots,n\}$. Consider the sum
of all planar graphs of genus~$1$ with trivalent vertices such that:
(i) all edges are black, (ii) the wheel contains at least two edges,
(iii) the wheel is oriented in the clockwise direction and (iv) when
we go around the graph in the counterclockwise direction, we meet the
input legs in the prescribed cyclic order.  There is, modulo an
overall sign, a unique choice of signs of these graphs such that
their sum is a $\partial_1$-closed element.

This $\partial_1$-closed element is exactly
what we denoted
in~Theorem~B(ii) by
$$
\begin{xy}
<0mm,0mm>*{\bullet},
<-5mm,-5mm>*{}**@{-},
<-2mm,-5mm>*{}**@{-},
<2mm,-5mm>*{}**@{-},
<5mm,-5mm>*{}**@{-},
<0mm,-7mm>*{_{\sigma(1)\sigma(2)\cdots\sigma(n)}},
\end{xy}
$$
where $\sigma\in\bS_n$ determines the cyclic order
(we identify  elements that differ by a
cyclic permutation).
The first few examples are (we omit the signs):
\begin{align*}
\begin{xy}
<0mm,0mm>*{\bullet},
<-3mm,-5mm>*{}**@{-},
<3mm,-5mm>*{}**@{-},
<-3mm,-7mm>*{_{\sigma(1)}},
<3mm,-7mm>*{_{\sigma(2)}},
\end{xy}
& =
\begin{xy}
<0mm,0mm>*{\bullet},
<-3mm,-5mm>*{}**@{-},
<3mm,-5mm>*{}**@{-},
<3mm,-5mm>*{\bullet};
<0mm,-10mm>*{}**@{-},
<0mm,-12mm>*{_{\sigma(2)}},
<-3mm,-7mm>*{_{\sigma(1)}},
(0,0)*{} \ar@{-}@(u,r) (3,-5)*{}
\end{xy}
& \forall\sigma\in\bS_2; \\
\begin{xy}
<0mm,0mm>*{\bullet},
<-3mm,-5mm>*{}**@{-},
<0mm,-5mm>*{}**@{-},
<3mm,-5mm>*{}**@{-},
<0mm,-7mm>*{_{\sigma(1)\sigma(2)\sigma(3)}},
\end{xy}
& =
\begin{xy}
<0mm,0mm>*{\bullet},
<-3mm,-5mm>*{}**@{-},
<3mm,-5mm>*{}**@{-},
<3mm,-5mm>*{\bullet};
<6mm,-10mm>*{}**@{-},
<0mm,-10mm>*{}**@{-},
<6mm,-10mm>*{\bullet};
<3mm,-15mm>*{}**@{-},
<3mm,-17mm>*{_{\sigma(3)}},
<0mm,-12mm>*{_{\sigma(2)}},
<-3mm,-7mm>*{_{\sigma(1)}},
(0,0)*{} \ar@{-}@(u,r) (6,-10)*{}
\end{xy}
+
\begin{xy}
<0mm,0mm>*{\bullet},
<-3mm,-5mm>*{}**@{-},
<3mm,-5mm>*{}**@{-},
<3mm,-5mm>*{\bullet};
<0mm,-10mm>*{}**@{-},
<0mm,-10mm>*{\bullet};
<-3mm,-15mm>*{}**@{-},
<3mm,-15mm>*{}**@{-},
<3mm,-17mm>*{_{\sigma(3)}},
<-3mm,-17mm>*{_{\sigma(2)}},
<-3mm,-7mm>*{_{\sigma(1)}},
(0,0)*{} \ar@{-}@(u,r) (3,-5)*{}
\end{xy}
+
\begin{xy}
<0mm,0mm>*{\bullet},
<-3mm,-5mm>*{}**@{-},
<3mm,-5mm>*{}**@{-},
<3mm,-5mm>*{\bullet};
<0mm,-10mm>*{}**@{-},
<0mm,-10mm>*{\bullet};
<-3mm,-15mm>*{}**@{-},
<3mm,-15mm>*{}**@{-},
<3mm,-17mm>*{_{\sigma(1)}},
<-3mm,-17mm>*{_{\sigma(3)}},
<-3mm,-7mm>*{_{\sigma(2)}},
(0,0)*{} \ar@{-}@(u,r) (3,-5)*{}
\end{xy}
+
\begin{xy}
<0mm,0mm>*{\bullet},
<-3mm,-5mm>*{}**@{-},
<3mm,-5mm>*{}**@{-},
<3mm,-5mm>*{\bullet};
<0mm,-10mm>*{}**@{-},
<0mm,-10mm>*{\bullet};
<-3mm,-15mm>*{}**@{-},
<3mm,-15mm>*{}**@{-},
<3mm,-17mm>*{_{\sigma(2)}},
<-3mm,-17mm>*{_{\sigma(1)}},
<-3mm,-7mm>*{_{\sigma(3)}},
(0,0)*{} \ar@{-}@(u,r) (3,-5)*{}
\end{xy}
& \forall\sigma\in\bS_4; \\
\end{align*}

Up to signs,
the formula for the differential of these new generators is almost
obvious . The only subtlety is that, if we want
$\partial_2$ to induce $\partial$ as in
equation~\eqref{wheeled-delta}, we must consider also binary graphs of
genus $1$ satisfying conditions (i), (iii), and (iv), but with only one
edge in the wheel.  However, it is easy to check that we need to add each
such a graph twice and with opposite signs. So, they do not affect the
result.

Theorem~B is therefore proved modulo signs.
We leave the sign issue to the reader.%
\qed

\vskip 1cm

\noindent
{\bf Acknowledgment.}
The first two authors would like to express their thanks to I.H.E.S.\
for a very stimulating and pleasant atmosphere, and acknowledge
the support of the European
Commission through its 6th Framework Programme ``Structuring the
European Research Area'' and the contract Nr.~RITA-CT-2004-505493 for
the provision of Transnational Access implemented as Specific Support.

\def\cprime{$'$}

\end{document}